\documentclass[10pt,reqno]{amsart}

\usepackage[utf8]{inputenc}
\usepackage[T1]{fontenc}
\usepackage{mlmodern}       % Modern Computer Modern font
\usepackage{bbm}            % Indicator function \mathbbm{1}
 \usepackage{tikz-cd}
\DeclareFontFamily{U}{wncy}{}
\DeclareFontShape{U}{wncy}{m}{n}{<->wncyr10}{}
\DeclareSymbolFont{mcy}{U}{wncy}{m}{n}
\DeclareMathSymbol{\Sh}{\mathord}{mcy}{"58}

\usepackage{amsmath}
\usepackage{mathtools}     % Adds support for \eqqcolon, \coloneq, etc.
\usepackage{amsfonts}
\usepackage{amssymb}
\usepackage{amsthm}
\usepackage{mathrsfs}      % Support for \mathscr

\numberwithin{equation}{section}

\DeclareFontFamily{U}{mathx}{} % Wide accents
\DeclareFontShape{U}{mathx}{m}{n}{<-> mathx10}{}
\DeclareSymbolFont{mathx}{U}{mathx}{m}{n}
\DeclareMathAccent{\widehat}{0}{mathx}{"70}
\DeclareMathAccent{\widecheck}{0}{mathx}{"71}

\usepackage{setspace}
\usepackage{enumitem}
\usepackage{float}
\usepackage{makecell}
\usepackage{tablefootnote}
\usepackage{caption}
\usepackage[labelformat=simple]{subcaption}
\usepackage{csquotes}      % Support for \begin{displayquote}

\restylefloat{table}
\newlist{enumlemma}{enumerate}{3}
\setlist[enumlemma]{label*={ (\alph*)}, ref={(\alph*)}}
\newlist{enumcount}{enumerate}{3}
\setlist[enumcount]{label*={ (\arabic*)}, ref={(\arabic*)}}

\usepackage{marginnote}
\usepackage[colorinlistoftodos]{todonotes}

\usepackage[most]{tcolorbox}
\tcbuselibrary{theorems, skins, breakable}

\tcbset{
    sharp/.style={
        sharp corners,
        breakable
    }
}

\newtcolorbox{sumbox}[1][]{
    colback=blue!5!white,
    colframe=blue!75!black,
    fonttitle=\bfseries,
    title={Summary: #1},
    breakable,
    sharp corners
}
\newtcolorbox{chkbox}[1][]{
    colback=yellow!10!white,
    colframe=yellow!50!black,
    fonttitle=\bfseries,
    title={Check: #1},
    breakable,
    sharp corners
}

\usepackage[colorlinks=true, linkcolor=blue, citecolor=blue, linktocpage=true, urlcolor=black]{hyperref}

\usepackage[alphabetic, backrefs, abbrev,nobysame,lite]{amsrefs}

\usepackage[capitalize, sort]{cleveref}

\crefformat{section}{\S#2#1#3}

\usepackage[toc, title]{appendix}
\theoremstyle{plain}
\newtheorem{theorem}{Theorem}[section]

\newtheorem*{CORO*}{Corollary}

\newtheorem{proposition}[theorem]{Proposition}
\newtheorem{lemma}[theorem]{Lemma}
\newtheorem{corollary}[theorem]{Corollary}

\newtheorem{claim}[theorem]{Claim}
\crefname{claim}{Claim}{Claims}

\theoremstyle{definition}
\newtheorem{definition}[theorem]{Definition}
\newtheorem{remark}[theorem]{Remark}

\newtheorem{question}[theorem]{Question}

\newtheorem{setting}[theorem]{Settings}
\crefname{setting}{Setting}{Settings}

\theoremstyle{remark}
\newtheorem*{remark*}{Remark}

\newtcbtheorem[use counter*=theorem, number within=section, crefname={exercise}{Exercises}]{kextcb}{Key exercise}{
    breakable,
    sharp corners,
    colback=red!1,
    colframe=red!35!black,
    fonttitle=\bfseries
}{ex}

\newtcbtheorem[use counter*=theorem, number within=section, crefname={exercise}{Exercises}]{extcb}{Exercise}{
    breakable,
    sharp corners,
    colback=yellow!1,
    colframe=yellow!35!red,
    fonttitle=\bfseries
}{ex}

\newtcbtheorem[use counter*=theorem, number within=section, crefname={example}{examples}]{egtcb}{Example}{
    breakable,
    sharp,
    colback=yellow!1,
    colframe=yellow!35!black,
    fonttitle=\bfseries
}{eg}

\newtcbtheorem[use counter*=theorem, number within=section, crefname={proposition}{propositions}]{proptcb}{Proposition}{
    breakable,
    sharp,
    colback=green!1,
    colframe=blue!20!black,
    fonttitle=\bfseries
}{prop}

\newtcbtheorem[use counter*=theorem, number within=section, crefname={proposition}{propositions}]{kproptcb}{Proposition}{
    breakable,
    sharp,
    colback=red!1,
    colframe=red!35!black,
    fonttitle=\bfseries
}{prop}

\newtcbtheorem[use counter*=theorem, number within=section, crefname={lemma}{lemmas}]{lemtcb}{Lemma}{
    breakable,
    sharp,
    colback=green!1,
    colframe=green!20!black,
    fonttitle=\bfseries
}{lem}

\newtcbtheorem[use counter*=theorem, number within=section, crefname={notations}{Notations}]{notatcb}{Notations}{
    breakable,
    sharp,
    colback=yellow!1,
    colframe=yellow!35!black,
    fonttitle=\bfseries
}{nota}

\newtcbtheorem[use counter*=theorem, number within=section, crefname={Theorem}{Theorems}]{theoremtcb}{Theorem}{
    breakable,
    sharp,
    colback=red!5,
    colframe=red!35!black,
    fonttitle=\bfseries
}{thm}

\newtcbtheorem[use counter*=theorem, number within=section, crefname={Corollary}{Corollaries}]{corotcb}{Corollary}{
    breakable,
    sharp,
    colback=red!5,
    colframe=red!35!black,
    fonttitle=\bfseries
}{coro}

\newtcbtheorem[use counter*=theorem, number within=section, crefname={Definition}{Definitions}]{defnlemtcb}{Definition/Lemma}{
    breakable,
    sharp,
    colback=red!5,
    colframe=red!35!black,
    fonttitle=\bfseries
}{defn}

\newtcbtheorem[use counter*=theorem, number within=section, crefname={definition}{Definitions}]{defntcb}{Definition}{
    breakable,
    sharp,
    colback=yellow!1,
    colframe=yellow!10!blue,
    fonttitle=\bfseries
}{defn}

\newtcbtheorem[use counter*=theorem, number within=section, crefname={exercise}{Exercises}]{keyrmktcb}{Key Remark}{
    breakable,
    sharp corners,
    colback=red!1,
    colframe=red!35!black,
    fonttitle=\bfseries
}{rmk}

\long\def\symbolfootnote[#1]#2{\begingroup\def\thefootnote{\fnsymbol{footnote}}\footnote[#1]{#2}\endgroup}

\DeclareMathOperator{\Aut}{Aut}

\DeclareMathOperator{\supp}{supp}

\DeclareMathOperator{\End}{End}
\DeclareMathOperator{\Hom}{Hom}

\DeclareMathOperator{\id}{Id}

\def\Im{\mathrm{Im}}
\def\im{\Im}
\def\Id{\id}

\providecommand{\restrict}{}
\renewcommand{\restrict}[1]{\raisebox{-.0ex}{$\upharpoonright$}_{#1}}

\newcommand{\op}{{\rm op}}

\newcommand{\Lie}{{\rm Lie}}

\newcommand{\tr}{{\rm tr}}

\newcommand{\Vol}{{\rm Vol}}

\renewcommand{\emph}[1]{{\bf #1}} %%%% <- \emph becomes \textbf
\renewcommand{\hat}{\widehat}
\renewcommand{\phi}{\varphi}

\renewcommand{\bar}{\overline}

\newcommand{\diff}{\mathrm{Diff}}
\newcommand{\Diff}{\diff}

\newcommand{\Aff}{\mathrm{Aff}}

\newcommand{\Sl}{\mathrm{SL}}
\def\SL{\Sl}
\newcommand{\Gl}{\mathrm{GL}}
\def\GL{\Gl}
\newcommand{\SO}{\mathrm{SO}}

\providecommand{\C}{}\renewcommand{\C}{\mathbb{C}}

\providecommand{\Rbb}{}\renewcommand{\Rbb}{\mathbb{R}}

\providecommand{\Hbb}{}\renewcommand{\Hbb}{\mathbb{H}}

\providecommand{\Cbb}{}\renewcommand{\Cbb}{\mathbb{C}}

\providecommand{\Zbb}{}\renewcommand{\Zbb}{\mathbb{Z}}

\newcommand{\Ical}{\mathcal I}

\newcommand{\Bcal}{\mathcal{B}}
\newcommand{\Ocal}{\mathcal{O}}

\newcommand{\Dcal}{\mathcal{D}}

\newcommand{\Ecal}{\mathcal{E}}
\newcommand{\Acal}{\mathcal{A}}
\newcommand{\Ccal}{\mathcal{C}}

\newcommand{\Tcal}{\mathcal{T}}
\newcommand{\Lcal}{\mathcal{L}}

\newcommand{\Fsc}{\mathscr{F}}

\newcommand{\klie}{\mathfrak{k}}
\newcommand{\hlie}{\mathfrak{h}}

\newcommand{\glie}{\mathfrak{g}}
\newcommand{\clie}{\mathfrak{c}}
\newcommand{\slie}{\mathfrak{s}}

\newcommand{\mlie}{\mathfrak{m}}
\newcommand{\plie}{\mathfrak{p}}

\newcommand{\nlie}{\mathfrak{n}}

\newcommand{\alie}{\mathfrak{a}}

\newcommand{\wtd}{\widetilde}
\newcommand{\what}{\widehat}

\newcommand{\Fr}{\textrm{Fr}}
\newcommand{\ad}{\textrm{ad}}

\newcommand{\Hbf}{\mathbf{H}}

\newcommand{\alge}{\textrm{alg}}

\title[Volume preserving Anosov $\Gamma$ actions]{On the classification of volume-preserving Anosov actions by higher-rank lattices}

\author[H.~Lee]{Homin Lee}
\address{HCMC, Korea Institute for Advanced Study, Seoul, South Korea}
\email{hominlee@kias.re.kr}

\author[A.~Maldonado]{Amadeus Maldonado}
\address{Northwestern University, Evanston, United States}
\email{amadeusmaldonado2027@u.northwestern.edu}

\begin{document}
%\symbolfootnote[0]{\itshape  Date: \today.}
\maketitle

\let\oldtocsection=\tocsection
\let\oldtocsubsection=\tocsubsection
\let\oldtocsubsubsection=\tocsubsubsection
\begin{abstract} Let $\Gamma$ be a lattice in the higher rank simple Lie group $G$ and let $M$ be a closed smooth manifold. Let $\alpha:\Gamma\to \textrm{Diff}_{\textrm{Vol}}^{2+\textrm{H\"older}}(M)$ be a $C^{2+\textrm{H\"older}}$-volume-preserving Anosov $\Gamma$-action on $M$, that is, there is a $\gamma_0\in\Gamma$ such that $\alpha(\gamma_0)$ is an Anosov diffeomorphism. We showed that $M$ is $C^{2+\textrm{H\"older}}$-diffeomorphic to an infra-nilmanifold and $\alpha$ is $C^{2+\textrm{H\"older}}$-conjugate to an affine action. If the action $\alpha$ is $C^\infty$, then $\alpha$ is $C^\infty$ conjugate to an affine action. \end{abstract}
\renewcommand{\tocsection}[2]{\hspace{0em}\oldtocsection{#1}{#2}}
\renewcommand{\tocsubsection}[2]{\hspace{1em}\oldtocsubsection{#1}{#2}}
\renewcommand{\tocsubsubsection}[2]{\hspace{2em}\oldtocsubsubsection{#1}{#2}}
\setcounter{tocdepth}{1}

\tableofcontents

\section{Introduction}

Let $M$ be a closed connected smooth manifold. A $C^1$ diffeomorphism $f\colon M\to M$ is \emph{Anosov} if there are a continuous $Df$-invariant splitting $TM=E^s\oplus E^u$,  a constant $C\geq1$, and a constant $0<\lambda<1$ such that, for all $x\in M$, $j\geq1$, $v^s\in E^s_x$, and $v^u\in E^u_x$,
\[ \|D_x f^j v^s\|\leq C\lambda^j\|v^s\|,\qquad\|D_x f^{-j} v^u\|\leq C\lambda^j\|v^u\|.\]
Compactness of $M$ makes this definition independent of the choice of Riemannian metric. An action
$\alpha\colon\Gamma\to\Diff^1(M)$ is called an \emph{Anosov action} if $\alpha(\gamma_0)$ is an Anosov diffeomorphism for some $\gamma_0\in\Gamma$.

The standard action of $\mathrm{SL}(n,\mathbb Z)$ on $\mathbb T^n=\mathbb R^n/\mathbb Z^n$
by toral automorphisms is the basic example. Every matrix in the group preserves Lebesgue
volume, and every matrix with no eigenvalues on the unit circle induces an Anosov
diffeomorphism. We ask to what extent these two properties determine the manifold and the
action. Fisher formulated this question as a conjecture for lattices in
$\mathrm{SL}(n,\mathbb R)$ \cite{Fisher11}. Brown, Rodriguez Hertz, and Wang restated it in
\cite{BRHW17}.  \cref{thm:main} gives an
affirmative answer for these conjectures. The following corollary of our main theorem states the case $\Gamma=\mathrm{SL}(n,\mathbb
Z)$.

\begin{CORO*}[$\SL(n,\Zbb)$-case of \cref{thm:main}]
Let $n\geq3$ and let $M$ be a closed connected smooth manifold. Suppose that $\alpha\colon \SL(n,\Zbb)\to\Diff_{\Vol}^{\infty}(M)$ is a smooth volume-preserving Anosov action. Then there is a $C^\infty$ diffeomorphism from $M$ to an infranilmanifold that conjugates $\alpha$ to an affine action.
\end{CORO*}
One of the main novelty of our theorem, \cref{thm:main}, is that we do not have any assumptions on the dimension or topology of the manifold. In particular, the dynamics determine the infra-nilmanifold structure on $M$.

Here an \emph{infranilmanifold} is a compact quotient $\Pi\backslash N$, where $N$ is a
connected, simply connected nilpotent Lie group and $\Pi$ is a discrete group of affine
transformations of $N$. The group $\Pi$ acts freely and properly discontinuously, and its
image under the linear-part homomorphism to $\operatorname{Aut}(N)$ is finite. An affine
transformation of $N$ has the form $x\mapsto aA(x)$, with $a\in N$ and
$A\in\operatorname{Aut}(N)$. If it normalizes $\Pi$, it induces an affine diffeomorphism of
$\Pi\backslash N$.

Margulis' superrigidity theorem relates finite-dimensional representations of a higher-rank
lattice to representations of its ambient Lie group \cite{Mar91}. Zimmer's cocycle
superrigidity applies this principle to the derivative of a volume-preserving action. The
derivative maps the tangent space at $x$ to the tangent space at $\alpha(\gamma)x$. After a
measurable change of frame and passage to a covering group when necessary, the derivative
cocycle is given by a representation of the ambient group together with a commuting
compact-valued cocycle \cite{Zimmer1980,Zimmerbook,FisherMargulis}. An Anosov element
supplies additional information at every point of $M$. It has continuous stable and unstable
bundles, with uniform contraction on the stable bundle and uniform expansion on the unstable
bundle. The problem is to use this information to recover a geometric structure from the
measurable description of the derivative.

Hurder formulated a related rigidity conjecture, as did Katok and Lewis
\cite{Hurder1994,KatokLewis}. 
Labourie also formulated an Anosov rigidity conjecture and proposed
using the stable and unstable bundles to recover an invariant geometric structure
\cite{LabourieICM}. A longstanding topological question asks whether every closed manifold
admitting an Anosov diffeomorphism is homeomorphic to an infranilmanifold
\cite{Franks}. 

As a companion with higher rank lattice actions, the Katok--Spatzier conjecture grew out of the rigidity theory of higher-rank abelian group $\Rbb^k\times \Zbb^l$-actions, $k+l\ge 2$. Its classical formulation predicted that a smooth transitive Anosov action of a higher-rank abelian group on a compact manifold, with no rank-one factor, would become conjugate to an algebraic action on a finite cover \cite[Conjecture~1.7]{SV}.  
\footnote{It is worth to mention that $\Rbb^k$ action case, Vinhage constructed smooth Anosov $\mathbb R^2$-actions without smooth rank-one factors for which this conclusion fails \cite{KurtCounterExample}. The example is given by time change, hence, for $\Rbb^k$-action case, the conjecture can be replaced to up-to-time change the action is algebraic.}

A central themes in higher-rank rigidity conjectures are two-fold: 
\begin{enumerate}
\item Construct the underlying homogeneous structure on the  manifold from the given (hyperbolic) dynamics.
\item Show that the given (non-linear) action on a (bi-)homogeneous space is conjugate to an algebraic action, e.g., affine action, on the (bi-)homogeneous space.
\end{enumerate}

Let us mention two important recent breakthroughs \cite{SV,DSVX} in this direction, especially, their works provide new ways to build a homogeneous structure from higher-rank group actions. In \cite{SV}, Spatzier and Vinhage classified totally Cartan $\Rbb^k$ actions. They developed a construction using paths along the one-dimensional coarse Lyapunov foliations. In their work, the relations among these paths determine the Lie group underlying the model. In \cite{DSVX} Damjanovi\'c, Spatzier, Vinhage, and Xu built on this idea for higher-dimensional coarse Lyapunov foliations. They showed how the groups acting along individual leaves fit together by proving that the action of their free product factors through a finite-dimensional Lie group. For smooth actions of semisimple Lie groups whose simple factors all have real rank at least two, this yields a bi-homogeneous model on a finite cover under their hypotheses of total partial hyperbolicity, accessibility, and volume preservation by a maximal real split Cartan subgroup. Earlier than those works, for certain totally Cartan $\mathbb R^k$-actions, Kalinin and Spatzier showed that vector fields tangent to the coarse Lyapunov foliations, together with the generators of the action, span a Lie algebra that integrates to a transitive group action on a finite cover \cite{CartanKS}.

The main theorem in this paper provides an infra-nilmanifold structure on manifold from given Anosov actions by a lattice $\Gamma$ in connected higher-rank simple Lie group $G$. Here, \emph{higher-rank} means that the real rank of the Lie algebra $\Lie(G)$ is at least two.

\begin{theorem}\label{thm:main}
Let $G$ be a connected higher-rank simple Lie group, let $\Gamma<G$ be a lattice, and let $M$ be a closed connected smooth manifold. Fix $\eta>0$, and suppose that the action $\alpha\colon\Gamma\to\Diff^{2+\eta}(M)$ preserves a smooth positive density and $\alpha(\gamma_0)$ is an Anosov diffeomorphism for some $\gamma_0 \in \Gamma$. Then, for some $\eta'>0$, there are an infranilmanifold $X$, an affine action $\rho$ of $\Gamma$ on $X$, and a $C^{2+\eta'}$ diffeomorphism $h\colon M\to X$ such that
\[ h\circ\alpha(\gamma)=\rho(\gamma)\circ h\qquad\text{for every }\gamma\in\Gamma.\]
If $\alpha$ is $C^\infty$, then $h$ can be chosen to be $C^\infty$.
\end{theorem}
The $C^{2+\eta}$ hypothesis enters in the construction of an invariant connection. The transformation law for a connection involves $D^2\alpha(\gamma)$, which is H\"older continuous under this hypothesis. The resulting developing map yields a $C^{2+\eta'}$ conjugacy on $M$. For the last $C^\infty$ conclusion, the same argument with the proof of the smooth rigidity theorem in  \cite{BRHW17,RHW2014,FKS2013} can be directly applied.

The volume hypothesis supplies the invariant measure needed to apply cocycle superrigidity. This raises the following question.

\begin{question}
Let $G$, $\Gamma$, and $M$ be as in Theorem~\ref{thm:main}. If a smooth action
of $\Gamma$ on $M$ contains an Anosov diffeomorphism, must the
whole action preserve a smooth positive density on $M$?
\end{question}

\subsection{Motivation and history}
Motivated in part by Zimmer's cocycle superrigidity theorem, the Zimmer program studies how the algebraic structure of higher-rank Lie groups and their lattices constrains their smooth actions on compact manifolds.

Among many other directions, the following three directions have seen substantial progress and are particularly relevant here.
\begin{enumerate}
    \item Actions on low-dimensional manifolds.
    \item Actions preserving geometric structures.
    \item Actions with uniform and non-uniform hyperbolic dynamics.
\end{enumerate}

The first asks in which dimensions every smooth action of $\Gamma$ has finite image. Zimmer
initiated this study for actions preserving a volume or a geometric structure
\cite{Zimmer1984,Zimmer1986}. Brown, Fisher, and Hurtado established dimension bounds for
actions of cocompact lattices, of $\mathrm{SL}(n,\mathbb Z)$, and of nonuniform lattices
\cite{BFHUniform,BFHSLmZ,BFHNonuniform}. Brown, Damjanovi\'c, and Zhang treated $C^1$
actions of cocompact lattices, and An, Brown, and Zhang established further cases for
nonsplit groups \cite{BDZ,ABZ}. In their work on actions on the line, Deroin and Hurtado
proved that an irreducible lattice in a higher-rank semisimple Lie group with finite center
admits no left-invariant total order \cite{DeroinHurtado}. In the case of minimal dimension, for instance, studying actions on $n$-dimensional manifolds by $\SL(n,\Rbb)$ or lattices in $\SL(n,\Rbb)$, there are recent works including \cite{FisherMelnick2024,BrownLee}. See also, \cite{Son2025,ArosemenaSerrato2025}.

A second direction asks what can be deduced from a geometric structure preserved by the
action. Zimmer and Goetze studied actions preserving connections, and Feres and Zeghib
obtained rigidity results for affine actions in low dimensions
\cite{Zimmer1986,Goetze1994,Feres1993,Zeghib1997,ZimmerAlgebraicHull}. Related work on rigid geometric
structures includes Melnick's and Pecastaing's results \cite{MelnickSurvey,PecastaingConformal,PecastaingProjective}. Particularly close to our question is the
theorem of Feres and Labourie for a lattice in $\mathrm{SL}(n,\mathbb R)$ acting on an
$n$-manifold. Under hyperbolicity and transitivity assumptions on the suspended action, they
identify the manifold as a flat torus and the lattice action as affine \cite{FeresLabourie}.
There are many related works in this direction, for instance, \cite{BaderFrancesMelnick,MelnickFrobenius}.

The third direction starts with hyperbolicity. Early results on local and global rigidity of
toral and affine lattice actions include work of Katok and Lewis, Hurder,
Katok--Lewis--Zimmer, Qian, and Fisher and Margulis
\cite{KatokLewis1991,KatokLewis,Hurder1992,Hurder1994,KLZ,Qian,FisherMargulisLocal}. Goetze
and Spatzier obtained a global classification for Cartan lattice actions under stronger
hyperbolicity hypotheses \cite{GoetzeSpatzier1997,GoetzeSpatzier1999}. For volume-preserving
Anosov lattice actions on a given nilmanifold, Margulis and Qian proved topological rigidity
assuming the action has a finite orbit \cite{MargulisQian}. Brown, Rodriguez Hertz, and Wang
proved topological and smooth global rigidity on nilmanifolds when a finite-index subgroup
acts on the universal cover by lifts of the given action. They also gave criteria ensuring the existence of
such a lift \cite{BRHW17}. See also \cite{LeeSandfeldt2024}. The first author treated weakly irreducible lattices in a broader class of
semisimple groups under an irreducibility hypothesis on a fully supported invariant measure
\cite{LeeLattice2024}.

It is worth to mention that, in the setting of our theorem, the hyperbolicity of the action supplies a geometric structure of the kind discussed above. We construct an affine connection preserved by the lattice.

In the present problem, the topology of $M$ is part of the conclusion. There are precedents for recovering the topology in the dimension of the defining
representation as well as recent works \cite{SV,DSVX} we already mentioned. For lattices in $\mathrm{SL}(n,\mathbb R)$ with $n\geq3$,
$\mathrm{Sp}(2n,\mathbb R)$ with $n\geq2$, and $\mathrm{SO}(m,m)$ with $m\geq5$, the first
author studied volume-preserving actions on manifolds of dimensions $n,2n,2m$, respectively
\cite{Lee}. Those results require only one element with a dominated splitting. In the
special linear case, the proof produces a codimension-one Anosov element and uses the
Franks--Newhouse classification. In the other two cases, the proof uses spectral pinching
and the Brin--Manning theorem to identify the topology. The affine-conjugacy statements in
those cases include a lifting hypothesis \cite{Franks,Newhouse,BrinManning,Lee}. Brown and
the first author subsequently removed the volume-preservation hypothesis for smooth Anosov
actions in the defining dimensions of $\mathrm{SL}(n,\mathbb R)$ for $n\geq4$, obtaining smooth affine conjugacy on an infratorus
\cite{BrownLee}. 

There are many algebraic example that does not admit codimension one Anosov diffeomorphims or even a quasi-conformal Anosov diffeomorphism in  arbitrary dimension. Hence, the way to recover the topology of manifold which is used previously cannot be extended to general dimension directly.

A higher-rank abelian subgroup can also supply rigidity for a lattice action. Katok and Spatzier used first-cohomology rigidity to study time
changes of standard higher-rank abelian Anosov actions \cite{KatokSpatzier}.
They later used abelian subactions to prove smooth local rigidity for
Anosov lattice actions by nilmanifold automorphisms \cite{KatokSpatzier1997}.
Recently, Rodriguez Hertz and Wang have studied related classification questions for
real-analytic lattice actions \cite{RHWvolume,RHWanalytic}.

The common aim of \cite{CartanKS,SV,DSVX} and the present paper is to recover the model space from the dynamics as we mentioned before.  
Our construction takes place on the universal cover $\wtd M$. Many earlier arguments use the suspension $(G\times M)/\Gamma$. For a nonuniform lattice $\Gamma$, this space is noncompact. Even for a cocompact lattice, deducing hyperbolicity of a flow on the suspension from the hyperbolicity of a single lattice element is a separate problem. See the discussion in \cite[Section~2.5.7]{DSVX}.

We directly identify the universal cover $\wtd M$ with a nilpotent Lie group $N$. In order to do that we construct a flat connection with parallel torsion. If $T$ denotes the torsion of the lifted connection and $X$ and $Y$ are parallel vector fields on $\wtd M$, then $[X,Y]=-T(X,Y)$. Thus the torsion determines the Lie bracket of $N$. Under
the resulting identification $\wtd M\simeq N$, the deck transformations act affinely with finite linear holonomy. The quotient of $N$ by this deck action is the required infranilmanifold.

\subsection{Outline of proof and standing assumptions}
On \cref{sec2}, we recall necessary facts we used in the proof. The proof of \cref{thm:main} occupies the entire remaining manuscript.  Since it is quite involved, we describe each step below.

Throughout the proof, we use two coverings, topological- and algebraic-universal covers of simple Lie group in order to apply Zimmer's cocycle superrigidity and utilize Zariski density of lattices in simple algebraic groups. 

Note that if $G_0$ is a connected real Lie group, $G$ its topological-universal cover and $\Gamma_0<G_0$ is a lattice, then the pre-image of $\Gamma_0$ in $G$ is still lattice. 
In the statement of \cref{thm:main}, hence, we can consider the topological universal cover of $G$ and pre-image of $\Gamma$.

\begin{setting}[Standing assumptions]\label{setting} Throughout the paper, we assume the following setting:
\begin{enumerate}
\item 
Let $G$ be a connected, topologically simply connected real Lie group with simple Lie algebra $\glie$. Let $r=\textrm{rank}_\Rbb\glie$ and assume that $r=\textrm{rank}_{\Rbb}\glie\ge 2$. Let $\Gamma$ be a lattice in $G$.  
\item For $d>0$, let $M$ be a smooth closed $d$-dimensional connected orientable manifold and fix a smooth Riemannian metric on it. 
\item For $\eta>0$, let $\alpha:\Gamma\to \Diff^{2+\eta}(M)$ be a $C^{2+\eta}$ group action preserving a smooth density $\omega$ and such that $\alpha(\gamma_0)$ is an Anosov diffeomorphism for some $\gamma_0 \in \Gamma$.
We normalize $\omega$ so that $\int_M \omega = 1$ and denote by $\mu$ its associated probability measure.  
\end{enumerate}
 We call $\mu$ a \emph{volume measure}. Let $V=\Rbb^d$ be a fixed vector space called the \emph{model vector space}. 
\end{setting}
We will prove \cref{thm:main} under the \cref{setting} in \cref{sec3}--\cref{sec7}.

 \textbf{\cref{sec3}. Invariant H\"older compact group reduction $P$ in the frame bundle:} This step is highly motivated by the method in \cite{KLZ} and the usage of it in the first author's paper \cite{Lee}.

 However, instead of comparing only stable and unstable distributions, we consider the non-commutative unital algebra generated by stable projections. This algebra is large enough so that we can, by comparing measurable data from Zimmer's cocycle superrigidity and continuous data from an Anosov diffeomorphism and its conjugations, deduce that there is a continuous compact group $C$ reduction $P\subset \Fr(TM)$ which is compatible with the group action, that is, linearize the action up to compact error.  Here, $\Fr(TM)$ is the frame bundle of $M$.

 We can choose adapted local frame $b$, that is, for a small enough open set $\Ocal$, $b:\Ocal \to P$ gives a local section of $\Fr(TM)\to M$ on $\Ocal$. Let $b_1\in P$ and $b_2\in P$ be two local adapted frames with domains $\Ocal_1$ and $\Ocal_2$, respectively. Then, for $x\in\Ocal_1\cap \Ocal_2$, $(b_1)_x$ may not agree with $(b_2)_x$; they can be differ by an element in $C$. 
In terms of adapted local frame, the continuous $C$-reduction can be written as follows: 

Let $V$ be a fixed $\Rbb$-vector space with $\dim_\Rbb V= \dim M$. There is a continuous representation $\rho:G\to \SL(V)$, compact group $C\subset \SL(V)$ that commutes with $\rho(G)$, and a continuous $C$-reduction $P$ such that for any adapted frame, $b_x\in P_x$ and $b_{\alpha(\gamma)(x)}\in P_{\alpha(\gamma)(x)}$, there is a unique $k\in C$ such that \[ D_x\alpha(\gamma)b_x =b_{\alpha(\gamma)(x)}\rho(\gamma)k. \]

The continuous compact group reduction $P$, combining with the fact that superrigidity representation can not have zero weight (because of the Anosov diffeomorphism), already allows us to find many Anosov diffeomorphisms: $\alpha(\gamma)$ is an Anosov diffeomorphism as long as $\rho(\gamma)$ is hyperbolic where $\rho$ is the cocycle superrigidity homomorphism for the derivative cocycle.

 \textbf{\cref{sec4}.  Invariant continuous torsion free connection $\nabla$:} Fix a smooth connection $\nabla^0$. Using Zimmer's cocycle superrigidity again, we can find a measurable $\alpha(\Gamma)$-invariant torsion free connection $\nabla^{m}=\nabla^0+S^{m}$ on $M$ where $S^{m}$ is a measurable section of the H\"older bundle $\Ecal=\Hom(\textrm{Sym}^{2} TM, TM)$. Recall that the difference between two torsion free connections is a section of H\"older bundle $\Ecal = \Hom(\textrm{Sym}^2 TM , TM)$ over $M$. On the other hand, we can choose an Anosov diffeomorphism $\alpha(\gamma)$ for some $\gamma\in\Gamma$ so that $\Ecal$ can be decomposed into uniformly contracting, uniformly expanding, and isometric subbundles. Using expanding and contracting dynamics, we can improve the regularity of $S^{m}$ on uniformly expanding and contracting subbundles. On the neutral subbundle, as long as the cocycle is H\"older, we can upgrade measurable conformal structure to H\"older conformal structure, \cite{Butler}. Note that this is the place we need $C^{2+\textrm{H\"older}}$ regularity of the action $\alpha$ in order to get the H\"older cocycle on $\Ecal$ as the action on $\Ecal$ is given by second derivatives. Consequently, we can find a H\"older $\alpha(\Gamma)$-invariant (torsion free) connection $\nabla$ that agrees with $\nabla^{m}$ almost everywhere.

 \textbf{\Cref{sec5}. Regularity of the compact reduction.} Recall that we can find a higher rank free abelian subgroup $\Delta\simeq \Zbb^{\textrm{rank}_{\Rbb}G}$ (we call it a lattice Cartan, \cref{CartanGamma}) in $\Gamma$ such that $\alpha(\delta)$ is Anosov diffeomorphism as long as $\rho(\delta)$ is hyperbolic using the continuous compact group $C$-reduction $P$. Fix the lattice Cartan subgroup $\Delta$.

\emph{Local Lipschitz regularity.} For simplicity, let us assume that the superrigidity homomorphism $\rho$ is irreducible. Otherwise, one can consider each irreducible components. Choose a restricted weight $\xi$ that is primitive in the weight lattice, i.e., if $t\xi$ is in the weight lattice then $t\ge 1$.
Let $c$ be a short curve of length $\ell<1$ tangent to $E^\xi$. Through adapted frames at end points, the parallel transport can be written as the map $V\to V$.  Since $\nabla$ is $\Gamma$-invariant, we can compare parallel transport along $c$ with transport along $\alpha(\delta)(c)$ for $\delta\in\Delta$. The same comparison on the diagonal shows that the entire
block-diagonal part of $Q_c$ is a single element of $C$. Shrinking the curve $c$, we will see that 
\begin{enumerate}
\item an off-diagonal block part from $V_\lambda$ to $V_\nu$ with $\lambda\neq \nu$ will not be zero only if $\nu-\lambda=t\xi$ for some $t>0$, hence $t\ge1$. Indeed, it $\xi$ and $\nu-\lambda$ are not positively proportional, we can choose a Cartan element that contracts the curve $c$ and expands the block under conjugation. Transport along the shrinking curves stays bounded, so the block must vanish. 
\item a block-diagonal part, from the weight space $V_\lambda$ to $V_\lambda$, of the transport belongs to $C$, hence bounded size by similar argument as above.
\end{enumerate} 
The off-diagonal block has size $O(\ell^t)$, hence, $O(\ell)$.  Note that $\ell <1, t\ge 1\Rightarrow \ell^t \le \ell$. Indeed, with $\nu-\lambda=t\xi$, we expand the curve $c$ by $\alpha(\delta)$ until the length of the image curve $c_n\coloneqq\alpha(\delta^n)(c)$ is comparable to a fixed scale. If this iterate multiplies length along $E^\xi$ by $s$, it multiplies the norm of the block by $s^t$. At the fixed scale, transport is uniformly bounded. Since $s$ is comparable to $\ell^{-1}$, the original block is $O(\ell^t)$. This shows local Lipschitzness along $E^\xi$. 
 
 Finally, we can find finite elements $\gamma_1,\dots,\gamma_l$ in $\Gamma$ so that $E_i\coloneqq D\alpha(\gamma_i)E^\xi$ span $TM$. Furthermore, using the lemma in \cite{GHS}, nearby points $x,y$ can be joined by curves tangent to $E_i$ with total length $O(d(x,y))$. This proves that $P$ is locally Lipschitz.

\emph{From locally Lipschitz to $C^{1+\beta}$.}
Choose a locally Lipschitz adapted frame $b$. Its weak covariant derivative, taken modulo $\clie=\textrm{Lie}(C)$, has
essentially bounded coefficients
\[ a_b(x)(u) =\left[b_x^{-1}(\nabla_{b_xu}b)_x\right]       \in\mathfrak{gl}(V)/\clie,\qquad a_b(x) \in \Hom(V,\mathfrak{gl}(V)/\clie)\] for all $u\in V$.
In the quotient, a difference from the change of adapted frame is vanished. Therefore, the invariance of $\nabla$ now gives an equivariance equation for $a_b$ under $\Gamma$ action. Essential boundedness implies that almost every $a_b(x)$ has bounded $\Gamma$-orbit, therefore, $a_b(x)$ is in the $G$-fixed subspace. This $G$-invariance will provide that $a_b$ is a measurable $\Gamma$-invariant section of an isometric cocycle over an Anosov element. Again, using \cite{Butler}, we can get a H\"older representative of $a_b$. This says that weak derivative of $P$ has a H\"older representative which proves that $P$ is $C^{1+\beta}$ for some $\beta>0$.

 \textbf{\Cref{sec6}. A flat connection with parallel torsion.}
Recall that for the connection that makes left-invariant vector fields parallel, the torsion evaluated on two such fields is the negative of their Lie bracket. Hence, we want to find a connection on $M$ for which this identity recovers a single Lie algebra in every tangent space.

The coefficient $a_b$ obtained in \cref{sec5} is fixed by $G$, and its remaining dependence on the adapted frame is through $C$. Its $C$-orbit is a continuous $\Gamma$-invariant function of the
base point, hence is constant by ergodicity and full support. Choosing frames in which the coefficient equals one representative of this orbit reduces $C$ to its stabilizer. Repeating this step finitely many times gives a reduction $P_*$ for which the
coefficient is constant modulo the Lie algebra of its own structure group $C_*$. An equivariant lift of that coefficient defines a tensor $S$. The modified connection $\widetilde\nabla=\nabla-S$ preserves $P_*$ and an associated Riemannian metric. Its torsion is parallel and has a constant coefficient $T_0$ in adapted frames.

It remains to show that $-T_0$ is a nilpotent Lie bracket. The possible obstruction is curvature. Since $\widetilde\nabla$ is H\"older, its curvature is initially a distribution. The first Bianchi identity gives a linear equation for this curvature in terms of $T_0$. We prove that the corresponding linear map is injective on curvature tensors with values in
$\clie_*=\textrm{Lie}(C_*)$. This determines a regular, parallel curvature tensor with constant adapted coefficient $R_0$. At this point the two Bianchi identities give
the Jacobi identity on $\clie_*\oplus V$, where the bracket of $u,v\in V$ has $V$-component $-T_0(u,v)$ and $\clie_*$-component $-R_0(u,v)$. A generic element of the split Cartan algebra acts on this Lie algebra by a derivation that vanishes on $\clie_*$ and is invertible on $V$. The nilradical lemma then identifies $V$ with the largest nilpotent ideal. In particular $[V,V]\subset V$, so $R_0=0$. Thus $\widetilde\nabla$ is flat and $-T_0$ is the required nilpotent Lie bracket.

\medskip
\noindent
\textbf{\Cref{sec7}. The infranilmanifold and the affine conjugacy.}
Let $N$ be the simply connected nilpotent Lie group with Lie algebra $(V,-T_0)$. Flatness of $\wtd\nabla$ gives a global $C^{1+\beta}$ parallel frame on the universal cover $\widetilde M$. Its parallel fields satisfy the same bracket relations as the left-invariant fields on $N$. There is also a global reason that these fields can be integrated for all time. Their lengths are constant in the lifted parallel metric, and this metric is complete because $M$ is compact. The resulting developing map $H:\widetilde M\to N$ is therefore a global diffeomorphism.

Under $H$, every deck transformation becomes affine, with linear part in the compact group of bracket-preserving transformations coming from $C_*$. The deck group acts freely, discretely and cocompactly, so the quotient is an infranilmanifold. Moreover, every lift of each $\alpha(\gamma)$ preserves the connection and therefore sends parallel fields to parallel fields. In the coordinates on $N$ it is an affine transformation normalizing the deck group. It follows that the same developing map conjugates the full $\Gamma$-action on $M$ to an affine action on the quotient.

Finally, if $e_x:V\to T_x\widetilde M$ is the parallel frame and
$\ell_z:V\to T_zN$ is the smooth left-invariant frame on $N$, then
\[
                         dH_x=\ell_{H(x)}e_x^{-1}.
\]
Starting with $H\in C^1$, this identity first gives $H\in C^{1+\beta}$ and then $H\in C^{2+\beta}$. This will be enough to get the $C^{2+\textrm{H\"older}}$-conjugacy statement. After this, $C^{\infty}$ statement is followed directly from the argument in \cite{BRHW17}. This shows \cref{thm:main}.

\subsection*{Acknowledgments}
H.L. was supported in part by a KIAS Individual Grant (HP104101) via the June E Huh Center for Mathematical Challenges at Korea Institute
for Advanced Study.

A.M. was supported in part by grants from the NSF (DMS-2235451) and Simons Foundation (MPS-NITMB-00005320) to the NSF-Simons National Institute for Theory and Mathematics in Biology (NITMB).

This research was initiated during the first author's visit to Northwestern University. The authors appreciate the hospitality of Northwestern University. The authors thanks to Aaron Brown, Danijela Damjanovi\'c, David Fisher, Ralf Spatzier, and Kurt Vinhage for useful discussions.

\subsection*{AI Disclosure}
Parts of the proof were developed and refined through discussions with ChatGPT (ChatGPT-6 Astra). It was also used for routine literature searches. The manuscript itself was written and organized by the authors, with limited assistance from ChatGPT.

\subsection*{Note on this version.}
 The authors believe it is useful to make the present version of this manuscript available at this stage. The authors are continuing to improve the exposition and presentation, including correcting grammatical and typographical errors. An updated version with improved readability will be posted in the near future.

\section{Preliminary}\label{sec2}
\subsection{Lie groups and lattices}\label{sec2.0}  
We collect a few facts about Lie groups and lattices which one can find in many standard literatures.

For a real semisimple Lie algebra $\glie$, its real rank $\textrm{rank}_{\Rbb}(\glie)$ is the largest dimension of an abelian subalgebra $\alie\subset\glie$ such that every $\textrm{ad}(H)$, $H\in\alie$ is diagonalizable over $\Rbb$.

We assume the standing assumptions \cref{setting}. 
\subsubsection{Algebraic model} We also need its algebraic model:

\begin{claim}\label{Borel} Let $\Hbf$ be the algebraically simply connected simple real algebraic group with Lie algebra $\glie$ and let $H=\Hbf(\Rbb)^{\circ}$. 
Denote covering homomorphims by $q:G\to H$. 
\begin{enumerate}
\item The kernel of $q$ is a  discrete central subgroup in $G$ and $q(\Gamma)$ is a lattice in $H$. 
\item  Every continuous finite-dimensional real representation $\rho:G\to\GL(V)$ is of the form \[ \rho=\rho_H\circ q,\] where $\rho_H$ is the restriction to $H$ of a real algebraic representation of $\Hbf$. 

\item $q(\Gamma)$ is Zariski dense in $\Hbf$.  
\end{enumerate}

\end{claim}
We call $\Hbf$ and $H$ the \emph{algebraic model of $G$}. Sometimes we wrote $G^\alge$ as the algebraic model of $G$. The real rank of the algebraic model $H$ or $G^\alge$ is $\textrm{rank}_{\Rbb}\glie=r\ge 2$.

\subsubsection{A free abelian subgroup of full real rank}
\label{sec:cartanlattice}
Let $G,\Gamma$ be groups as in \cref{setting}.
Let $q\colon G\to H=\Hbf(\Rbb)^\circ$ be the covering homomorphism
to the algebraic model, and put $r=\operatorname{rank}_{\Rbb}H$.
Recall that $\ker q$ is discrete and central and that $q(\Gamma)$
is a lattice in $H$. We shall use a higher rank abelian subgroup of $\Gamma$.

\begin{theorem}[\cite{PrasadRaghunathan}]\label{CartanGamma}
There exist a maximal $\Rbb$-split algebraic torus $\mathbf A\subset\Hbf$
and a subgroup $\Delta\subset\Gamma$ isomorphic to $\Zbb^r$ with the
following properties. Write $A=\mathbf A(\Rbb)^\circ=\exp(\alie)$.
There is a compact torus $T_c\subset Z_H(A)$ such that the restriction
of $q$ maps $\Delta$ isomorphically onto a subgroup
$\bar\Delta\subset q(\Gamma)\cap AT_c$.
For every $\delta\in\Delta$ there is a unique factorization
\[
 q(\delta)=a_\delta m_\delta,\qquad a_\delta\in A,\quad m_\delta\in T_c.
\]
The map $\delta\mapsto\log a_\delta$ identifies $\Delta$ with
a full lattice in $\alie=\operatorname{Lie}(A)$.
\end{theorem}

We fix $A$, $T_c$, and $\Delta$ as in \cref{CartanGamma}. For
$\delta\in\Delta$, put $f_\delta=\alpha(\delta)$ and
$H_\delta=\log a_\delta$. The set $\{H_\delta:\delta\in\Delta\}$
is a full lattice in $\alie$.  Let $\mathcal U\subset\alie$ be a nonempty open cone, if $\chi_1,\ldots,\chi_s\in\alie^*$ are nonzero,
we can choose $\delta$ with $H_\delta\in\mathcal U$ and
$\chi_j(H_\delta)\ne0$ for every $j$.

\subsubsection{Applications of \cref{Borel}} Let $G$ and $\Gamma$ be as in \cref{setting}.
\begin{lemma}\label{lem:bdd}
let
$\tau:G\to\GL(W)$ be a finite-dimensional continuous representation. If $\{\tau(\gamma)(w):\gamma\in\Gamma\}$ is bounded, then $w\in W^G$, that is $\tau(g)(w)=w$. \end{lemma}
\begin{proof}
Let $W_0$ be the linear span of the orbit and choose a basis $w_i=\tau(\gamma_i)w$, $1\le i\le m$, from the orbit. For each $i$, the vectors $\tau(\gamma)w_i$ belong to the same bounded orbit. The operator norms of $\tau(\gamma^{\pm})\restrict{{W_0}}$ are therefore bounded uniformly in $\gamma$. Hence the smallest closed subgroup containing all \{$\tau(\gamma)\restrict{W_0}$ is compact in $\GL(W_0)$, hence, there is a positive definite inner product on $W_0$ preserved by $\Gamma$.

Using \cref{Borel}, we know that $\tau$ factors through $H$ and $W_0$ is $G$-invariant. Furthermore, the inner product is $G$-invariant. Thus, the resulting homomorphism from $G$ to $\GL(W_0)$ has compact image, hence trivial as $G$ is non-compact simple connected group. This implies that $G$ action on $W_0$ is trivial, especially $\tau(g)$ fixes $w\in W_0$.\end{proof}

\begin{lemma}[Real multiplicative Jordan decomposition]\label{jordan} Let $H$ be the algebraic model of $G$. For $h\in H$, write $h=h_e h_h h_u$ for its commuting elliptic, diagonalizable with positive real eigenvalue, and unipotent factors.  In every real algebraic representation $\tau$ of $\mathbf H$,  if $\tau(h)$ is hyperbolic, the stable projections of $\tau(h)$ and $\tau(h_h)$ are equal.
 \end{lemma}

\subsection{Representation theory}\label{sec2.1} 
In this subsection, we collect well known facts in representation theory. 

For finite dimensional vector spaces $V$ and $W$, we denote \[\Lcal(V,W)\coloneqq\{T:V\to W \mid T \textrm{ is $\Rbb$-linear}\},\] \[\End(V)\coloneqq\Lcal(V,V),\quad \GL(V,W)\coloneqq\{T\in\Lcal(V,W)\mid T\textrm{ is invertible}\}.\] For any two irreducible $G$-representations $(V,\rho_V)$ and $(W, \rho_W)$, we define \[\Lcal_G (V,W)=\{ T\in \Lcal(V,W):T\circ \rho_V(g)=\rho_W(g)\circ T, \forall g\in G\}\] and  \[\End_{G}(U)=\{T\in \End(L): T\circ \rho(g)=\rho(g)\circ T\quad\forall g\in G\}.\] If we consider $\Cbb$ or $\Hbb$ instead of $\Rbb$, then we put the superscript on the above notations.

\subsubsection{Completely irreducibility and Schur's lemma}
Let $V$ be a finite dimensional real vector space and $\rho:G\to \GL(V)$ be a real representation. Since $G$ is simple, we have the following: %  (Weyl's theorem
\begin{claim}[Complete irreducibility] $V$ can be written as a direct sum of irreducible representations, i.e., \[ V=\bigoplus_{i=1}^{m} d_iU_i \] where $U_i$ is a real irreducible $G$-space and $d_i$ is the number of multiplicity for all $i=1,\dots,m$. That is, $d_iU_i=U_i\oplus\dots\oplus U_i$ is the direct sum of $d_i$ copies of the representation $U_i$.
\end{claim}
Let $U$ be real irreducible sapce. Note that, since $U$ is irreducible, if a linear map $T:U\to U$ commutes with $G$-action, $T$ is invertible. This leads the following:
\begin{claim}[Schur's lemma]
Let $U$ be $G$-irreducible over $\Rbb$. Then, $\End_{G}(U)$ is finite dimensional real division algebra, hence, $\End_{G}(U)\simeq \Rbb, \Cbb,$ or $\Hbb$.
\end{claim}

\subsubsection{A decomposition of $G$-rep into irreducible representations}\label{sec2.1.2}
Let $(V,\rho:G\to\GL(V))$ be a real finite dimensional $G$-representation space and let $V=\bigoplus_{i=1}^{m} d_iU_i$ be the decomposition into irreducible components. For each $i$, let $\rho_i:G\to \GL(U_i)$ be the irreducible representation.

Let $L_i=\Lcal_G(U_i,V).$ Thus an element $e\in L_i$ is a real-linear map $e:U_i\to V$ satisfying \[e(\rho_i(g)u)=\rho(g)e(u).\]
Let $D_i=\End_{G}(U_i)^{\op}$, that is, the algebra with the same vectors as in $\End_{G}(U_i)$ and with reverse multiplication. For $d\in D$ let $\bar d\in \End_{G}(U_i)$ be its corresponding vector. If we set \[ud\coloneqq\bar d(u),\] then $U_i$ is right $D_i$-module, i.e., $(ud_1 )d_2 = u(d_1 d_2)$ for all $d_1,d_2\in D_i$.  
On the other hand, we define left $D_i$-action on $L_i\coloneqq \Lcal_{G}(U_i, V)$ by \[(dT)(u)\coloneqq T(ud) \qquad \forall d\in D_i, \quad \forall T\in L_i, \quad \forall u\in U_i.\] One can easily check that it is well-defined. Now, we can define tensor produce $U_i\otimes_{D_i}L_i$ and evaluation map $u\otimes T\mapsto T(u)\in V$, hence, the decomposition can be naturally identified by 
\[\textrm{ev}:\bigoplus_i U_i\otimes_{D_i}L_i\to  V, \qquad \textrm{ev}(u\otimes e)=e(u).\]

This leads the following decomposition into irreducible components: 
\[V= \bigoplus_{i=1}^{m} U_{i}\otimes_{D_{i}} L_i, \qquad \rho(g)=\bigoplus_{i=1}^{m}\rho_{i}(g)\otimes I_{L_{i}}\] for all $g\in G$ and $d_i = \dim_{D_{i}}L_{i}$.

\subsubsection{Complexification}\label{sec2.2.3} Here, we need to deal with not only $\Rbb$ but also complex field $\Cbb$ and quaternion algebra $\Hbb$. Let $(U,\rho_U)$ be an real irreducible $G$-representation. As before, let $D=\End_{G}(U)^{\op}$. Set \[\End^D(U)\coloneqq \{A\in \End(U):A(ud)=A(u)d\quad\forall u\in U, \forall d\in D\} =\{A\in\End(U):A\bar d =\bar d A\}.\] Then $\End^D(U)$ is an real associative algebra.
The complexification $U_\Cbb=U\otimes_{\Rbb}\Cbb$ may not be irreducible. Nevertheless, we have the following:
\begin{enumerate}
\item When $\End_{G}(U)\simeq \Rbb$, $U_{\Cbb}$ is complex irreducible and \[\End^{D}(U)\otimes_{\Rbb}\Cbb\simeq\End^{\Cbb}(U_{\Cbb}).\]
\item When $\End_{G}(U)\simeq \Cbb$, $U_{\Cbb}\simeq S\oplus \bar S$ for a complex irreducible representation $S$ (and an inequivalent irreducible representation $\bar S$) and \[\End^{D}(U)\otimes_{\Rbb}\Cbb\simeq\End^\Cbb(S)\oplus \End^{\Cbb}(\bar S).\]
\item When $\End_{G}(U)\simeq \Hbb$, $U_\Cbb\simeq S\oplus S \simeq S\otimes \Cbb^2$ for a complex irreducible representation $S$ and 
\[\End^{D}(U)\otimes_{\Rbb}\Cbb \simeq \End^{\Cbb}(S)\otimes I_{\Cbb^2}\] 
\end{enumerate}

\subsection{Root systems}\label{sec2.3} Let $G$ be as in the \cref{setting}. Without loss of generality, we assume $G$ is same with the algebraic model $H$, $G=\Hbf(\Rbb)^0$ in this subsection, especially $G$ has a finite center. For general $G$ in \cref{setting} and a finite dimensional representation $\rho$ of $G$, we can replace $G$ to its algebraic model.
\subsubsection{Restricted roots}
Let $\kappa(X,Y)=\tr(\ad(X)\circ\ad(Y))$ be the Killing form on $\glie$ and choose a Cartan involution $\theta$ on $\glie$. Then, we have the Cartan decomposition $\glie=\klie\oplus \plie$, where $\klie=\{X\in\glie: \theta X=X\}$ and $\plie=\{Y\in\glie:\theta Y= -Y\}$.

Fix a maximal abelian subspace $\alie \subset \plie$ and $A=\exp\alie$. Then $\dim_\Rbb\alie= \textrm{rank}_{\Rbb}\glie$. Furthermore, $A$ is the identity component of a maximal $\Rbb$-split torus. Let $K<G$ be the maximal compact subgroup. 

A linear function $\lambda\in\alie^*$ is called \emph{restricted roots} if \[\lambda\neq 0, \quad \glie_\lambda=\{X\in\glie: [H,X]=\lambda(H)X \textrm{ for all } H\in \alie\}\neq \{0\}.\] We denote the set of all restricted roots \[\Sigma(G,A)=\{\lambda\in\alie^*:\lambda\neq 0, \glie_\lambda\neq \{0\}\}.\] Then, we have a \emph{restricted root decomposition} as follows: 
\[\glie=\mlie\oplus \alie \bigoplus_{\lambda\in\Sigma(G,X)} \glie_{\lambda}.\] Here $\mlie=Z_{\klie}(\alie)$. Especially, $\glie_0=\mlie\oplus \alie$.

We choose an order on $\alie^*$, hence define the set $\Sigma^+(G,A)$ of \emph{positive roots} and the set $\Delta(G,A)$ of \emph{simple roots}.

Let, for $\lambda\in\alie^*$, $T_\lambda\in\alie$ be the element characterized by $\kappa(T_\lambda,H)=\lambda(H)$ for all $H\in \alie$. Note that the form $\kappa\restrict{\alie}$ is positive definite, hence $T_\lambda$ is well defined. Then, define an inner product and a norm on $\alie^*$ by $(\lambda_1,\lambda_2)\coloneqq\kappa(T_{\lambda_{1}},T_{\lambda_{2}})$ and $\|\lambda\|^2=(\lambda,\lambda)$, respectively.

For $\lambda\in\Sigma$, its \emph{coroot}, denoted by $\lambda^\vee$ or $H_{\lambda}\in\alie$, is defined by \[\lambda^\vee=H_{\lambda}=\frac{2T_{\lambda}}{\|\lambda\|^2}.\] Thus, for $\eta,\lambda\in\alie^*$,
\begin{equation}\label{eqcoroot}
\eta(\lambda^\vee)=\eta(H_\lambda)=\frac{2(\eta,\lambda)}{\|\lambda\|^2},\qquad \lambda(H_\lambda)=2.
\end{equation}

\subsubsection{Complexification}
Put $\glie_{\mathbb C}=\glie\otimes_{\Rbb}\Cbb,  \tau(X\otimes z)=X\otimes\overline z$.
Thus $\tau$ is an antilinear Lie algebra involution and its fixed real Lie algebra is the embedded copy of $\glie$. Choose a maximal abelian subalgebra $\mathfrak t\subset\mathfrak m=Z_{\mathfrak k}(\alie)$
and set
\[\hlie_{\Rbb}:=\mathfrak t\oplus\alie,\qquad \hlie:=(\hlie_{\Rbb})_{\Cbb} \subset\glie_{\Cbb}. \]
The space $\hlie$ is a complex Cartan subalgebra of $\glie_{\mathbb C}$ and contains $\alie_{\mathbb C}$.
For $\zeta\in\hlie^*$, where the dual is over $\mathbb C$, put
\[(\glie_{\mathbb C})_\zeta=\{X\in\glie_{\mathbb C}:[H,X]=\zeta(H)X \text{ for all }H\in\hlie\}.\]
The ordinary root system is $\Phi=\Phi_{\Cbb}(\glie,\hlie)=\{\zeta\in\hlie^*\setminus\{0\}: (\glie_{\mathbb C})_\zeta\ne0\}.$
Let $\kappa_{\mathbb C}$ be the complex-bilinear extension of
$\kappa$. For $\zeta\in\Phi$, define $T_\zeta\in\hlie$
and the ordinary coroot $\zeta^\vee=H_\zeta$ by
\[\kappa_{\mathbb C}(T_\zeta,H)=\zeta(H)\quad(H\in\hlie), \qquad H_\zeta=\frac{2T_\zeta}{\kappa_{\mathbb C}(T_\zeta,T_\zeta)}.\]
In particular $\zeta(H_\zeta)=2$. The real span of the ordinary coroots is \[\alie_0:=\textrm{span}_{\Rbb} \{H_\zeta:\zeta\in\Phi\}.\]

The complex Lie algebra $\glie_{\mathbb C}$ has the root decomposition
\[\glie_{\Cbb}=\hlie\oplus\bigoplus_{\zeta\in\Phi}(\glie_{\Cbb})_\zeta, \qquad \dim_{\Cbb}(\glie_{\Cbb})_\zeta=1. \]
The restriction of $\kappa_{\Cbb}$ to $\hlie$ is nondegenerate, and the ordinary coroot normalization is well defined.
One has
\[\alie_0=\alie\oplus i\mathfrak t,\qquad \hlie=\alie_0\oplus i\alie_0.\]
The latter is a direct sum of real vector spaces. Every ordinary root is real on $\alie_0$, and
\[\alie_0=\{Z\in\hlie:\zeta(Z)\in\Rbb \text{ for all }\zeta\in\Phi\}. \]

\subsection{Measurable and H\"older cocycles}
\subsubsection{Cocycle}
Let $L$ be a discrete countable group and $X$ be a space. Assume that $L$ acts on $X$, $(l,x)\mapsto l.x$. For a locally compact second countable topological group $H$, a map $\Acal:L\times X\to H$ is called a \emph{cocycle} over $L$-action on $X$, if it satisfies the following cocycle equation: \begin{equation}\label{coceq} \Acal(l_1 l_2 , x)=\Acal(l_1, (l_2).x) \Acal (l_2, x)\end{equation} for $x\in X$ and $l_1 ,l_2 \in L$. If $(X,\mu)$ is a standard probability space and $L$ acts non-singularly (e.g., $\mu$ is $L$-invariant) then $\Acal$ is called \emph{measurable cocycle} if $\Acal$ measurable, using Borel $\sigma$-algebra on $H$ and the $\sigma$-algebra $2^L$ on $L$, \eqref{coceq} satisfies for all $l_1,l_2\in L$ and $\mu$-almost every $x\in X$. Since we will only consider countable group $L$, we do not need to worry about the cocycle is strict or not (cf., \cite{Zimmerbook}). 

If $(X,d)$ is metric space, we can consider \emph{continuous or H\"older cocycle} $\Acal$ as long as $\Acal$ is continuous or H\"older and \eqref{coceq} is satisfied for all $x\in X$ and for all $l_1, l_2\in L$, respectively.

\subsubsection{Cocycle superrigidity}
Let us assume the standing assumption in \cref{setting}. A version of Zimmer's cocycle superridigity can be stated as follows:
\begin{theorem}[Zimmer's cocycle superrigidity, \cite{FisherMargulis}]\label{thm:cocsup} Let $(X,\mu)$ be a standard probability space. Assume that $\Gamma$ acts on $(X,\mu)$, $(\gamma,x)\mapsto \gamma.x$, and $\mu$ is $\Gamma$-invariant ergodic. Let $V$ be a finite dimensional real vector space. Let $\mathbf{L}$ be an algebraic group defined over $\Rbb$ and $L=\mathbf{L}(\Rbb)$.
Let $\Bcal:\Gamma\times X\to L$ be a measurable cocycle. Assume the following integrability condition: \begin{equation}\label{L2}\log\left|\left|\Bcal(\gamma,\cdot)\right|\right|\in L^{1}(X,\mu)\end{equation} Then, there exists 
\begin{enumerate}
\item a measurable map $\phi:X\to L$, 
\item a continuous representation $\rho:G\to L$,
\item a compact group $K<L$, and
\item a measurable cocycle $\kappa:\Gamma\times X\to K$,
\end{enumerate}
such that $[K,\rho(G)]=1$ and \[\Bcal(\gamma,x)=\phi(\gamma.x)^{-1}\rho(\gamma)\kappa(\gamma,x) \phi(x)\] for all $\gamma\in \Gamma$ and $\mu$-a.e. $x\in X$.
\end{theorem}
Note that, in the \cref{setting}, since $M$ is compact, any continuous cocycle, including the derivative cocycle, over $\alpha$ satisfies the integrability condition \eqref{L2}.

\subsection{Anosov diffeomorphisms}\label{sec:Anosov}
Let $M$ be a closed connected smooth manifold. Let $f:M\to M$  be a $C^{1+\eta}$ Anosov diffeomorphism with $\eta>0$.  Assume that $f$ preserves a volume measure $\mu$. 
Since we assume that $f$ has $C^{1+\eta}$ regularity, we know that $f$ is ergodic with respect to volume. 
\begin{theorem}
$f$ is ergodic with respect to $\mu$.
\end{theorem}

\subsubsection{H\"older solutions of H\"older cocycles over Anosov diffeomorphisms}
Let $E\to M$ be a finite dimensional H\"older real vector bundle and $F:E\to E$ be an affine H\"older cocycle. That is, we have the followings: 
\begin{enumerate} 
\item $F$ preserves fibers projects to $f$, $F:E_x \to E_{f(x)}$, 
\item for all $x\in M$, $F_{x}:E_{x}\rightarrow E_{f(x)}$ is given by $F_{x}(v)=L_{x}(v)+a_{x}$ for all $v\in E_x$, where $L_x \in \GL(E_x , E_{f(x)})$ and $a_x \in E_{f(x)}$, and
\item $L_{x}$ and $a_{x}$ are all H\"older in local coordinates.
\end{enumerate}
For $n\ge1$, the iterates are the compositions
\[
 F_x^{(n)}=F_{f^{n-1}x}\circ\cdots\circ F_x,
 \qquad
 F_x^{(-n)}=(F_{f^{-n}x}^{(n)})^{-1}.
\]
In particular $F_x^{(n)}:E_x\to E_{f^n x}$, and
$F_x^{(n+m)}=F_{f^m x}^{(n)}\circ F_x^{(m)}$.
For the linear cocycle generated by $L$ we use the same
superscript convention. 

An affine cocycle $F$ is called \emph{linear cocycle} if $a_x=0$ for all $x\in M$.

We can upgrade measurable conformal structure to H\"older conformal structure over $f$. More precisely; 

\begin{theorem}[{\cite[Theorem ~1.6]{Butler}}]\label{Butleroriginal} Let $\Acal$ be a H\"older continuous linear cocycle on rank $d$ vector bundle $\Ecal\to M$ of $M$ over a volume preserving $C^{1+H\"older}$ Anosov diffeomorphism $f:M\to M$. Assume that $\Acal$ preserves a measurable conformal structure $\eta:M\to \SL(d,\Rbb)/\SO(d)$, i.e., $\Acal_x[\eta_x]=\eta_{f(x)}$ for all $x\in M$. Then there is a H\"older continuous conformal structure $\wtd\eta:M\to \SL(d,\Rbb)/\SO(d)$ such that $\eta=\wtd\eta$ almost everywhere and $\Acal$-invariant.
\end{theorem}
We have variants of \cref{Butleroriginal} as follows:
\begin{proposition}\label{prop:conf}
Let $f:M\to M$  be a transitive $C^{1+\eta}$ Anosov diffeomorphism with $\eta>0$.  Assume that $f$ preserves a smooth measure $\mu$. Let $\Acal\to M$ be a finite dimensional H\"older real vector bundle and $F:\Acal\to \Acal$ be a linear H\"older cocycle. If there is a $F$-invariant measurable positive definite bilinear form $g$ on $\Acal$, that is, \[g_{f(x)}(F_x u, F_x v) =g_{x} (u,v)\] for all $u,v\in \Acal_x$ and $x\in M$, then there is a H\"older positive definite bilinear form $\wtd{g}$ on $\Acal$ such that $\wtd{g}=g$ $\mu$-a.e. and $\wtd{g}_{f(x)}(F_x u, F_x v)=\wtd{g}_{x}(u,v)$ for all $u,v\in \Acal_x$

\end{proposition}
\begin{proposition}[Regularity of a specified measurable invariant section]
\label{prop:sectionconf}
Let $f:M\to M$ be a $C^2$ transitive Anosov diffeomorphism of a
closed connected manifold, and let $\mu$ be the equilibrium measure
of a H\"older potential. Let $E\to M$ be a finite-dimensional
H\"older real vector bundle, equipped with a H\"older positive
definite metric $h$. Suppose that a H\"older linear bundle map
$F:E\to E$ covers $f$ and satisfies
\[
 h_{f(x)}(F_xu,F_xv)=h_x(u,v).
\]
If a measurable section $s$ of $E$ satisfies
\[
 s(f(x))=F_xs(x)\qquad\text{for $\mu$-almost every }x,
\]
then this same section agrees almost everywhere with a H\"older
section $\widehat s$. Its invariant-section equation holds everywhere.
No triviality of the bundle $E$ is assumed.
\end{proposition}

\subsection{Differential geometry}
Let $M$ be a smooth closed manifold of dimension $d$, and set $V=\Rbb^d$. We collect a few facts in Riemannian geometry. One can find it in standard references, for instance, \cite{KN,JohnLee}.

\subsubsection{Frames, reductions, and local representatives}

\begin{definition}[Frames and the frame bundle]
A frame at $x\in M$ is a linear isomorphism $b_x:V\to T_xM$.
The frame bundle and its projection are
\[ \textrm{Fr}(TM)=\coprod_{x\in M}     \Lcal(V,T_xM), \qquad \pi(b_x)=x.\]
 The right action is $b_x\cdot A=b_x\circ A$ for
$A\in\textrm{GL}(V)$.
\end{definition}

\begin{definition}[A principal reduction]
Let $C\subset\textrm{GL}(V)$ be compact. Write
\[ q_C:\textrm{Fr}(TM)\to\textrm{Fr}(TM)/C,\qquad \bar\pi([b_x]_C)=x.\] 
 A principal $C$-reduction is a locally trivial subbundle $P\subset\textrm{Fr}(TM)$ with each fiber a single right orbit $P_x=b_xC$. An adapted local frame is a section $b:\mathcal O\to P$ satisfying $\pi\circ b=\textrm{id}_{\mathcal O}$.
\end{definition}

Principal $C$-reductions correspond to sections $s:M\to\textrm{Fr}(TM)/C$ of $\bar\pi$ by
\[ s(x)=[b_x]_C\ (b_x\in P_x),\qquad P_x=q_C^{-1}(s(x)),\qquad P=q_C^{-1}(s(M)). \]
The reduction has continuous, H\"older, locally Lipschitz, or $C^{1+\beta}$ regularity when $s$ does. This is equivalent to the existence of adapted local frames of that regularity near every point.

\subsubsection{Affine connections and their coordinate coefficients}

\begin{definition}[A continuous affine connection]
A continuous affine connection is a local $\Rbb$-bilinear operation $(X,Y)\mapsto\nabla_XY$ from $C^1$ vector fields to continuous vector fields satisfying
\[ \nabla_{fX}Y=f\nabla_XY,\qquad \nabla_X(fY)=X(f)Y+f\nabla_XY \] for $C^1$ scalar functions $f$. Its coefficients in smooth frames are required to be continuous. H\"older or higher regularity refers to these coefficients. \end{definition}

A $C^2$ diffeomorphism $f$ preserves a continuous connection if
\[ D_xf((\nabla_XY)_x)=(\nabla_{f_*X}(f_*Y))_{f(x)}\] for $C^1$ fields $X,Y$ on their domains in $M$. Here $(f_*X)_{f(x)}=D_xf(X_x)$. The pushforward connection is defined by $(f_*\nabla)_XY=f_*\left(\nabla_{f_*^{-1}X}(f_*^{-1}Y)\right)$; invariance means $f_*\nabla=\nabla$. An action preserves $\nabla$ if each acting diffeomorphism does.

In terms of coordinates, we have the following: 
Let $q=(q^1,\ldots,q^d):\mathcal O\to U\subset\mathbb R^d$ be
a smooth chart, and set $V=\mathbb R^d$. Write
$e_i=\partial/\partial q^i$ for the coordinate vector fields
on $\mathcal O$ and $e_x=(dq_x)^{-1}:V\to T_xM$ for the
coordinate frame.

The Christoffel symbols $\Gamma^k_{ij}$ are given by
\[
 \nabla_{e_i}e_j=\sum_{k=1}^d\Gamma^k_{ij}e_k.
\]
For each $i$, collect these coefficients into the matrix
$A_i:\mathcal O\to\mathfrak{gl}(V)$ and form the matrix-valued
one-form $A$ by setting
\[
 (A_i)^k{}_j=\Gamma^k_{ij},
 \qquad
 A=\sum_{i=1}^d A_i\,dq^i.
\]
Thus $A_x(e_i(x))=A_i(x)$. In particular, the bilinear
coefficient $\Gamma_q(x):V\times V\to V$ is given by
\[
 \Gamma_q(x)(v,w)=A_x(e_xv)w,
 \qquad
 \bigl(\Gamma_q(x)(v,w)\bigr)^k
   =\sum_{i,j=1}^d\Gamma^k_{ij}(x)v^iw^j.
\]

For $C^1$ vector fields $X=\sum_i X^i e_i$ and $Y=\sum_j Y^j e_j$
on $\mathcal O$, the product rule gives
\begin{equation}\label{conncoeff}
 \nabla_XY
   =\sum_{k=1}^d
      \left(X(Y^k)+\sum_{i,j=1}^d\Gamma^k_{ij}X^iY^j\right)e_k.
\end{equation}
Equivalently, let $v:\mathcal O\to V$ be a $C^1$ function and
write $ev$ for the vector field $x\mapsto e_xv(x)$. Then for
every vector field $Z$ on $\mathcal O$,
\begin{equation}\label{connd+A}
 (\nabla_Z(ev))_x
   =e_x\bigl(dv_x(Z_x)+A_x(Z_x)v(x)\bigr).
\end{equation}
In the frame $e$, we abbreviate \eqref{connd+A} as $\nabla=d+A$.

\subsubsection{Torsion, curvature, and covariant derivatives}

Let $\nabla$ be an affine connection on the tangent bundle of a smooth $d$-dimensional manifold $M$. We assume that its coefficients in smooth local frames are continuous. All regularity assertions below are local.

\begin{definition}
For $C^1$ vector fields $X,Y$ on $M$, the torsion of $\nabla$ is \[ T^\nabla(X,Y)=\nabla_XY-\nabla_YX-[X,Y].\] The connection is torsion-free if $T^\nabla=0$. A $C^1$ vector field $Y$ is parallel if $\nabla_ZY=0$ for every smooth vector field $Z$. For a $C^1$ tangent-valued bilinear tensor $T$ and a $C^1$ Riemannian metric $g$, their covariant derivatives are given by
\begin{align*}
 (\nabla_ZT)(X,Y)&=\nabla_Z\left(T(X,Y)\right)-T(\nabla_ZX,Y)-T(X,\nabla_ZY),\\
 (\nabla_Zg)(X,Y)&=Z\left(g(X,Y)\right)-g(\nabla_ZX,Y)-g(X,\nabla_ZY).
\end{align*}
We call $T$ parallel if $\nabla T=0$. We call $\nabla$ a metric connection for $g$ if $\nabla g=0$.
\end{definition}

The torsion is a continuous section of $\wedge^2T^*M\otimes TM$. It is $C^\beta$ if the connection coefficients are $C^\beta$. If $X$ and $Y$ are parallel vector fields on an open set, then \[ [X,Y]=-T^\nabla(X,Y). \]

Let $P_c:T_xM\to T_yM$ denote parallel transport along a piecewise $C^1$ path $c$ from $x$ to $y$. For a $C^1$ tangent-valued bilinear tensor $T$, the condition $\nabla T=0$ is equivalent to \[ T_y(P_cu,P_cv)=P_c\left(T_x(u,v)\right)\] for $u,v\in T_xM$ and for every such path. Similarly, $\nabla g=0$ is equivalent to $g_y(P_cu,P_cv)=g_x(u,v)$. In particular, if $T^\nabla$ is $C^1$, then $\nabla T^\nabla=0$ is a classical equation even when the connection coefficients are only continuous.

\begin{definition}
For a $C^1$ connection and smooth vector fields $X,Y,Z$, fix the
curvature convention
\[
 R^\nabla(X,Y)Z
  =\nabla_X\nabla_Y Z-\nabla_Y\nabla_X Z-\nabla_{[X,Y]}Z.
\]
Fix smooth coordinates $(x^1,\ldots,x^d)$ on $U\subset M$ and a smooth
local frame $e_x:\mathbb R^d\to T_xM$. Write $e_b(x)$ for its $b$th
column. Define the matrices $A_i:U\to\mathfrak{gl}(d,\mathbb R)$ by
\[
 \nabla_{\partial_i}e_b=\sum_a(A_i)^a{}_b\,e_a,\qquad
 A=\sum_i A_i\,dx^i.
\]
Thus $A_x(v)$ is an endomorphism of $\mathbb R^d$ for $v\in T_xM$. If $Z=e z$, where $z:U\to\mathbb R^d$ is the column of coefficients of $Z$, then $\nabla_{\partial_i}Z=e(\partial_i z+A_i z)$, as $\nabla=d+A$. In particular, a $C^1$
frame $b=eB$ is parallel precisely when $\partial_iB=-A_iB$ for all
$i$. Denote by $F_e$ the matrix
of curvature in this frame, so that
$R^\nabla(X,Y)Z=e\left(F_e(X,Y)z\right)$.
\end{definition}

In the frame $e$, the curvature form is \[ F_e=dA+A\wedge A=\sum_{i<j}F_{ij}\,dx^i\wedge dx^j,\qquad F_{ij}=\partial_i A_j-\partial_j A_i+[A_i,A_j],\] where $[A_i,A_j]=A_iA_j-A_jA_i$.
\begin{definition}[Distributional curvature and flatness]
For continuous coefficients $A_i$, define $F_{ij}$ by the same formula,
with $\partial_iA_j$ taken in the sense of distributions. Explicitly,
in the coordinate domain, for $\varphi\in C_c^\infty(U)$,
\[
 \langle F_{ij},\varphi\rangle
  =-\int_U A_j\,\partial_i\varphi\,dx
    +\int_U A_i\,\partial_j\varphi\,dx
    +\int_U[A_i,A_j]\varphi\,dx .
\]
Here $dx$ denotes coordinate Lebesgue measure, and the integrals are
taken entry by entry. The connection is \emph{flat} if this matrix is
zero for all $i,j$ and every test function in every coordinate chart.

For the intrinsic definition under the weaker assumption of locally
square-integrable connection coefficients, see
\cite[Definition~3.3]{LeFlochMardare}.
\end{definition}

\begin{proposition}[Change of frame]
Let $e'=eC$ with $C:U\to\GL(d,\mathbb R)$ smooth. The
connection matrices $A,A'$ and the curvature matrices $F_e,F_{e'}$
satisfy
\[
 A'=C^{-1}AC+C^{-1}dC,\qquad F_{e'}=C^{-1}F_eC.
\]
The second identity holds in distributions for continuous $A$.
Consequently, distributional curvature is a distributional section of
$\Lambda^2T^*M\otimes\End(TM)$, and flatness is
independent of the chosen smooth frame.
\end{proposition}
\begin{proof}
The first identity follows by applying $\nabla$ to $eC$ column by
column. For smooth $A$, the second follows by expressing the same
curvature endomorphism in the two frames. For continuous $A$, mollify
its coefficients on a smaller coordinate ball. The resulting forms
$A^\epsilon$ converge locally uniformly to $A$. Hence
$dA^\epsilon\to dA$ in distributions and
$A^\epsilon\wedge A^\epsilon\to A\wedge A$ locally uniformly.
The transformed forms $C^{-1}A^\epsilon C+C^{-1}dC$ also converge
locally uniformly to $A'$. Applying the smooth identity and passing
to the limit gives the asserted formula.
\end{proof}

\begin{remark}[Curvature in an adapted frame]
The last formula also applies to the $C^{1+\beta}$ adapted frames used
later. Indeed, the test-function formula for $F_{ij}$ extends
continuously to $C^1$ test functions with support in a fixed compact
set. Thus $C^{-1}F_eC$ is defined entry by entry for a $C^1$ change
of frame $C$, by multiplying the test functions by the corresponding
entries of $C^{-1}$ and $C$. On a smaller ball approximate $C$ in
$C^1$ by smooth invertible matrices $C_n$. The transformed connection
forms converge locally uniformly, and the multiplied test functions
converge in $C^1$. Passing to the limit in the smooth-frame change
formula gives $F_{eC}=C^{-1}F_eC$ in distributions.
\end{remark}

\subsection{Differential equations}
\label{sec27}

\subsubsection{Mollification and distributional derivatives}

\begin{definition}[Mollification]
Choose $\varphi\in C_c^\infty(\mathbb R^d)$ satisfying
\[
 \varphi\geq0,\qquad
 \textrm{supp}\varphi\subset B(0,1),
 \qquad
 \int_{\mathbb R^d}\varphi=1.
\]
Set
\[
 \varphi_\varepsilon(x)
 =\varepsilon^{-d}\varphi(x/\varepsilon).
\]
For a locally integrable matrix function $A_i$ on an open set
$U\subset\mathbb R^d$, its mollification is defined entry by entry:
\[
 A_{i,\varepsilon}(x)
 =\int\varphi(z)A_i(x-\varepsilon z)\,dz,
 \qquad
 \textrm{dist}(x,\partial U)>\varepsilon.
\]
\end{definition}

\begin{lemma}
\label{mollif}
On every $K\Subset U$, the following assertions hold.
\begin{enumerate}
\item
If $A_i$ is continuous, then $A_{i,\varepsilon}\to A_i$
uniformly on $K$.

\item
If $A_i$ is $\beta$-H\"older on a neighborhood of $K$, then
\[
 \|A_{i,\varepsilon}-A_i\|_{C^0(K)}
 \leq C\varepsilon^\beta.
\]

\item
If $A_i\in L^p_{\mathrm{loc}}$, $1\leq p<\infty$, then
$A_{i,\varepsilon}\to A_i$ in $L^p(K)$.

\item
Distributional differentiation commutes with convolution:
\[
 \partial_jA_{i,\varepsilon}
 =(\partial_jA_i)*\varphi_\varepsilon.
\]
\end{enumerate}
Uniform convergence of derivatives is not asserted for a merely
continuous or H\"older function.
\end{lemma}

\subsubsection{Existence of integral curves and uniqueness}

\begin{theorem}[Peano and Picard--Lindel\"of]
\label{ode:existence}
Let $\mathcal U\subset\Rbb\times\Rbb^m$ be open and let
$f:\mathcal U\to\Rbb^m$ be continuous. For every
$(t_0,z_0)\in\mathcal U$, the problem
\[
 z'(t)=f(t,z(t)),\qquad z(t_0)=z_0
\]
has a $C^1$ solution on a neighborhood of $t_0$.

If $f$ is locally Lipschitz in $z$, (locally uniformly in $t$), this solution is unique on its interval of existence. A sufficient condition is that $\partial_z f$ exists and is continuous.
\end{theorem}

\subsubsection{Gronwall's inequality}
\begin{lemma}[Integral Gronwall inequality] \label{gronwallineq}
Suppose $u:[0,T]\to[0,\infty)$ is continuous, $b:[0,T]\to[0,\infty)$ is integrable, and $a:[0,T]\to[0,\infty)$ is nondecreasing. If \[ u(t)\le a(t)+\int_0^t b(s)u(s)\,ds,\] then \[ u(t)\le a(t)\exp\left(\int_0^t b(s)\,ds\right).\] More precisely, if $a_0\ge0$, $g\ge0$ is integrable, and
\[ u(t)\le a_0+\int_0^t\left(b(s)u(s)+g(s)\right)\,ds,\]
then
\begin{equation}
 u(t)\le a_0\exp\left(\int_0^t b(s)\,ds\right) +\int_0^t g(s)\exp\left(\int_s^t b(r)\,dr\right)\,ds. \label{gronwallf}
\end{equation}
\end{lemma}

\section{Continuous compact group reduction $P$ and the adapted metric}\label{sec3} In this section, we will show the following: 
\begin{theorem}[Continuous compact reduction]\label{thm:reduction}
Under the \cref{setting}, there exists an algebraic representation $\rho:G\to\SL(V)$ with $\dim V= \dim M$, a compact group \[C\subset \SO(V)\cap \mathcal{Z}_{\SL(V)}(\rho(G)),\]  and a global H\"older principal $C$-reduction $P$ of $\Fr(TM)$ such that \begin{equation}\label{eq:superrigid}D_{x}\alpha(\gamma)P_{x}=P_{\alpha(\gamma)(x)}\rho(\gamma)\end{equation} for all $\gamma\in \Gamma$ and $x\in M$.
Moreover, $\rho$ has no zero restricted weight and there is an H\"older metric $h$ on $M$ that gives a H\"older Lyapunov splitting which is same with weight space decomposition of $V$ with respect to $\rho$ through $P$ for a lattice Cartan subgroup $\Delta$ in $\Gamma$.
\end{theorem}
Recall that a lattice Cartan subgroup $\Delta$ is fixed in \cref{CartanGamma}. Equation \eqref{eq:superrigid} is the equivariance relation. It can be read that, for all $b_x\in P_x$ and $b_{\alpha(\gamma)(x)}\in P_{\alpha(\gamma)(x)}$, there is a unique $k\in C$ such that
\[ D_x\alpha(\gamma)b_x =b_{\alpha(\gamma)(x)}\rho(\gamma)k. \]
It also has a precise quotient formulation. Let $s_P:M\to \Fr(TM)/C$ be the section corresponding to the reduction $P$. For $b_x\in P_x$, define
\[\widehat\alpha_\gamma(b_x)=D_x\alpha(\gamma)\circ b_x\circ\rho(\gamma)^{-1}:V\to T_{\alpha(\gamma)(x)}M.\]
Because $C$ commutes with $\rho(\Gamma)$, $\widehat\alpha_\gamma(bc)=\widehat\alpha_\gamma(b)c$, this defines an action on the quotient bundle $\Fr(TM)/C$. Our identity \eqref{eq:superrigid}  is exactly equivalent to $s(\gamma x)=\widehat\alpha_\gamma(s(x))$.

\subsection{Cocycle superrigidity and measurable frame} Zimmer's cocycle superrigidity \cref{thm:cocsup} already gives the measurable version of \cref{thm:reduction}
\begin{proposition}\label{prop:msblframe}
Under the \cref{setting}, there exists an algebraic representation $\rho:G\to\SL(V)$ with $\dim V= \dim M$, a compact group \[K\subset \mathrm{SO}(V)\cap \mathcal{Z}_{\SL(V)}(\rho(G)),\]   a measurable family of linear isomorphisms  $B_x:V\to T_{x}M$ with $|\det B_x| =1$,  and a measurable cocycle $\kappa:\Gamma\times M\to K$ such that \[D_{x}\alpha(\gamma)B_{x}=B_{\alpha(\gamma)(x)}\rho(\gamma)\kappa(\gamma,x)\] for all $\gamma\in \Gamma$ and $\mu$ almost every $x\in M$. 

\end{proposition}

Since we have the Anosov diffeomorphism $\alpha(\gamma_0)$ let $TM=E^s_0\oplus E^u_0$ be the decomposition into H\"older stable and unstable distributions with respect to $\alpha(\gamma_0)$. Then, the following is directly from \cref{prop:msblframe}:
\begin{claim}\label{claim:matching}
$\rho(\gamma_0)$ is hyperbolic, that is, there is no eigenvalues of modulus $1$. Let $V=V^s_0\oplus V^u_0$ be the $\rho(\gamma_0)$-invariant decomposition into subspaces $V^s_0$ and $V^u_0$ where $V^s_0$ is the direct sum of generalized eigenspaces that has eigenvalue with modulus smaller than $1$ and $V^u_0$ is the direct sum of generalized eigenspaces that has eigenvalue with modulus larger than $1$. Then, \[B_{x}V^s_0=E^s_0(x), \qquad B_x V^u_0=E^u_0(x)\] for $\mu$-almost every $x\in M$.
\end{claim}

The $\rho$ in \cref{thm:reduction} will be the same $\rho$ with \cref{prop:msblframe}. Now, from an Anosov diffeomorphism, we can easily show that $\rho$ does not have zero weight:
\begin{proposition}\label{prop:nozero} The superrigidity homomorphism $\rho$ does not have zero restricted weight. \end{proposition}
\begin{proof}
Since $\alpha(\gamma_0)$ is Anosov diffeomorphism. By \cref{claim:matching}, $\rho(\gamma_0)$ is hyperbolic. Using Jordan decomposition,  $\gamma_0=\gamma_e\gamma_h\gamma_u$, where $\gamma_e$, $\gamma_h$, $\gamma_u$ commutes, $\gamma_e$ is elliptic (that is, semisimple with all eigenvalues of modulus $1$), $\gamma_h$ is $\Rbb$-diagonalizable with positive real eigenvalues, and $\gamma_u$ is unipotent (that is, all its eigenvalues are equal to $1$). Note that this decomposition is not in $\Gamma$, i.e., $\gamma_h$ may not be in $\Gamma$. However, the representation $\rho$ is defined on $G$, hence we can still consider $\rho(\gamma_h)$.  After conjugation, we may assume $\gamma_h\in A$ for a maximal $\Rbb$-split torus $A$. Note that this decomposition is preserved by the algebraic representation $\rho$. Hence, eigenvalues of $\rho(\gamma_u)$ and $\rho(\gamma_e)$ are all modulus $1$. Let \[V_0=\{v\in V : \rho(\exp X)v = v,\quad \forall X\in \alie\}\] be the zero weight subspace. Assume that $V_0\neq\{0\}$. Since $\gamma_h\in A$,  \[ V_0\subset \ker(\rho(\gamma_h)-I_V),\] so $\ker (\rho(\gamma_h)-I_V) \neq\{0\}$. Both $\rho(\gamma_e)$ and $\rho(\gamma_u)$ preserve $\ker(\rho(\gamma_h)-I_V)$, because they commute with $\rho(\gamma_h)$. Let $W=\ker(\rho(\gamma_h) -I_V)$. On $\ker(\rho(\gamma_h)-I_V)$, \[\rho(\gamma_0)\restrict{W}=\rho(\gamma_e)\restrict{W}\rho(\gamma_u)\restrict{W}.\]
Since $\rho(\gamma_e)$ and $\rho(\gamma_u)$ commutes, their product has still eigenvalues with modulus $1$. This contradicts the hyperbolicity of $\rho(\gamma_0)$. Thus, $V_0=\{0\}$.
Hence, there is no zero weight.
\end{proof}

In the following sections, we use the decomposition into irreducible representations introduced in \cref{sec2.1.2}.
\[V= \bigoplus_{i=1}^{m} U_{i}\otimes_{D_{i}} L_i, \qquad \rho(g)=\bigoplus_{i=1}^{m}\rho_{i}(g)\otimes I_{L_{i}}\] for all $g\in G$ and $d_i = \dim_{D_{i}}L_{i}$. Let $\dim U_i = r_i$.

\subsection{Outline of proof of \cref{thm:reduction}} Since the proof is quite involved, let us explain the outline of the proof here.

Let us focus on the simple case.  
Assume that there is no compact error $\kappa$, that is $K=\{id\}$, and there are only two inequivalent real irreducible components in $V$, that is,  \[V=U_{1}\otimes_{\Rbb} \Rbb \oplus U_{2}\otimes_{\Rbb} \Rbb.\]  Let $\dim U_i=r_i$. In this case, we indeed construct H\"older continuous (global) $\Gamma$-invariant frame $b:M\to \Fr(TM)$.

Fix basis of $U_1$ and $U_2$ as $\{e_{1,1},\dots, e_{1,r_1}\}$ and $\{e_{2,1},\dots,e_{2,r_2}\}$, respectively. 
\cref{thm:samealg} says that we can recover the projection into $\Rbb e_{i,a}$ for all $i,a$ along other basis from stable projections of $\rho(\gamma\gamma_0\gamma^{-1})$ for all $\gamma\in\Gamma$. Here stable projection means, for instance, for $\rho(\gamma_0)$, the projection into $V_0^s$ along $V_0^u$. Using \cref{claim:matching} and the fact that stable and unstable distributions of $\alpha(\gamma_0)$ are H\"older, we can produce H\"older sections $E_{i,a}:M\to \End(TM)$ where $i=1,2$ and $a\in \{1,\dots, r_i\}$ that are same with the projection into the line $\Rbb e_{i,a}$ in $V$ almost everywhere, i.e., \[E_{i,a}=B_{x}\textrm{Proj}_{e_{i,a}}B_{x}^{-1}\] for $\mu$ almost every $x\in M$, where $\textrm{Proj}_{e_{i,a}}$ is the projection into $e_{i,a}$ along other basis. In \cref{thm:multiplicityspace}, we also separately produce the H\"older orthonormal basis on each line bundle $\Im E_{i,a}\subset TM$, for all $i=1,2$ and $a=1,\dots, r_i$, using \cref{prop:conf}. Basically, \cref{ConfAcal,prop:conf} asserts the existence of H\"older conformal structure on each line bundle, hence, the existence of H\"older orthonormal basis. Here is the one of places we used the existence of Anosov diffeomorphism. This is enough to provide a H\"older frame $\wtd{B}_{x}$. Since we consistently construct our H\"older objects measurably (through $B_x$) same with the corresponding objects in the model vector space $V$,  we have $\wtd B_x = B_x$ $\mu$-almost every $x$. As $B_x$ comes from the cocycle superrigidity,  for all $\gamma\in\Gamma$,  we have \begin{equation}\label{frameeq} D_{x}\alpha(\gamma)\wtd B_x =\wtd B_{\alpha(\gamma)(x)} \rho(\gamma)\end{equation} for $\mu$-almost every $x\in M$, hence,  by continuity of $\wtd B$ and $\supp\mu=M$. After choosing a H\"older frame $\wtd B:M\to \Fr(TM)$ with \cref{frameeq}, we can define a H\"older metric from an inner product on $V$ using $\wtd B$.

Note that, if there is a multiplicity, that is, $d_i>1$, we can not distinguish two copies of sub-irreducible representations $U_i$, hence, we only can get the simultaneous projection into the multiplicity number of copies of (real, complex, or quoternionic) lines. This is the reason why we consider the tensor decomposition we introduced in \cref{sec2.1.2} instead of direct sum of same copies and  why we verify H\"older conformal structure on the multiplicity space $E_i$ in \cref{thm:multiplicityspace}.

\subsection{Recovering Irreducible component $U_i$}
We start with the following property in representation theory. Let $\rho:G\to \GL(V)$ be the representation in \cref{prop:msblframe}. Recall that we have the decoposition \[V=\bigoplus_{i=1}^{m}U_{i}\otimes_{D_{i}}L_{i},\quad D_{i}=\End_{G}(U_{i})^{\op}\]  into irreducible representations as in \cref{sec2.1.2}. Note that, for $\gamma\in\Gamma$, $\rho(\gamma)$ can be identified with $\bigoplus_i \rho_i (\gamma) \otimes I_{L_{i}}$ in $\End(V)$.  Since $\alpha(\gamma_0)$ is Anosov, $\rho(\gamma_0)$ is hyperbolic, i.e., there is no complex eigenvalues of modulus $1$. Let $V=V^s_0\oplus V_0^u$ be the decomposition into $\rho(\gamma_0)$ invariant subspaces with eigenvalues of modulus smaller than $1$ (the \emph{stable space}) and eigenvalues of modulus bigger than $1$ (the \emph{unstable space}). Let $P_0$ be the projection onto $V_0^s$ along $V_0^u$. In general, for a hyperbolic matrix $A\in\GL(V)$, we can define the projection to the stable space along the unstable space. We call this projection the \emph{stable projection}.  

We can characterize what is the algebra generated by conjugations of $P_0$.
\begin{proposition}\label{thm:samealg}
Under the above setting, let \[\Acal=\textrm{alg}_{\Rbb}\{I_V, \rho(\gamma)P_0 \rho(\gamma)^{-1} : \gamma\in \Gamma\}\] be the smallest finite dimensional real unital algebra generated by $I_V$ and $\rho(\gamma) P_0 \rho(\gamma)^{-1}$ for all $\gamma\in \Gamma$. Let \[\Bcal=\bigoplus_{i=1}^{m}\End_{D_{i}}(U_i) \otimes I_{L_{i}}.\]   Then \[\Acal=\Bcal\subset \End (V).\]

\end{proposition}
Recall that $\rho$ factors through the algebraic model of $G$. 
\subsubsection{Proof of \cref{thm:samealg}}
Firstly, we show that $\Acal \subset \Bcal\subset \End V$. For a hyperbolic matrix $A$, let stable space $A^s$ be the direct sum of generalized eigen-spaces with eigenvalue modulus smaller than 1. Similarly define the unstable space $A^u$.  Since stable and unstable spaces preserved by $A$, the characteristic polynomial is a product of two coprime polynomials corresponding characteristic polynomials on restricted matrix for each spaces. Using Cayley--Hamilton theorem, we have the following:
\begin{claim} For a hyperbolic matrix $A\in\GL(V)$, the stable projection of $A$ is a real polynomial of $A$. Especially, the projection $P_0$ is a real polynomial in $\rho(\gamma_0)$.
\end{claim}
As $\rho(\gamma)$ is identified with $\bigoplus_i\rho_i(\gamma)\otimes I_{L_i}$ for all $\gamma\in\Gamma$, the above claim implies $\Acal\subset \Bcal$.

Recall the notations in \cref{sec2.0,sec2.1}.
Set $V_{\Cbb}=V\otimes_{\Rbb}\Cbb$. The complexification of the action on $V$ is the representation
\[ \rho_{\Cbb}:G\longrightarrow\GL_{\Cbb}(V_{\Cbb}), \qquad \rho_{\Cbb}(g)(v\otimes z)=\rho(g)v\otimes z. \]
The algebraic representation $\rho_{\Hbf}$ also extends by scalar extension to $\rho_{\Hbf,\Cbb}:\Hbf(\Cbb)\to\GL_{\Cbb}(V_{\Cbb})$.
These representations are related by  $\rho_{\Cbb}=\rho_{\Hbf,\Cbb}\circ\iota\circ q,$ where $\iota:\Hbf(\Rbb)^\circ\hookrightarrow\Hbf(\Cbb)$ is inclusion.

Recall that, for a group $K$, $\Lcal_K^{\mathbb F}(W_1,W_2)$ denotes the space of $\mathbb F$-linear intertwining maps between the indicated $K$-representations.

\begin{claim}\label{distinctreal}
If $U_i$ and $U_j$ are nonisomorphic real irreducible $G$-representations, then $(U_i)_{\Cbb}$ and $(U_j)_{\Cbb}$ have no common complex irreducible components.
\end{claim}

\begin{proof}
A complex-linear map intertwines the actions of $G$ through $\rho_\Cbb$ if and only if it intertwines the actions of $\Hbf(\Cbb)$ through $\rho_{\Hbf,\Cbb}$. One may see this in the lie algebra level.

Moreover, in real bases, a complex-linear map between the complexifications has the form $A+iB$, where $A$ and $B$ are real-linear maps. It intertwines the $G$-actions if and only if both $A$ and $B$ intertwine the real $G$-actions. Hence \[\Lcal_{\Hbf(\Cbb)}^{\Cbb}((U_i)_{\Cbb},(U_j)_{\Cbb})\simeq\Lcal_G^{\Rbb}(U_i,U_j)\otimes_{\Rbb}\Cbb=0\] where the last equality follows from the fact that $U_i$ is inequivalent to $U_j$. 

Suppose that a complex irreducible representation $S$ occurred in both complexifications. Complete reducibility would give an equivariant projection $p:(U_i)_{\Cbb}\to S$ and an equivariant inclusion $j:S\longrightarrow (U_j)_{\Cbb}$. The composition $j\circ p$ would be a nonzero intertwining map, hence, contradiction.
\end{proof}

Recall the real irreducible decomposition and the definition of $\Bcal$:
\[V\cong\bigoplus_{i=1}^{m} U_i\otimes_{D_i}L_i,\qquad D_i=\End_G^{\Rbb}(U_i)^{\mathrm{op}},\qquad\Bcal=\bigoplus_{i=1}^{m}\End_{D_i}(U_i)\otimes I_{L_i}.\]
We regard
\[\Bcal_{\Cbb}:=\Bcal\otimes_{\Rbb}\Cbb\subset\End_{\Rbb}(V)\otimes_{\Rbb}\Cbb\cong\End_{\Cbb}(V_{\Cbb}).\]

Using complete reducibility for $V_\Cbb$, \cref{sec2.2.3,distinctreal}, we can conclude the following:

\begin{claim}\label{cpxalg}
Let $S_1,\ldots,S_t$ represent all the inequivalent complex
irreducible components of $V_{\Cbb}$. There are positive
integers $m_1,\ldots,m_t$ and identifications
\[V_{\Cbb}\cong\bigoplus_{j=1}^tS_j\otimes_{\Cbb}\Cbb^{m_j},\qquad\Bcal_{\Cbb}=\bigoplus_{j=1}^t\End^{\Cbb}(S_j)\otimes I_{\Cbb^{m_j}}\simeq \bigoplus_{j=1}^{t} \End^{\Cbb}(S_j).\]
\end{claim}

Note that if we restrict $P_0$ to each irreducible component $U_i$ then the restricted projection $P_0\restrict{U_i}$ is not zero and is non-trivial as $\det \rho_i(\gamma_0)=1$. Similarly, on every complex irreducible component $S_j$ of $V_{\C}$, the restricted stable projection $(P_0)_{\Cbb}\restrict{S_j}$ is different from zero and identity. % Indeed, $
\begin{claim}
On every complex irreducible components $S_j$ of $V_{\C}$, the restricted stable projection is different from zero and identity.
\end{claim}

Furthermore, by the definition of $\Acal$ and Zariski density \cref{Borel}, $\Acal$ is preserved by conjugation by $\rho(G)$;
\begin{claim}
For all $g\in G$, we have  $\rho(g)\Acal\rho(g)^{-1}=\Acal$ and, for all  $X\in\glie$, $[d\rho(X),\Acal] \subset\Acal$. After complexification, $[d\rho_\Cbb(Z),\Acal_\Cbb]\subset \Acal_\Cbb\subset \End V_\Cbb=\End_\Cbb V$ for all $Z\in \glie_\Cbb$.
\end{claim}

Let us fix a complex irreducible component $S$ in $V_\Cbb$. Let $S\otimes \Cbb^m$ be the summand in $V$ corresponding $S$. Let $\Dcal_S\subset \End_\Cbb S$ be the image of $\Acal_\Cbb$ into $\End_\Cbb S$ through the $G$-equivariant map $\End_\Cbb V_\Cbb\twoheadrightarrow \End_\Cbb S \otimes I_{\Cbb^{m}}\simeq \End_\Cbb S$ where $G$ acts by conjugation.  We choose the connected simply connected compact group $K$ whose Lie algebra is the compact real form of $\glie_\Cbb$. Then, by integrating $d\rho$, we can consider $V$ as a $K$-representation and $S$ as an $K$-irreducible component of $V$. Furthermore, $\Acal_\Cbb$ is preserved by a conjugation by $K$, hence, $\Dcal_S$ is preserved by a  conjugation by $K$.

We will use the following adapted form of \cite[Thm.~ 4.3, Cor.~4.4]{Sage} in our setting. 
\begin{theorem}[{\cite[Thm.~ 4.3, Cor.~4.4]{Sage}}]\label{sage}
Let $K=K_1\times\cdots\times K_l$ be a product of connected simply connected compact groups with simple Lie algebras. Let $S$ be an irreducible $K$-representation. Let us write $S=S_1\otimes_{\C}\cdots\otimes_{\C}S_l$ where each $S_j$ is an irreducible complex representation of $K_j$, possibly trivial. Let $I=\{j:S_j\text{ is nontrivial}\}$. The unital complex associative subalgebras of $\End_{\C}(S)$ which is preserved by conjugation by $K$ are exactly
\[
 \bigotimes_{j=1}^l B_j,\qquad  B_j=\begin{cases}\End_{\Cbb}(S_j),&j\in J,\\  \Cbb I_{S_j},&j\notin J, \end{cases} \quad J\subset I.
\]
\end{theorem} Recall that every $K$-irreducible representation $S$ can be written as a tensor product of $K_j$-irreducible representations as in the above theorem (\cite[Ex.~2.36]{FultonHarris}). 
Note that, in our case, $K$ has at most $2$ simple factors; If $\glie_\Cbb$ is simple then $K$ itself is simple. If $\glie_\Cbb=\slie\oplus \bar\slie$, i.e., when we consider complex group as a real lie group, then $K$ has at most $2$ simple factors.

\begin{proposition}\label{twofactor}
$\mathcal D_S=\End_{\C}(S)$ for every irreducible component $S$.
\end{proposition}
\begin{proof}[Proof of \cref{twofactor}] We divided into two cases.

\emph{Case I. When $ \glie_\Cbb$ is simple:}  In this case, $K$ is simple by itself.  \cref{sage} gives exactly two unital $K$-invariant subalgebras in $\End_\Cbb S$:
\[ \Cbb I_S,\qquad\End_\Cbb(S).\]

The algebra $\Dcal_S$ contains the stable projection $P_S=(P_0)_\Cbb\restrict{S}\neq 0, I_S$. Also, $P_S^2=P_S$. A scalar idempotent is zero or identity, so the scalar algebra is excluded. Consequently \[\Dcal_S=\End_\Cbb(S).\]

\emph{Case II. When $\glie_\Cbb$ is not simple:} 
Let $G^\alge$ be the algebraic model of $G$.

Recall that $\glie=\klie\oplus \plie$ is the Cartan decomposition and $\alie\subset \plie$ is the maximal abelian subalgebra. 

$\rho(\gamma_0)$ is hyperbolic since $\alpha(\gamma_0)$ is Anosov. Using Jordan decomposition, we have $\gamma_0=\gamma_e\gamma_h\gamma_u$, where $\gamma_e$, $\gamma_h$, $\gamma_u$ commutes, $\gamma_e$ is elliptic $\gamma_h$ is $\Rbb$-diagonalizable with positive real eigenvalues, and $\gamma_u$ is unipotent. As we saw before in the proof of \cref{prop:nozero}, $\rho(\gamma_h)$ is hyperbolic and $\rho(\gamma_h)$ and $\rho(\gamma_0)$ have same stable projection as $\gamma_e,\gamma_h$, and $\gamma_u$ commutes. After a  conjugation (which preserves $\Dcal_S$ and two intermediate subalgebras), we may assume $\gamma_h\in A$ for a maximal $\Rbb$-split torus $A$, hence $Z_\gamma\coloneqq \log \gamma_h\in\alie$. 

Recall a few notations introduced in \cref{sec2.3}. Let $\hlie \subset \glie_\Cbb$ be the Cartan subalgebra defined by $(\alie\oplus \mathfrak{t})_{\Cbb}=\hlie$. Here, $\mathfrak{t}\subset\mlie=Z_\klie(\alie)$ be a maximal abelian subalgebra. Let  $\tau: \glie_\Cbb\to \glie_\Cbb$ be the conjugation map, $\tau(X\otimes z)=X\otimes\bar z$. Then, we are in the case of $\glie_{\Cbb}=\slie\oplus \tau(\slie)$ where $\slie$ and $\tau(\slie)$ are all complex simple ideals. 
Let $\slie_1\coloneqq \slie$ and $\slie_2\coloneqq \tau(\slie)$.  By the definition of $\tau$, $\glie$ be the set of all $\tau$-fixed vectors.

We can characterize Cartan subalgebras of $\slie_j$ as follows: 
Let $\hlie_1=\textrm{pr}_1\hlie$ and $\hlie_2=\textrm{pr}_2\hlie$ where $\textrm{pr}_j:\glie_\Cbb\to \slie_j$ is the projection along the decomposition. Similarly, $\alie_{0,j}=\textrm{pr}_j (\alie)$. 
Then, $\hlie_j=\alie_{0,j} \oplus i\alie_{0,j}$ and $\alie_{0,j}$ is the real span of coroots of $(\slie_j,\hlie_j)$ for $j=1,2$. 
With these notations, for every $Z\in \alie$ can be written as $Z=H+\tau H$ where $H\in \alie_{0,1}$. Conversely, for every $H\in\alie_{0,1}$, $H+\tau H\in \alie$. 

Set $Z_\gamma=H_\gamma+\tau H_\gamma$ for $H_\gamma\in \alie_{0,1}$.

Now we have $\glie_\Cbb=\slie_1\oplus\slie_2$ and $K=K_1\times K_2$. Write $S=S_1\otimes_\Cbb S_2$ such that  $K_i$-(as well as $\slie_i$-)irreducible $S_i$. Let $\pi_i :\slie_i \to \End^\Cbb(S_i)$ be Lie algebra representations, $d\rho_\Cbb=\pi_1\otimes I_{S_2}\oplus I_{S_{1}}\otimes \pi_2$.

If one factor is trivial, then we can use the same argument with simple $\glie_\Cbb$ case.  If both are nontrivial,  \cref{sage} gives that $\Dcal_S$ is one of the following four subalgebras:
\[ \Cbb I_S,\qquad \End_\Cbb(S_1)\otimes I_{S_2},\qquad I_{S_1}\otimes\End_\Cbb(S_2),\qquad \End_\Cbb(S). \]

The scalar algebra $\Cbb I_S$ can be excluded as before. We will exclude intermediate cases, $\End_\Cbb(S_1)\otimes I_{S_2}$ and $I_{S_1}\otimes\End_\Cbb(S_2)$.

Choose a Lie algebra involution $c:\slie_1 \to \slie_1$ which fixes $\alie_{0,1}$ pointwise and $c(zX)=\bar{z}c(X)$ for $z\in \Cbb$ and $X\in \slie_1$. For instance, choose a basis of $\slie$ consists of simple coroots and normalized root vectors with real bracket coefficients. Then, one can define $c$ as a map take the complex conjugation of coefficients. 

Define $\phi=\tau\restrict{\slie_1}\circ c:\slie_1\to\slie_2$. As $\tau$ and $c$ are all anti-linear, $\phi$ is $\Cbb$-linear. Also, $c$ fixes $\alie_{0,1}$, $\phi\restrict{\alie_{0,1}}=\tau$. Let \[\sigma_1\coloneqq \pi_1 : \slie_1 \to \End^\Cbb(S_1) \quad \textrm{ and } \quad \sigma_2\coloneqq \pi_2 \circ \phi :\slie_1\to \End^\Cbb(S_2).\] Both $\sigma_i$ are all well-defined and non-trivial complex Lie algebra representations of $\slie_1$. 

For every $Z=H+\tau H\in\alie$, $d\rho_\Cbb(Z)\restrict{S}$ can be written as \[d\rho_\Cbb (Z)\restrict{S} = \sigma_1 (H)\otimes I_{S_2}+I_{S_1}\otimes \sigma_2(H).\]
Hence, $\exp(d\rho_\Cbb(Z)\restrict{S})$ is given by  $e^{\sigma_1 H} \otimes e^{\sigma_2 H}$.  

Recall that $Z_\gamma=H_\gamma+\tau H_\gamma$ for $H_\gamma\in\alie_{0,1}$. Let $T=\exp(d\rho_\Cbb(Z_\gamma)\restrict{S})$. Let $\Lambda_j\subset \hlie_1^*$ be the set of weights of $\sigma_j$ for $j=1,2$. Here $\hlie_1^*$ is the set of $\Cbb$-linear functionals. However, every $\lambda\in\Lambda_j$ takes integer values on the coroots (using $\mathfrak{sl}_2$-triple). Since $\alie_{0,1}$ is their real span, for every $H\in\alie_{0,1}$, $\lambda(H)\in\Rbb$ for all $\lambda \in \Lambda_j$.

Hence, we can find $\eta_1$ and $\eta_2$ in $\Lambda_2$ so that $\eta_1(H_\gamma)\le \eta(H_\gamma) \le \eta_2(H_\gamma)$ for all $\eta\in\Lambda_2$. We claim the followings:
\begin{claim}\label{eta<0} $\eta_1(H_\gamma)<0<\eta_2(H_\gamma)$. Furthermore, for every $\beta\in\Phi_{\Cbb}(\slie_1,\hlie_1)$,  
\[|\beta(H_\gamma)|\le\eta_2(H_\gamma)-\eta_1(H_\gamma).\] \end{claim}

\begin{claim}\label{lambda1mid} There is $\lambda\in\Lambda_1$ such that $-\eta_2(H_\gamma)<\lambda(H_\gamma) <-\eta_1(H_\gamma)$.\end{claim}

Let us postpone the proof of \cref{eta<0,lambda1mid} and show that how we can exclude intermediate possibilities of $\Dcal_S$ first. Choose non-zero vectors $v_\lambda\in (S_1)_{\lambda}$, $w_1\in (S_2)_{\eta_1}$, and $w_2\in (S_2)_{\eta_2}$.
\cref{lambda1mid} implies that $\lambda(H_\gamma)+\eta_1(H_\gamma)<0$ and $\lambda(H_\gamma)+\eta_2(H_\gamma)>0$.  As \[T(v_\lambda\otimes w_1)=e^{\lambda(H_\gamma)+\eta_1(H_\gamma)}v_\lambda\otimes w_1,\qquad T(v_\lambda\otimes w_2)=e^{\lambda(H_\gamma)+\eta_2(H_\gamma)}v_\lambda\otimes w_2,\] we have \begin{equation}\label{P_T}P_T(v_\lambda\otimes w_1)= v_\lambda\otimes w_1,\qquad P_T(v_\lambda\otimes w_1)=0.\end{equation}

Now if $\Dcal_S=\End^\Cbb(S_1)\otimes I_{S_{2}}$, then $P_T=Q\otimes I_{S_{2}}$ for some $Q\in \End^\Cbb(S_1)$. Then, $Q(v_\lambda)=v_\lambda$ from the first equation in \eqref{P_T} and $Q(v_\lambda)=0$ from the second equation in \eqref{P_T}, hence, contradiction. Therefore, $\Dcal_S\neq \End^\Cbb(S_1)\otimes I_{S_{2}}$. Similarly, using the same argument with changing the role of $\slie_1$ and $\slie_2$, $\Dcal_S\neq I_{S_{1}}\otimes \End^\Cbb(S_2)$. Therefore, $\Dcal_S=\End^\Cbb(S)$. 
 \end{proof}

 Let us prove \cref{eta<0,lambda1mid} used in the proof.
 \begin{proof}[Proof of \cref{eta<0}]
 Since $\rho(\gamma_h)$ is hyperbolic, $H_\gamma\neq 0$ as well as $\sigma_2(H_\gamma)\neq 0$. Indeed, $\ker \sigma_2\subset \slie_1$ is a ideal but $\slie_1$ is simple and $\sigma_2$ is non-trivial, hence $\ker\sigma_2=\{0\}$ so that $\sigma_2(H_\gamma)\neq0$. Note that $\tr \sigma_2(H_\gamma) =0$ as $\slie_1=[\slie_1,\slie_1]$ and the trace vanishes on the commutator. Since $\sigma_2(H_\gamma)$ is real diagonalizable and its real eigenvalues are $\eta(H_\gamma)$ for $\eta\in\Lambda_2$, $\eta_1(H_\gamma)<0<\eta_2(H_\gamma)$.
 
 Now, let us show the second assertion. Choose $\beta\in \Phi_\Cbb (\slie_1,\hlie_1)$ and a root vector $X_\beta\in \slie_1$. As $\ker\sigma_2=\{0\}$, $\sigma_2(X_\beta)\neq 0$. Hence, there exists $\eta_0\in\Lambda_2$ such that $\sigma_2(X_\beta)(S_{\eta_0})$ is non-zero. Since $\sigma_2(X_\beta)(S_{\eta_0}) \subset S_{\eta_0+\beta}$ whenever $\eta_0+\beta\in \Lambda_2$ or $\sigma_2(X_\beta)(S_{\eta_0})=0$, we can deduce that $\eta_0+\beta\in\Lambda_2$. Therefore, $\eta_1(H_\gamma)\le \eta_0(H_\gamma), (\eta_0 +\beta)(H_\gamma)\le \eta_2(H_\gamma)$ so that $|\beta(H_\gamma)|=|(\eta_0+\beta)(H_\gamma)-\eta_0(H_\gamma)| \le \eta_2(H_\gamma)-\eta_1(H_\gamma)$. 
  \end{proof}

 \begin{proof}[Proof of \cref{lambda1mid}]
 If $0\in\Lambda_1$ then we can just choose $\lambda=0$. 
 
Suppose now that $0\notin\Lambda_1$. Fix a Weyl group invariant inner product  on $\hlie_1$. Choose a weight $\delta\in\Lambda_1$ with a minimum norm and in the positive Weyl chamber. 

Then, for all positive roots $\eta\in \Phi_\Cbb^+(\slie_1,\hlie_2)$, $\delta(\eta^\vee)$ is a positive integer where $\eta^\vee$ is the coroot of $\eta$. Fix a positive simple root $\eta\in\Phi_\Cbb^+(\slie_1,\hlie_2)$. If $\delta(\eta^\vee)\ge 2$, then, using $\mathfrak{sl}_2$ triple, $\delta-\eta$ is again a non-zero weight. However, \[\|\delta-\eta\|^2-\|\delta\|^2=\|\eta\|^2(1-\delta(\eta^\vee))<0\] gives the contradiction to the minimality of $\|\delta\|$.

Let $m=\eta_1(H_\gamma)$ and $M=\eta_2(H_\gamma)$. 
Consider the finite Weyl group orbit $\Ocal=W\cdot\delta\subset \Lambda_1$ where $W$ is the Weyl group of $(\slie_1,\hlie_1)$. 
Then; \begin{enumerate}
\item for $\nu\in\Ocal$, $\nu=\omega\cdot\delta$ and $\nu(\eta^\vee)=\delta(\omega^{-1}\cdot\eta^\vee)$. Hence, $\nu(\eta^\vee)\in\{-1,0,1\}$ for all simple roots $\eta$.
\item 
$\sum_{\nu\in\Ocal}\nu=0$ as $\sum_{\nu\in\Ocal}\nu$ is $W$-invariant.
\end{enumerate}

We now prove that some $\lambda\in\Ocal$ satisfies the desired condition, i.e.,  \[-M\le\lambda(H_\gamma)\le-m.\]
For the sake of the contradiction, let us assume that there is no $\lambda\in\Ocal$ satisfies the inequality. Set
\[\Ocal_-=\{\nu\in\Ocal:\nu(H_\gamma)<-M\},\qquad\Ocal_+=\{\nu\in\Ocal:\nu(H_\gamma)>-m\}.\]
Then, by our assumptions, $\Ocal$ is a disjoint union of non-empty set $\Ocal_\pm$, $\Ocal=\Ocal_-\sqcup \Ocal_+$. Indeed,  $m<0<M$ implies $-M<0<-m$, so if $\Ocal_+$ was empty, then $\lambda(H_\gamma)<0$ for all $\lambda\in \Ocal$. But, it is not possible since $\sum_{\nu\in\Ocal}\nu=0$. Similarly, $\Ocal_-$ is not empty.

Choose $\nu_-\in\Ocal_-$ and $\nu_+\in\Ocal_+$. There is a sequence of reflections $s_1,\dots,s_l$  with respect to hyperplane perpendicular to simple roots $\alpha_1,\dots,\alpha_l$, respectively, such that $\nu_+=s_l\dots s_1 \nu_-$ since $W$ is generated by such $s_i$. Let $\nu_j=s_j\dots s_1 \nu_-$. \[\nu_j - \nu_{j-1}=s_{j-1}\nu_{j-1}-\nu_{j-1}=-\nu_{j-1}(\alpha_{j-1}^{\vee})\alpha_{j-1}.\] 
Since  $\nu(\eta^\vee)\in\{-1,0,1\}$ for all simple roots $\eta$, $\nu_j-\nu_{j-1}$ is one of $\alpha_{j-1},0,$ or $-\alpha_{j-1}$. Especially, $|\nu_j(H_\gamma)-\nu_{j-1}(H_\gamma)|\le M-m$ for all $j=1,\dots, l$ by \cref{eta<0}.

On the other hand, choose the smallest $j_0$ such that $\nu_{j_0}\in\Ocal_+$, i.e., \[\nu_{j_0}(H_\gamma)>-m.\] Such $1< j_0\ge l$ exists since $\nu_+=\nu_l \in\Ocal_+$ and $\nu_-=\nu_1\in\Ocal_-$. Then, $\nu_{j_0 - 1} \in \Ocal_-$, i.e., \[\nu_{j_0-1}(H_{\gamma})<-M.\] Hence, \[\nu_{j_0}(H_{\gamma})-\nu_{j_0 -1} (H_{\gamma})>-m+M\] which gives the contradiction with the previous bound $|\nu_j(H_\gamma)-\nu_{j-1}(H_\gamma)|\le M-m$ for all $j$ from \cref{eta<0}. Hence, there is $\lambda\in\Ocal$ such that \[-M\le \lambda(H_{\gamma})\le -m.\]
Finally, if $\lambda(H_\gamma)=-M$, choose a weight $\eta'\in\Lambda_2$ with $\eta'(H_\gamma)=M$. For nonzero weight vectors $v_\lambda\in(S_1)_\lambda$ and $w_{\eta'}\in(S_2)_{\eta'}$, 
\[T(v_\lambda\otimes w_{\eta'})=e^{\lambda(H_\gamma)+\eta'(H_\gamma)}v_\lambda\otimes w_{\eta'}=v_{\lambda}\otimes w_{\eta'}.\]
This contradicts to the fact that $T$ is hyperbolic. Similarly, $\lambda(H_\gamma)\neq -n$. Hence, $-M<\lambda(H_\gamma)<-m$.
 \end{proof}

The kernels of $p_j:\Acal_\Cbb \to \End^\Cbb(S_j)$, $\ker p_j=I_j$, are maximal two-sided ideals, as there is no proper two sided ideal of $\End_\Cbb(S_j)$ as a full matrix algebra. Also, if $i\neq j$, then $S_i$ and $S_j$ are inequivalent $\slie_1$ representation, hence, we have the following:
\begin{claim}
For $i,j\in \{1,\dots, q\}$, if $I_i=I_j$ then $i=j$.
\end{claim}
Indeed, if $I_i=I_j$, there is a isomorphism $\phi:\End^\Cbb(S_i)\to\End^\Cbb(S_j)$ with  $\phi\circ p_i=p_j$. The map $\phi$ provides a linear isomorphism $L$ between $S_i$ and $S_j$ that intertwines the representation.

 The kernels are thus distinct and pairwise comaximal. The Chinese remainder theorem gives
\[ \Acal_\C\simeq \prod_{j=1}^{q} \End^{\Cbb}(S_j)\simeq \Bcal_\Cbb.\]
Especially, $\dim_\Cbb \Acal_\Cbb=\dim_\Cbb\Bcal_\Cbb$. As $\Acal_\Cbb\subset \Bcal_\Cbb$, $\Acal_\Cbb=\Bcal_\Cbb$.
Consequently, $\dim_\Rbb\Acal=\dim_\Rbb\Bcal$, hence, 
\[\Acal=\Bcal=\bigoplus_i\End_{D_i}(U_i)\otimes I_{E_i}\] as $\Acal\subset \Bcal$.

\subsubsection{H\"older family of maps $\Phi_x:\Bcal \to \End_{\Rbb}(TM)$}

Let us retain the setting in \cref{thm:samealg}.  Let $P_{\gamma}=\rho(\gamma)P_0 \rho(\gamma)^{-1}$. Since $\Acal$ is finite dimensional algebra, we can find a finitely many $\gamma_1,\dots,\gamma_l \in \Gamma$ so that $\{I_V\}\cup \{P_{\gamma_j}\}_{j=1,\dots,l}$ generates $\Acal$ as a real algebra. Define \[P_j = P_{\gamma_j}\in \Bcal=\Acal,\qquad j=1,\dots,l.\] 
  Note that $\alpha(\gamma\gamma_0\gamma^{-1})$ are all Anosov diffeomorphisms for all $\gamma\in\Gamma$, hence, $\rho(\gamma\gamma_0\gamma^{-1})$ is hyperbolic. Therefore, $P_j$ is the projection into stable subspaces along unstable subspaces with respect to $\rho(\gamma_j\gamma_0\gamma_j^{-1})$.

 Recall that $TM=E^s_{0}\oplus E^u_{0}$ is the decomposition into H\"older stable and unstable distributions with respect to the Anosov diffeomorphism $\alpha(\gamma_0)$. Similarly, for each Anosov diffeomorphisms $\alpha(\gamma_j\gamma_0\gamma_j^{-1})$, we define the decomposition as $TM=E^s_j\oplus E^u_j$ for $j=1,\dots,l$. For $x\in M$, let \[Q_j(x)\in\End(T_xM)\] be the projection into $E^s_j (x)$ along $E^u_j$. Using the cocycle superrigidity and the fact that stable and unstable distributions are H\"older, we have the following: \begin{claim} Each $Q_j$ is everywhere defined H\"older endomorphism on $TM$, i.e., $Q_j$ is a H\"older section of $\End(TM)\to M$. Moreover, \[Q_j(x)=B_x P_j B_x^{-1}\] $\mu$-almost every $x\in M$.\end{claim}
 Since $\Bcal=\Acal$ and $\Acal$ is generated by $P_1,\dots,P_l$, evaluating polynomials on the $Q_j$, we have the following:
 \begin{claim}\label{claim:Phi_x}
 For each $T\in\Bcal$, let $F^T$ be a real finite degree polynomial with $l$-variables such that $F^T(P_1,\dots,P_l)=T$. Let us define $\Phi_x: \Bcal\to \End(T_x M)$ as \[\Phi_x(T)\coloneqq F^T(Q_1,\dots,Q_l)\] for all $x\in M$. Then, $\Phi_x$ is well-defined everywhere, 
 \begin{enumerate}
 \item for all $x\in M$, $\Phi_x$ is a faithful unital algebra homomorphism,
\item $\Phi_x$ is continuous with respect to $x$, and 
\item $\Phi_x(T)=B_x T B_x^{-1}$ for $\mu$-almost every $x\in M$.
\end{enumerate}
We denote the family of maps $\Phi_x$ by $\Phi:\Bcal\to \End(TM)$.
\end{claim}
\begin{proof}[Proof of \cref{claim:Phi_x}]
Note that every polynomial relations between $P_1,\dots, P_l$ holds between $Q_1(x),\dots, Q_l(x)$ for almost every $x$, hence, every $x$ using continuity of $Q_j$ and $\supp \mu=M$. This shows that $\Phi_x$ is well-defined everywhere and (3) holds.  Continuity comes from continuity of $Q_j$ and every algebra identities hold almost everywhere (through $B_x$ and $Q_j$), hence, hold everywhere. This shows well-definedness, items (1) -(3) except injectivity, i.e., $\Phi_*$ is the continuous family of an unital algebra homomorphism satisfies the item (3). 

It remains to prove that $\Phi_x$ is injective for every $x\in M$.
Write $\Bcal=\prod_i\Bcal_i$, where
$\Bcal_i\simeq\End_{D_i}(U_i)$ is a simple real algebra.
Let $e_i\in\Bcal$ be the identity in the $i$th factor and zero in
all other factors. As an endomorphism of $V$, $e_i$ is the
projection onto the $i$th irreducible component. Hence
$\tr_V(e_i)=\dim_{\Rbb}(e_iV)>0$.

By item~(3) and invariance of trace under conjugation,
\[
 \tr_{T_xM}\bigl(\Phi_x(e_i)\bigr)
 =\tr_{T_xM}\bigl(B_xe_iB_x^{-1}\bigr)
 =\tr_V(e_i) 
\] for almost every $x$.
The left-hand side is continuous in $x$ by item~(2).
Since $\operatorname{supp}\mu=M$, it equals $\tr_V(e_i)$
everywhere. In particular, $\Phi_x(e_i)\ne0$ for every $i$ and $x$. Fix $x\in M$ and suppose that $T\in\ker\Phi_x$ is nonzero.
Choose $i$ such that the $i$th component of $T$ is nonzero.
Then $e_iT$ is a nonzero element of
$\ker\Phi_x\cap e_i\Bcal$, since $\Phi_x$ is an algebra
homomorphism. This intersection is a two-sided ideal of the
simple algebra $e_i\Bcal\simeq\Bcal_i$. It must therefore equal
$e_i\Bcal$ and in particular contain $e_i$. This contradicts
$\Phi_x(e_i)\ne0$. Thus $\ker\Phi_x=\{0\}$ for every $x\in M$.
\end{proof}

Indeed, the family $\{\Phi_x\}$ of maps $\Phi_x$ is the family of faithful algebra embeddings whose elements are H\"older in smooth local frames of $TM$.

\begin{claim}\label{claim:Tcomm}
For all $T\in \Bcal\subset \End(V)$, $T$ commutes with $K \subset \GL(V)\subset \End(V)$.
\end{claim}
\begin{proof}[Proof of \cref{claim:Tcomm}]
Recall that $K$ commutes with every $\rho(g)$ for all $g\in G$. Note that the model stable projection $P_0$ is a real polynomial in the model matrix $\rho(\gamma_0)$ in $\End(V)$. Hence, $K$ also commutes with that projection $P_0$ as well as $P_j$. More directly, for any $k\in K$,  $k\rho(\gamma_0)=\rho(\gamma_0)k$ gives, for all $v_s\in V^s_0$, $k\rho(\gamma_0^n)(v_s)=\rho(\gamma_0^n)(kv_s) \to 0$ as $n\to \infty$. Hence, $KV^s_0\subset V^s_0$. Similarly, $KV^u_0\subset V^u_0$. This gives $P_0k(v_s+v_u)=P_0(kv_s+kv_u)=kv_s=kP_0(v_s+v_u)$ for all $v_s\in V^s_0$, $v_u\in V_0^u$, and $k\in K$. Therefore, $K$ commutes with $P_0$. Same argument works for $P_j$, too.  Since $\Acal=\Bcal$, $T$ commutes with $K$. 
\end{proof}
As a corollary, we also have the following:
 \begin{claim}\label{relationtoD} 
For every $\gamma\in\Gamma$ and $T\in \Bcal$, we have \begin{equation}\label{conjD}D_{x}\alpha(\gamma)\Phi_{x}(T) \left(D_{x}\alpha(\gamma)\right)^{-1} = \Phi_{\alpha(\gamma)(x)} (\rho(\gamma)T\rho(\gamma)^{-1})\end{equation} Here we consider $\rho(\gamma)\in \Bcal$ as $\bigoplus_i \rho_i(g)\otimes I_{L_{i}}$
 \end{claim}
 \begin{proof}[Proof of \cref{relationtoD}]
For $\mu$-almost every $x$, 
\begin{align*}
D_x\alpha(\gamma)\Phi_x(T)(D_x\alpha(\gamma))^{-1} &=B_{\alpha(\gamma)(x)}\rho(\gamma)\kappa(\gamma,x) T \kappa(\gamma,x)^{-1}\rho(\gamma)^{-1}B_{\alpha(\gamma)(x)}^{-1}\\
&=B_{\alpha(\gamma)(x)}\rho(\gamma)T\rho(\gamma)^{-1}B_{\alpha(\gamma)(x)}^{-1}\\
&=\Phi_{\alpha(\gamma)(x)}\!\left(\rho(\gamma)T\rho(\gamma)^{-1}\right).
\end{align*} as $\kappa(\gamma,x)\in K$.
Since both sides are all continuous, \cref{conjD} holds everywhere. \end{proof}

 \subsubsection{Elementary matrices} Note that there is a real algebra isomorphism \[\Bcal=\bigoplus_{i=1}^{m} \End^{D_{i}}(U_{i})\otimes_{D_{i}} I_{L_{i}} \simeq \bigoplus_{i=1}^{m}\End^{D_{i}}(U_{i}), \quad\bigoplus_{i}T_{i}\otimes I_{L_{i}}\mapsto \bigoplus_{i}T_{i}.\]
 Let $u_{i,1},\dots, u_{i,r_i}$ be a right $D_i$-basis of $U_i$. For $q\in D_i$, define $e_{i,ab}(q)\in \Bcal$ by \[e_{i,ab}(q)(u_{j,c}d)\coloneqq\begin{cases} 0 & i\neq j \\ \delta_{bc}u_{i,a}qd & i=j \end{cases}\] for all $d\in D_i$. Then, \[e_{i,ab}(q)e_{j,cd}(r)=\begin{cases} 0 & i\neq j \\ \delta_{bc}e_{i,ad}(qr) & i=j\end{cases}\] for all $i,j=1,\dots,m$, $a,b\in \{1,\dots,r_i\}$, $c,d\in\{1,\dots,r_j\}$, $q\in D_i$, and $r\in D_j$. Note that the real basis of $\Bcal$ is obtained by $\{e_{i,ab}(q)\}$ allowing $q$ to range over $\{1\}$, $\{1,i\}$ or $\{1,i,j,k\}$.  $e_{i,ab}(q)$ will play the role of elementary matrices.

\subsubsection{Reconstruct the tangent bundle}

Let \[E_{i,ab}(q)\coloneqq \Phi(e_{i,ab}(q))\in \End(TM) ,\qquad E_{i,ab}=E_{i,ab}(1)\in \End(TM)\] for all $i=1,\dots,m$ and $a,b\in\{1,\dots,r_i\}$. Note that $E_{i,ab}(q)$ and $E_{i,ab}$ are all continuous sections of $\End(TM)\to M$. Consider $B_x$ first, we know desired algebra identities hold almost every $x$. Hence, by continuity, all algebraic identities hold among $e_{i,ab}(q)$ hold for all $x\in M$. For instance, \[E_{i,ab}\circ E_{j,cd}=\delta_{ij}\delta_{bc}E_{i,ad},\quad \sum_{i,a}E_{i,aa}=I_{TM}, \quad \tr E_{i,aa}=(\dim_{\Rbb}D_i) d_i=\dim_{\Rbb}L_i\] holds almost everywhere (through $B_x$), hence, hold everywhere. Set the continuous real vector bundle \[E_i =\im E_{i,11} \subset TM,\]  that is, $E_{i}(x)=\im E_{i,11}(x)$ for $x\in M$. Since $E_{i,11}(x)$ is an idempotent, $E_{i,11}(x)^2 =E_{i,11}(x)$, the rank of the vector bundle $E_{i}$ is $\tr E_{i,11}=(\dim_{\Rbb}D_{i})d_{i}=\dim_{\Rbb}L_i$.

Let us define left $D_i$-action on $E_i(x)$ for each $x\in M$. For $q\in D_i$ and $v\in E_{i}(x)$, we define \[q\cdot v\coloneqq E_{i,11}(q)(x)(v).\] Note that \[E_{i,11}\circ E_{i,11}(q)=E_{i,11}(1)\circ E_{i,11}(q)=E_{i,11}(q)\] holds everywhere, as $e_{i,11}(q)=e_{i,11}(1)e_{i,11}(q)$, hence the map is well-defined. One can easily verify that $(qr)\cdot v=q\cdot (r\cdot v)$ and $1\cdot v=v$ for all $q,r\in D_i$ and $v\in E_{i}(x)$. Hence, it gives left $D_i$-action.

\begin{remark}
When $D_i=\Cbb$, this is same with having the H\"older section $J$ of $\End(TM)\to M$ so that $J^2=-I_{TM}$. For $D_i=\Hbb$ case, this gives H\"older sections $J_1,J_2,J_3$ such that $J_{1}J_{2}=J_{3}=-J_{2}J_{1}$ and $J_{1}^2=J_2^2=J_3^2=-I_{TM}$.
\end{remark}

Define the isomorphism which one can consider it as the reconstruction of the tangent bundle as
\begin{equation}\Psi_x :\bigoplus_{i}U_i \otimes_{D_{i}} E_{i}(x) \to T_{x} M , \qquad \Psi_{x}(u_{i,a}\otimes v)\coloneqq E_{i,a1}(x)(v)\end{equation} for all $x\in M$. Recall that every element $w\in \bigoplus_{i} U_{i}\otimes_{D_{i}}E_{i}$ can be written uniquely $w=\sum_{i,a} u_{i,a}\otimes v_{i,a} $, so, $\Psi_x$ is well-defined. The inverse map can be written as \[X \mapsto \sum_{i,a}u_{i,a}\otimes\left( E_{i,1a}(x)(X)\right)\] where $X\in T_x M$. Note that, as $E_{i,1a}=E_{i,11}\circ E_{i,1a}$, $E_{i,1a}(x)(X)\in E_{i}(x)$. Note that $E_{i}(x)$ is H\"older vector bundle, $\Psi_x$ is H\"older continuous on $x$.

\subsection{Recovering multiplicity space $L_i$}\label{sec3.1}
We retain the setting and notations in the previous subsections.
\subsubsection{Identify $\kappa$ in $\Bcal$}
Since $\rho(G)$ commutes with $K$, we have the following: 
\begin{claim}\label{kforminV}
For the element $k\in K$, $k$ has the form $\bigoplus_{i} I_{U_{i}}\otimes k_i$ in $\End(V)$. Denote the compact group $\{k_i: k\in K\}$ by $K\restrict{L_i}$.
\end{claim}
\begin{proof}[Proof of \cref{kforminV}]
Recall that we can canonically identify the multiplicity space as $L_i=\Lcal_G(U_i,V).$ Thus an element $e\in L_i$ is a real-linear map $e:U_i\to V$ satisfying \[e(\rho_i(g)u)=\rho(g)e(u).\]
Recall that its left $D_i$-action is $(d\cdot e)(u)=e(ud).$
For $k\in K$, define
\[k_i:L_i\to L_i ,\qquad k_i(e)=k\circ e.\]
\begin{claim}\label{3.12} $k\circ e$ is $G$-equivariant, hence, $k_i(e)\in L_i$. Furthermore, $k_i$ is left $D_i$-linear and invertible with the inverse $e\mapsto k^{-1}\circ e$. \end{claim}
\begin{proof}[Proof of \cref{3.12}] $k\circ e$ is $G$-equivariant since, for $u\in U_i$ and $g\in G$,
\[(k\circ e)(\rho_i(g)u)=k\left(\rho(g)e(u)\right)=\rho(g)k(e(u))=\rho(g)(k\circ e)(u).\]
Hence $k_i(e)$ does belong to $L_i$. Also, $k_i$ is left $D_i$-linear. For $d\in D_i$, evaluating at any $u\in U_i$ gives
\[(k_i(d\cdot e))(u)=k\left((d\cdot e)(u)\right)=k(e(ud))=(k_i e)(ud)=(d\cdot(k_i e))(u).\]
Therefore, $k_i(d\cdot e)=d\cdot k_i(e)$. Its inverse is $e\mapsto k^{-1}\circ e$, so $k_i$ is invertible.\end{proof}

Furthermore, recall that the decomposition is given by 
\[\textrm{ev}:\bigoplus_i U_i\otimes_{D_i}L_i\to  V, \qquad \textrm{ev}(u\otimes e)=e(u).\]
Hence, for $u\in U_i$ and $e\in L_i$, \[k\textrm{ev}(u\otimes e)= k(e(u))= (k_i e)(u) = \textrm{ev}(u\otimes k_i e).\]
Thus, under this identification,
\[k=\bigoplus_i I_{U_i}\otimes_{D_i}k_i.\]
The formula $k_i(e)=k\circ e$ also proves uniqueness. \end{proof}
\subsubsection{Recover multiplicity space}

Set \[C_{\gamma}(x)=\Phi_{\alpha(\gamma)(x)}(\rho(\gamma)^{-1})D_{x}\alpha(\gamma) : T_{x}M \to T_{\alpha(\gamma)(x)}M.\] As $\Phi$ and $D\alpha(\gamma)$ are H\"older, $C_{\gamma}(x)$ is H\"older on $x$. 
Using \cref{conjD}, we can deduce \begin{equation} C_{\gamma}(x)\Phi_{x}(T)=\Phi_{\alpha(\gamma)(x)}(T)C_{\gamma}(x)\end{equation} for all $T\in \Bcal$.
Especially, when $T=e_{i,11}(1)$,  \[C_{\gamma}(x)E_{i,11}(x) =E_{i,11}(\alpha(\gamma)(x))C_{\gamma}(x).\] This implies that \[C_{\gamma}(x)\restrict{E_{i}(x)} :E_{i}(x)\to E_{i}(\alpha(\gamma)(x))\] is well-defined. Let us define \[\Acal_i(\gamma,x) :E_{i}(x)\to E_{i}(\alpha(\gamma)(x)),\qquad \Acal_i(\gamma,x)\coloneqq C_{\gamma}(x)\restrict{E_{i}(x)}.\] 
$\Acal_i$ also preserves left $D_i$ actions; 
\begin{align*}
\Acal_i(\gamma,x)(q\cdot v) &=C_\gamma(x) E_{i,11}(q)(x)v\\
&= E_{i,11}(q)(\alpha(\gamma)(x))C_\gamma(x)v\\ &=q\cdot \Acal_i(\gamma,x)v.
\end{align*}

Combining above discussions, we have the following:
\begin{claim}\label{Holdermult} We have the following equation;  \begin{equation}\label{Psiconj} \Psi_{\alpha(\gamma)(x)}^{-1}D_{x}\alpha(\gamma)\Psi_x =\bigoplus_i \rho_{i}(\gamma)\otimes A_{i}(\gamma,x)\end{equation} for all $\gamma\in\Gamma$ and for all $x\in M$. Especially, $\Acal_{i}$ is H\"older cocycle over the H\"older bundler $E_{i} \to M$. \end{claim}
\begin{proof}[Proof of \cref{Holdermult}] 
For a fixed $u_{i,a}$, let \[\rho_i(\gamma)u_{i,a} =\sum_b u_{i,b}c^\gamma_{ba}.\]
This implies that, as $\rho(\gamma)=\bigoplus_i \rho_i(g) \otimes I_{L_{i}}$,
\[\rho(\gamma)e_{i,a1}(1)=\sum_b e_{i,b1}(c^\gamma_{ba}).\] Indeed, evaluating either side on $u_{i,c}d$, with $d\in D_i$, gives
\[ \delta_{1c}\sum_b u_{i,b}(c^\gamma_{ba}d). \] Therefore, 
\begin{align*}
\Phi_{\alpha(\gamma)(x)}(\rho(\gamma)) E_{i,a1}(\alpha(\gamma)(x)) &=\sum_b E_{i,b1}(c^\gamma_{ba})(\alpha(\gamma)(x))\\
&=\sum_b  E_{i,b1}(1)(\alpha(\gamma)(x)) \circ E_{i,11}(c^\gamma_{ba})(\alpha(\gamma)(x)).
\end{align*}
Now $\Acal_i(\gamma,x)v\in E_i(\alpha(\gamma)(x))$, and the left scalar action on this space is, by definition, \[c\cdot \Acal_i(\gamma,x)= E_{i,11}(c)(\alpha(\gamma)(x))\Acal_i (\gamma,x).\]
Then,
\begin{align*} D_{x}\alpha(\gamma)\Psi_x (u_{i,a}\otimes v) & =D_x \alpha(\gamma) E_{i,a1}(x)v =D_x \alpha(\gamma)\Phi_x(e_{i,a1}(1))v\\ 
&= \Phi_{\alpha(\gamma)(x)}(\rho(\gamma)e_{i,a1}(1)\rho(\gamma)^{-1})D_{x}\alpha(\gamma)v\\
&=\Phi_{\alpha(\gamma)(x)}(\rho(\gamma))E_{i,a1} (\alpha(\gamma)(x))C_{\gamma}(x)v\\
&=\Phi_{\alpha(\gamma)(x)}(\rho(\gamma))E_{i,a1} (\alpha(\gamma)(x))\Acal_i(\gamma,x)(v)\\
&=\sum_b  E_{i,b1}(1)(\alpha(\gamma)(x))  \left(c^\gamma_{ba}\cdot \Acal_i(\gamma,x)v\right)\\
&=\Psi_{\alpha(\gamma)(x)}\left(\sum_bu_{i,b}\otimes(c^\gamma_{ba}\cdot \Acal_i(\gamma,x)v) \right)\\
&=\Psi_{\alpha(\gamma)(x)}\left(\sum_b(u_{i,b}c^\gamma_{ba})\otimes\Acal_i(\gamma,x)v\right)\\
&=\Psi_{\alpha(\gamma)(x)}\left(\rho_i(\gamma)u_{i,a}\otimes \Acal_i(\gamma,x)v\right).
\end{align*}
 This shows the \eqref{Psiconj}. Since \[\Psi_{\alpha(\gamma'\gamma)(x)}^{-1}D_{x}\alpha(\gamma')\Psi_{\alpha(\gamma)(x)}  \Psi_{\alpha(\gamma)(x)}^{-1}D_{x}\alpha(\gamma)\Psi_x   =\Psi_{\alpha(\gamma'\gamma)(x)}^{-1}D_{x}\alpha(\gamma'\gamma)\Psi_x, \] this gives the desired cocycle equation \[\Acal(\gamma'\gamma,x)=\Acal(\gamma',\alpha(\gamma)(x))\Acal(\gamma,x)\] for all $x\in M$ and $\gamma,\gamma'\in\Gamma$.
\end{proof}

Set the measurable family of maps \[J_{i}(x):L_i \to E_{i}(x), \qquad J_{i}(x)(v)\coloneqq B_{x}(u_{i,1}\otimes v)\] for all $i=1,\dots,m$. Since, using \cref{kforminV}, the compact error $\kappa(\gamma,x)$ is of the form $\bigoplus_i I_{U_{i}} \otimes \kappa_i (\gamma,x)$, one can deduce \begin{equation}\label{msbleAcal} \Acal_i (\gamma,x)J_{i}(x)=J_{i}(\alpha(\gamma)(x))\kappa_i(\gamma,x)\end{equation}

Fix a positive real inner product $h_i$ on $L_i$ that is invariant under the compact group $K\restrict{L_i}$ and under unit scalars in $D_i$, that is, $\{\pm 1\}$, $U(1)$, or the unit quaternions. We can find it by averaging since $K\restrict{L_i}$ and the unit scalars in $D_i$ are compact and commute, i.e., $k_i(q\cdot e)=q\cdot k_i (e)$ for $q\in D_i$ and $e\in L_i=\Lcal_{G}(U_i,V)$. Let us define the measurable family of metric \[g_i (x) (J_i(x)v,J_i(x)v)=h_i(v,w)\] for all $v,w\in L_i$ and $\mu$-almost every $x$. Then \[g_i (\alpha(\gamma)(x))(\Acal_i (\gamma,x)(v),\Acal_i (\gamma,x)(w))=g_i(x)(v,w)\] for all $v,w\in E_i(x)$,$\gamma\in \Gamma$ and $\mu$-almost every $x$. Using \cref{prop:conf}, there is a H\"older representative of $g_i$ such that $\Acal_i$ preserves it. In summary, 
\begin{claim}\label{ConfAcal}\label{thm:multiplicityspace}
There is a H\"older family of metric (positive inner product) $g_i(x)$ on $E_i(x)$ such that \[ g_{i}(\alpha(\gamma)(x))(\Acal(\gamma,x)(v),\Acal(\gamma,x)(w))=g_{i}(x)(v,w)\] for all $v,w\in E_i(x)$, $x\in M$, $\gamma\in \Gamma$, and $i=1,\dots, m$.
\end{claim}

\subsection{Continuous compact group reduction}

Put $m_i=\dim_{D_i}L_i$. On a sufficiently small neighborhood $\Ocal_x$ of $x\in M$, choose a H\"older local $D_i$-frame $w_1,\ldots,w_{m_i}$ of $E_i$. Our choices of the metrics $h_i$, hence $g_i$, are compatible with multiplication by unit elements of $D_i$. We can apply Gram--Schmidt to entire $D_i$-spans to get a  H\"older local frame on $\Ocal_x$. 

More explicitly, set
$S_{\Rbb}=\{1\}, S_{\Cbb}=\{1,\mathrm i\}, S_{\Hbb}=\{1,\mathrm i,\mathrm j,\mathrm k\},$
and define inductively
\[ v_j=w_j-\sum_{l<j}\sum_{a\in S_{D_i}} g_i(w_j,a\cdot q_l)\,a\cdot q_l, \qquad q_j=\frac{v_j}{\sqrt{g_i(v_j,v_j)}}.\]
The scalar compatibility of $g_i$ implies that the vectors $a\cdot q_j$, for $a\in S_{D_i}$ and $1\le j\le m_i$, are real orthonormal. The denominators are positive and, after shrinking the neighborhood, bounded away from zero. Hence we may assume that the resulting fields $q_j$ are H\"older on the neighborhood $\Ocal_x$.

Mapping a fixed $D_i$-orthonormal basis of $(L_i,h_i)$ to the $q_j(x)$ gives a H\"older family of $D_i$-linear isometries \[ \ell_i(x):(L_i,h_i)\to(E_i(x),g_i(x)). \] Let 
\[P_x^0= \left\{ \Psi_x\circ\bigoplus_i(I_{U_i}\otimes_{D_i}\ell_i) \middle| 
\ell_i:(L_i,h_i)\to(E_i(x),g_i(x)) \text{ is a $D_i$-linear isometry.}\right\}.\]
Let $\textrm{O}^{D_i}(L_i,h_i)$ denote the group of left $D_i$-linear isometries of the real inner product $h_i$, and set
\[C_0=\left\{\bigoplus_i I_{U_i}\otimes_{D_i}c_i \middle| c_i\in\textrm{O}^{D_i}(L_i,h_i) \right\} \subset\textrm{GL}(V).\]
Each factor is a closed subgroup of a real orthogonal group. Thus $C_0$ is compact. Moreover, $C_0$ commutes with $\rho(G)$, because it acts on the multiplicity factors.

\begin{claim}
Every fiber $P_x^0$ is nonempty and is a single free right $C_0$-orbit. $P^0\subset Fr(TM)$ is the set of H\"older local sections.
\end{claim}

\begin{proof}
The local isometries $b_i$ constructed above showed that every fiber $P_x^0$ is non-empty. If $p,p'\in P_x^0$ correspond to isometries $\ell_i,\ell_i'$, then
\[p^{-1}p'=\bigoplus_i I_{U_i}\otimes_{D_i}(\ell_i^{-1}\ell_i') \in C_0. \] Conversely, $pc\in P_x^0$ whenever $p\in P_x^0$ and $c\in C_0$.
The right action is free because every $p$ is invertible.
\end{proof}

\begin{claim}
For every $x\in M$ and $p\in P_x^0$, one has $|\det p|=1$.
\end{claim}

\begin{proof}
Since $C_0$ is compact, every $c\in C_0$ satisfies $|\det c|=1$.
Consequently \[\nu(x)=|\det p|,\qquad p\in P_x^0,\] is independent of the choice of $p$. Local H\"olderness of $P^0$ shows that $\nu$ is a positive H\"older function. Recall that the measurable maps $J_i(x):L_i\to E_i(x)$ are $D_i$-linear isometries for $h_i$ and $g_i(x)$ almost everywhere,
and that \[B_x=\Psi_x\circ\bigoplus_i(I_{U_i}\otimes_{D_i}J_i(x)). \] Thus $B_x\in P_x^0$ almost everywhere. 
$|\det B_x|=1$ implies $\nu(x)=1$ almost everywhere. Continuity and $\textrm{supp}\mu=M$ give $\nu=1$ everywhere.
\end{proof}

Define $P_x=P^0_x,C=C_0$.
Finally, we can show that $P$ is the desired continuous reduction: 
\begin{proposition}[Equivariant H\"older reduction $P$]\label{linearization}
For every $x\in M$ and $\gamma\in\Gamma$,
\[D_x\alpha(\gamma)P_x =P_{\alpha(\gamma)(x)}\rho(\gamma). \]
For every $b_x\in P_x$ and $b_{\alpha(\gamma)(x)}\in P_{\alpha(\gamma)(x)}$, there is a unique $k\in C$ such that \[ D_x\alpha(\gamma)b_x =b_{\alpha(\gamma)(x)}\rho(\gamma)k. \]
\end{proposition}

\begin{proof}
Write \[ p_x=\Psi_x\circ\bigoplus_i(I_{U_i}\otimes_{D_i}\ell_i)\in P_x. \] Each composition $ \Acal_i(\gamma,x)b_i$ is a $D_i$-linear isometry from $(L_i,h_i)$ to $(E_i(\alpha(\gamma)(x)),g_i(\alpha(\gamma)(x)))$. Hence \[ \wtd p_{\alpha(\gamma)(x)} \coloneqq\Psi_{\alpha(\gamma)(x)}\circ\bigoplus_i  \left(I_{U_i}\otimes_{D_i} \Acal_i(\gamma,x)\ell_i\right) \in P_{\alpha(\gamma)(x)}^0.\]
\eqref{Psiconj} gives
\begin{align*}
D_x\alpha(\gamma)p_x &=\Psi_{\alpha(\gamma)(x)}  \left(\bigoplus_i\rho_i(\gamma)\otimes_{D_i}\Acal_i(\gamma,x)\right)\left(\bigoplus_i I_{U_i}\otimes_{D_i}\ell_i\right)\\
&=\widetilde p_{\alpha(\gamma)(x)}\rho(\gamma).
\end{align*}
Since the action preserves the chosen positive density and $|\det p_x|=|\det\rho(\gamma)|=1$, we have $|\det\widetilde p_{\alpha(\gamma)(x)}|=1$. Thus $\widetilde p_{\alpha(\gamma)(x)}\in P_{\alpha(\gamma)(x)}$, proving one inclusion. Applying the same argument to $\gamma^{-1}$ proves the other inclusion, hence proves the equality.

Write $\widetilde p_{\alpha(\gamma)(x)}=b_{\alpha(\gamma)(x)} c$ for some $c\in C$. Since $C$ commutes with $\rho(G)$,
\[ D_x\alpha(\gamma)b_x =b_{\alpha(\gamma)(x)} c\rho(\gamma) =b_{\alpha(\gamma)(x)}\rho(\gamma)c. \]
The desired $k$ is therefore $k=c$: \[k=\rho(\gamma)^{-1}(b_{\alpha(\gamma)(x)})^{-1}D_x\alpha(\gamma)b_x \] and it is unique from the formula.
\end{proof}
 Combining \cref{prop:nozero,linearization}, we can deduce \cref{thm:reduction}.

\subsection{The adapted metric $h$ and Lyapunov splitting}\label{sec:adaptmetric}
Recall that we fix a maximal real split torus $A$, a compact torus $M$ commutes with $A$, and a lattice Cartan $\Zbb^r\simeq\Delta<\Gamma$ in \cref{sec:cartanlattice}. 

Let $\Sigma_{V}\subset \alie^*$ be the finite set of restricted weights of $\rho$. Write  \[V=\bigoplus_{\lambda\in \Sigma_{V}} V_{\lambda}, \quad \rho(\exp (X))\restrict{V_{\lambda}}=e^{\lambda(X)}I_{V_{\lambda}}\] for all $X\in \alie$. Recall that, in general, for a representation $(W,\rho_W)$ of $G$, we denote the space of restricted weights of $\rho_W$ by $\Sigma_{\rho_W}$.

Since both compact groups $C$ and $\rho(M)$ preserves every $V_\lambda$ and commutes, we can find (after averaging) an inner product $\langle\cdot, \cdot\rangle_V$ (hence the norm $||\cdot||_V$) on $V$ that is invariant under $C$ and $\rho(M)$ simultaneously. For a local section $b\in P$, define a metric $h$ on $TM$ by \[h_x (u,v)=\langle b^{-1}_x u,b^{-1}_x v\rangle_V\] for all $x\in M$ and for all $u,v\in T_xM$. A different choice of $b'\in P$ can be written as $b'=bc$ for some $c\in C$ so that $h$ does not depend on the choice of $b$. Let $E^{\lambda}_{x}=b_{x}(V_{\lambda})$. Since $C$ preserves $V_\lambda$, $E^{\lambda}_{x}$ also does not depend on the choice of $b$.

Since $P$ is H\"older reduction, we can choose H\"older local section, hence, $h$ is H\"older metric and $E^\lambda$ is H\"older bundle. With respect to this metric $h$, we have the followings:

\begin{proposition}[H\"older metric and H\"older Lyapunov splitting]\label{prop:metric}
There are a H\"older metric $h$ and a $D\alpha(\Delta)$-invariant H\"older splitting \[TM= \bigoplus_{i}\bigoplus_{\lambda\in\Sigma_V} E^{\lambda}\]  such that, for all $x\in M$, $v\in E^{\lambda}_x$, and $\delta\in \Delta$,  \[||D_{x}\alpha(\delta)v||_{h}=e^{\lambda(\log a_{\delta})}||v||_{h}.\] Especially, $\alpha(\delta)$ is Anosov diffeomorphism if and only if $\lambda (\log a_\delta)\neq 0$ for all $\lambda\in \Sigma_V$. 

 In addition, for any finite collection $\Ccal \subset \alie^*$,  there is $\delta \in\Delta$ such that $\xi (\log a_\delta) \neq 0$ for all $\xi \in \Ccal$. Hence, there is $\delta_0\in \Delta$ such that $\alpha(\delta_0)$ is Anosov diffeomorphism.
\end{proposition}

\begin{remark} This gives H\"older Lyapunov splitting of $\alpha(\Delta)$-action. \end{remark}
\begin{remark} One can show that, for general $\gamma\in\Gamma$, $\rho(\gamma)$ is hyperbolic if and only if $\alpha(\gamma)$ is Anosov diffeomorphism. \end{remark}

\section[H\"older torsion-free invariant connection]{H\"older torsion-free $\Gamma$-invariant connection $\nabla$}
\label{sec4}

In this section, we construct H\"olcer torsion free-invariant connection, \cref{thm:invconn}. Zimmer's cocycle superrigidity supplies a measurable invariant torsion-free connection. Using Anosov diffeomorphisms, we regularize the same connection. 

\begin{theorem}[A H\"older invariant torsion-free connection]
\label{thm:invconn}
In the \cref{setting}, In the standing $C^{2+\eta}$ lattice-action setting, there is a H\"older torsion-free $\Gamma$-invariant affine connection $\nabla$.
\end{theorem}
Using the application of cocycle superrigidity \cref{msblecon}, we first find a a measurable torsion-free affine connection $\nabla^m$ on $M$ that is $\Gamma$-invariant. After that we improve the regularity of $\nabla^m$ to $\nabla$ using a carefully chosen Anosov diffeomorphism. Especially, for $\mu$ almost everywhere, \cref{thm:invconn} will be equal to $\nabla^m$.

\subsection{Improving regularity of measurable sections}

We retain the notations in \cref{sec:Anosov}. Let $E$ be a H\"older finite rank real vector bundle over $M$. An invertible affine cocycle over $f$ is generated by affine bijections
\[ F_x:E_x\longrightarrow E_{f(x)},\qquad F_x(v)=L_xv+a_x, \]
where $L_x:E_x\to E_{f(x)}$ is a linear isomorphism and $a_x\in E_{f(x)}$.

Note that $a$ is a H\"older section of the pullback bundle $f^*E$, not a section of $E$.
For $n\ge1$ set \[L_x^{(n)}=L_{f^{n-1}x}\cdots L_x:E_x\to E_{f^n x},\qquad L_x^{(-n)}=(L_{f^{-n}x}^{(n)})^{-1}:E_x\to E_{f^{-n}x}, \] and put $L_x^{(0)}=I_{E_x}$.
The following will improve the regularity of $\nabla^m$ to H\"older later. 
\begin{proposition}\label{prop:conn} Let $M$ be a closed smooth manifold and let $f:M\to M$ be a volume preserving $C^{1+\eta}$ Anosov diffeomorphism, where $\eta>0$. Let $\mu$ be the normalized probability measure induced by the volume. Let $E\to M$ be a finite-dimensional H\"older real vector bundle, and let $F_x(v)=L_xv+a_x:E_x\to E_{f(x)}$ be an invertible affine
H\"older cocycle. Assume the followings:
\begin{enumerate}
\item $E=E^+\oplus E^0\oplus E^-$ is a H\"older splitting and $L_xE_x^\sigma=E_{f(x)}^\sigma$ for $\sigma\in\{+,0,-\}$.
\item There is a H\"older metric on $E$, some $C\ge1$ and $0<\lambda<1$, and every $x\in M$, $n\ge1$, \[ \|L_x^{(n)}|_{E_x^-}\|\le C\lambda^n, \qquad \|L_x^{(-n)}|_{E_x^+}\|\le C\lambda^n.\]
\item There is a measurable positive definite metric $g^0$ on $E^0$ such that, for $\mu$ almost every $x$ and all $u,v\in E_x^0$, \[ g^0_{f(x)}(L_xu,L_xv)=g^0_x(u,v). \]
\end{enumerate}
If a $F$-invariant measurable section $s$, i.e., $ s(f(x))=L_xs(x)+a_x$, for $\mu$-almost everywhere, then it agrees $\mu$-almost everywhere with a $F$-invariant H\"older section $\widetilde s$, i.e.,   \[ \widetilde s(f(x))=L_x\widetilde s(x)+a_x\] for every $x\in M$.  \end{proposition}

\begin{proof}[Proof of \cref{prop:conn}] The idea is that we can regularize $s$ on each $E^\sigma$, $\sigma\in\{+,0,-\}$, separately. On $E^{\pm}$, we can use the uniform contraction (or expansion) dynamics. On $E^0$, we can use \cref{prop:conf,prop:sectionconf}. Let $\pi_x^\sigma:E_x\to E_x^\sigma$ be the projections along the splitting. Define, for each $\sigma\in\{+,0,-\}$, 
\[ s^\sigma(x)=\pi_x^\sigma s(x),\qquad a_x^\sigma=\pi_{f(x)}^\sigma a_x. \]
The projections are H\"older and $\pi_{f(x)}^\sigma L_x=L_x\pi_x^\sigma$. Therefore, each component is $F$-invariant. 

\emph{On $E^{\pm}$:} As $a$ is a H\"older section of $f^*E$, there is a constant $C_a$ such that $||a_x||<C_a$ for all $x\in M$. Uniform contraction gives a uniformly convergent series in $E_x^-$:
\begin{equation}\label{eq:affine-stable-series} t^-(x)=\sum_{j=0}^{\infty} L_{f^{-j}x}^{(j)}a^-_{f^{-(j+1)}x}\in E_x^-.\end{equation}
For $j$, $a^-_{f^{-(j+1)}x}$ lies in $E^-_{f^{-j}x}$, hence $L_{f^{-j}x}^{(j)}a^-_{f^{-(j+1)}x}\in E_x^-$. Similarly, uniform expansion (backward contraction) on $E^+$ gives
\begin{equation}\label{eq:affine-unstable-series} t^+(x)=-\sum_{j=0}^{\infty}(L_x^{(j+1)})^{-1}a^+_{f^j x}\in E_x^+.\end{equation}
The summands in \eqref{eq:affine-stable-series} are bounded by $C\lambda^j\max_{x\in M}\{\|a^-_x\|\}<C\lambda^j C_a$, and the summand in \eqref{eq:affine-unstable-series} is bounded by
$C\lambda^{j+1}\max_{x\in M}\{\|a^+_x\|\}<C\lambda^{j+1}C_a$.
Thus, as $0<\lambda<1$, both sums uniformly converge. Hence, $t^\sigma$ is continuous and bounded. Furthermore, we have \[ t^\sigma(f(x))=L_xt^\sigma(x)+a_x^\sigma,\qquad \sigma\in\{+,-\}\]
The fact that $t^\sigma$ is H\"older (possibly, smaller exponent than $\eta$) is now followed by a standard argument. See, for instance, \cite{Katok--Hass}.

Now, put $r^-=s^--t^-$. Then $r^-(f(x))=L_xr^-(x)$ almost everywhere. Iteration gives \[ \|r^-(x)\|\le C\lambda^n\|r^-(f^{-n}x)\|\] on $\mu$ almost everywhere. Consequently, for every $\delta>0$, measure invariance implies, as $n\to \infty$
\[
 \mu\{\|r^-\|>\delta\}
 \le\mu\{\|r^-\|>\delta/(C\lambda^n)\}
 \longrightarrow0.
\]
The limit follows because $r^-$ is finite almost everywhere. Hence $s^-=t^-$, for $\mu$ almost everywehere. Similarly, $s^+=t^+$, for $\mu$ almost everywhere.

\emph{On $E^0$:} After we embed $\Aff(E^0)$ into $\GL(E^0\oplus \Rbb)$, we can consider a linear cocycle. More precisely, let us consider trivial $\Rbb$-bundle $M\times \Rbb$ on $M$. On $D=E^0\oplus(M\times \Rbb)$, consider the H\"older linear cocycle \[A_x(u,t)=((L\restrict{E^0})_xu+t a_x^0,t).\] The given measurable solution $s^0$ and metric $g^0$ define \[H_x((u,t),(v,r))=g_x^0(u-t s^0(x),v-r s^0(x))+tr.\] As its value on $((u,t),(u,t))$ is $g_x^0(u-t s^0(x),u-t s^0(x))+t^2$, $H$ is positive definite measurable metric on $D$. 
The identity \[(L\restrict{E^0})_xu+t a_x^0-t s^0(f(x))=(L\restrict{E^0})_x(u-t s^0(x))\] and invariance of $g^0$ imply \[H_{f(x)}(A_x(u,t),A_x(v,r))=H_x((u,t),(v,r)) \] $\mu$-almost everywhere.

Using \Cref{prop:conf}, for the cocycle $A$ and a measurable invariant metric $H$, we can get a H\"older metric $\what H$ that equals $H$ $\mu$-almost everywhere.   Define the H\"older positive definite map $J_x^0:E_x^0\to(E_x^0)^*$ and H\"older $\ell_x\in(E_x^0)^*$ by \[ (J_x^0u)(v)=\widehat H_x((u,0),(v,0)),\qquad \ell_x(u)=\widehat H_x((u,0),(0,1)). \]
Then, $\mu$-almost everywhere, $J_x^0=g_x^0, \ell_x=-g_x^0(\,\cdot\,,s^0(x))$. Therefore,  \[t^0(x):=-(J_x^0)^{-1}\ell_x=s^0(x) \quad\text{$\mu$-almost everywhere}.\]
Furthermore, $t^0$ is H\"older continuous.

Finally set $\widetilde s=t^-+t^0+t^+$. This section is H\"older and equals $s$ almost everywhere. Its affine invariance equation holds almost everywhere (through $s$) and both sides are continuous sections of $f^*E$. Since $\mu$ has full support, the equation holds everywhere.
\end{proof}

\subsection{The affine action on torsion-free connections}

Fix a smooth torsion-free connection $\nabla^0$ on $TM$, for example the Levi-Civita connection of a smooth Riemannian metric. Define the smooth bundle $\Ecal$ by
\[\Ecal=\Hom(\textrm{Sym}^2TM,TM), \qquad W=\Hom(\textrm{Sym}^2V,V).\] We call $\Ecal$ a \emph{connection bundle}.
An element of $\mathcal E_x$ is a symmetric bilinear map $T_xM\times T_xM\to T_xM$. The difference of two torsion-free connections is a section of $\mathcal E$. Conversely, a section $S$ of $\Ecal$ defines the torsion-free connection $\nabla^0+S$ by \[(\nabla^0+S)_X Y=\nabla^0_XY+S(X,Y),\] where $X,Y$ are vector fields on $M$. 
The regularity (for instance, measurable or H\"older) of the section $S$ will give the same regularity of $\nabla^0+S$.

Recall that, for a $C^1$ diffeomorphism $f:M\to M$, the pushforward connection is defined by \[ (f_*\nabla)_X Y =f_*(\nabla_{f_*^{-1}X}(f_*^{-1}Y)).\]
Here $X,Y$ are vector fields on the range of $f$, and $f_*^{-1}X$ denotes their pullback vector field on the domain. Furthermore, for a section $S$ of $\Ecal$, we have the identity \[f_*(\nabla^0+S)=f_*\nabla^0+f_*S=\nabla^0+(f_*\nabla^0-\nabla^0)+f_*S\] where \[(f_*S)_{f(x)}(u,v)= D_xf(S_x((D_xf)^{-1}u,(D_xf)^{-1}v))\] for $x\in M$ and $u,v\in T_{f(x)}M$.

Direct calculation with respect to a smooth chart coordinate shows the following:
\begin{lemma}\label{lem:conaff}
Let $f:M\to M$ be a $C^{2+\eta}$ diffeomorphism. For $x\in M$ define
\[ L_{f,x}:\Ecal_x\to  \Ecal_{f(x)},\quad
 (L_{f,x}S)(u,v)=D_xf(S((D_xf)^{-1}u,(D_xf)^{-1}v)),\]
where $(D_xf)^{-1}$ is the inverse linear map $(D_xf)^{-1}:T_{f(x)}M\to T_xM$, and $u,v\in T_{f(x)}M$.
Set \[a_{f,x}=(f_*\nabla^0-\nabla^0)_{f(x)}\in\mathcal E_{f(x)}.\]
Then $f_*\nabla=\nabla$ for $\nabla=\nabla^0+S$ if and only if
\begin{equation}
 S(f(x))=L_{f,x}S(x)+a_{f,x}.
 \label{eq:consel}
\end{equation}

\end{lemma}
Note that for two diffeomorphisms $f,g$, pushforward is functorial:
$(f\circ g)_*=f_*g_*$. Hence \[L_{f\circ g,x}=L_{f,g(x)}L_{g,x},\qquad a_{f\circ g,x}=L_{f,g(x)}a_{g,x}+a_{f,g(x)}.\]
These precisely give the cocycle identity. We call the linear cocycle given by $L_{f,x}$ \emph{connection linear cocycle over $f$} and the affine cocycle given by $L_{f,x}+a_{f,x}$ \emph{connection affine cocycle over $f$}.

\subsection{Regularizing the measurable invariant connection}
In this subsection, we retain \cref{setting}. Furthermore, in \cref{thm:reduction}, we already show that 
\begin{enumerate}
\item there is the previously obtained H\"older compact reduction $P\subset\textrm{Fr}(TM)$ has structure group $C$ centralizing $\rho(G)$ and satisfies 
\[D\alpha(\gamma)_xP_x=P_{\alpha(\gamma)x}\rho(\gamma);\]
\item there is a H\"older metric $h$ satisfies \cref{prop:metric}.
\end{enumerate}

\subsubsection{Existence of the measurable invariant connection}
One important application of \cref{thm:cocsup} is that there is an $\alpha(\Gamma)$-invariant measurable connection $\nabla^m$. This is from standard argument, for instance, see \cite{ZimmerAlgebraicHull}. For the reader's convenient, let us summary the argument.

\begin{proposition}[A measurable invariant torsion-free connection] \label{msblecon}
Let $G$ and $\Gamma$ be a  higher rank simple group and a lattice in $G$, respectively, as in \cref{setting}.  Suppose that $\Gamma$ acts by $C^2$ diffeomorphisms on a closed smooth manifold $M$ and preserves an ergodic probability measure $\mu$, $\alpha:\Gamma\to \Diff^{2}_{\mu}(M)$. There is a measurable torsion-free affine connection $\nabla^m$ on $M$ that is $\Gamma$-invariant almost everywhere.
\end{proposition}

\begin{proof}
Fix a linear identification $V\simeq\Rbb^d$. In this identification the model tensor space is
\[ W=\Hom(\textrm{Sym}^2\Rbb^d,\Rbb^d),\qquad \textrm{Aff}(W)=\textrm{GL}(W)\ltimes W. \]
The multiplication is $(A,t)(B,u)=(AB,Au+t)$, corresponding to composition of $w\mapsto Aw+t$. This is a real algebraic group, realized as
\[ (A,t)\longmapsto \begin{bmatrix}A&t\\0&1\end{bmatrix}  \quad\hbox{in }\textrm{GL}(W\oplus\Rbb). \]

Choose a smooth background Riemannian metric on $M$ and a finite smooth local orthonormal frame cover for this metric,
with all coordinate calculations performed in slightly larger coordinate neighborhoods. 

We can find a measurable frame $e_x:\Rbb^d\to T_xM$ such that the norm and conorm is essentially bounded. This induces measurable identification; 
\[ \beta_x:W\longrightarrow\mathcal E_x, \qquad (\beta_xw)(u,v)=e_xw(e_x^{-1}u,e_x^{-1}v). \]
In these coordinates the affine connection cocycle is
\[\Fsc(\gamma,x)(w)=\beta_{\alpha(\gamma)x}^{-1}(L_{\alpha(\gamma),x}\beta_xw+a_{\alpha(\gamma),x}). \]
For each fixed $\gamma\in\Gamma$, the first and second derivatives of $\alpha(\gamma)$, its inverse first derivatives, and the reference coefficients are bounded. Therefore $\Fsc(\gamma,x)$ and its inverse have essentially bounded matrix norms. Hence, this verifies the integrability condition \eqref{L2} in \cref{thm:cocsup}.

Now apply \cref{thm:cocsup} to the target $\Fsc$. There exists a measurable map $\phi:M\to\textrm{Aff}(W)$, a continuous homomorphism $\pi:G\to\textrm{Aff}(W)$, and a measurable cocycle $\kappa$ with values in a compact subgroup $K$ centralizing $\pi(G)$, such that \begin{equation} \label{eq:consup} \Fsc(\gamma,x)=\phi({\alpha(\gamma)x})^{-1}\pi(\gamma)\kappa(\gamma,x)\phi(x) \quad\mu\textrm{-almost everywhere}.\end{equation}

When we write $\pi(g)w=A(g)w+t(g)$ for $w\in W$, the representation $\pi$ on $W\oplus\Rbb$ is realized by  \[ (w,r)\longmapsto(A(g)w+r\,t(g),r).\] Especially, the subspace $W\oplus0$ is invariant and its one-dimensional quotient is a trivial representation. Complete reducibility gives an invariant complementary line. The vector on that line whose last coordinate is $1$, say $(w_0,1)$, is fixed by $G$. Therefore $w_0$ is fixed by the affine action $\pi(G)$. The other way to see this is using the fact $H^1(G,W)=0$.

Since $K$ centralizes $\pi(G)$, the entire orbit $Kw_0$ consists of fixed points of $\pi(G)$. Define \[ w_* = \int_K k w_0\,dm_K(k),\] where $m_K$ is normalized Haar measure on $K$. Affine maps preserve averages with total mass one. Consequently $w_*$ is fixed by $\pi(G)$, and left invariance of $m_K$ shows that $w_*$ is also fixed by every element of $K$.

Set \[ s^m(x)=\phi(x)^{-1}w_*,\qquad S^m_x=\beta_xs^m(x)\in\Ecal_x, \qquad \nabla^m=\nabla^0+S^m. \]
The \eqref{eq:consup} gives
\[ \Fsc(\gamma,x)s^m(x)=\phi(\alpha(\gamma)x)^{-1}\pi(\gamma)k(\gamma,x)w_* =s^m(\alpha(\gamma)x).\]
After applying $\beta_{\alpha(\gamma)x}$, this is exactly \eqref{eq:consel} for $S^m$. Thus $\nabla^m$ is invariant. It is torsion free because $S^m_x$ is symmetric. Since $\Gamma$ is countable and preserves $\mu$, one can obtain a single invariant conull set on which all these equivariance equations hold.
\end{proof}

\subsubsection{Choosing an Anosov element for the tensor cocycle}
We would like to decompose the connection bundle $\Ecal$ into uniform contracting, neutral, and uniform expanding subbundles with respect to an Anosov diffeomorphism in order to use \cref{prop:conn}. However, the original Anosov element $\alpha(\gamma_0)$ may not give an (measurable) isometric neutral subbundle from cocycle superrigidity since we are considering the connection bundle with fiber $W$ instead of the tangent bundle with the derivative cocycle. Nevertheless, the weight of the induced representation of $G$ on $W$ is a finite set, hence, we can find a n element in $\Delta\subset \Gamma\cap MA$, such that $\alpha(\delta)$ is Anosov diffeomorphism (by \cref{prop:metric}) and this gives the desired Anosov diffeomorphism:

\begin{proposition}\label{propconntor}
There exists $\delta\in\Gamma$ such that $f=\alpha(\delta)$ is Anosov and the connection linear cocycle $L_{f,x}$ on $\Ecal$ satisfies all three linear hypotheses of Proposition~\ref{prop:conn}. In fact the invariant metric on its neutral subbundle is already H\"older.
\end{proposition}
\begin{proof}
Choose $A,M,\Delta$ as in \cref{CartanGamma}. As before, we write the weight space decomposition
\[ V=\bigoplus_{\lambda\in\Sigma_V}V_\lambda,\qquad \rho(\exp H)v=e^{\lambda(H)}v\quad(v\in V_\lambda),\]
where $\Sigma_V$ is the finite set of actual restricted weights, and $0\notin\Sigma_V$. The compact groups $C$ and $\rho(M)$ commute and preserve each $V_\lambda$. Choose an inner product on $V$ by averaging inner products on these spaces over $C\rho(M)$ and making different weight spaces orthogonal. It induces an inner product $\langle\cdot,\cdot,\rangle_W$ and a metric $||\cdot||_W$ on $W=\Hom(\textrm{Sym}^2V,V)$.

For $R\in\textrm{GL}(V)$ and $S\in W$, define $\sigma:G\to\GL(W)$, \[ (\sigma(R)S)(u,v)=R\,S(R^{-1}u,R^{-1}v) \] for $u,v\in V$. This gives the induced representation on $W$ by $\sigma:G\to \GL(W)$.  The induced inner product makes $\sigma(C\rho(M))$ orthogonal. Let us write the weight space decomposition (with respect to $A$) of $(W,\sigma)$ is
\[W=\bigoplus_{\chi\in\Sigma_W}W_\chi, \qquad  \sigma(\rho(\exp H))\restrict{W_\chi}=e^{\chi(H)}I\] for all $H\in\alie$.
By the definition of $\sigma$, every weight $\chi$ has the form $\nu-\lambda-\beta$ for tangent weights $\lambda,\beta,\nu$: on a bilinear map from $V_\lambda\times V_\beta$ to $V_\nu$, the output contributes $e^{\nu(H)}$ and the two inverse inputs contribute $e^{-\lambda(H)}e^{-\beta(H)}$.

Let $\Ccal=\{\nu-\lambda-\beta: \beta,\lambda,\nu\in \Sigma_V\}\cup\Sigma_V$. Since $\Ccal$ is finite, there is a $\delta\in \Delta$ such that  $H=\log a_\delta$ avoiding the kernels of all tangent weights
and all nonzero tensor weights.

For $b\in P_x$ define
\[ \mathcal I_{b_x}:W\to\Ecal_x,\quad (\mathcal I_{b_x} S_x)(X,Y)=b_xS_x(b_x^{-1}u,b_x^{-1}y), \quad\forall u,v\in T_xM.\]
Then, for $c\in C$, $\mathcal I_{bc}=\mathcal I_b\sigma(c)$. Since $\sigma(c)$ preserves each $W_\chi$ and is orthogonal,
\[ \Ecal_x^\chi\coloneqq\Ical_b(W_\chi),\qquad\|\Ical_bS\|_x\coloneqq\|S\|_W\] define, independently of $b$, H\"older tensor subbundles and a H\"older metric on $\mathcal E$. 
Their regularity follows by choosing local H\"older adapted sections of $P$.

Choose $b_x\in P_x$ and $b_{f(x)}\in P_{f(x)}$. There is $k(x)\in C$ with \[D_xf\,b_x=b_{f(x)}\rho(a_\delta)\rho(m_\delta)k(x).\]
Direct substitution into the definition of $L_{f,x}$ gives
\[ L_{f,x}\mathcal I_{b_x}=\mathcal I_{b_{f(x)}}\sigma(\rho(a_\delta))\sigma(\rho(m_\delta)k(x)).\]
The $\sigma(\rho(m_\delta)k(x))$ factor is orthogonal and preserves each $W_\chi$. Consequently, for every $x$ and $S\in\mathcal E_x^\chi$,
\[\|L_{f,x}S\|_{f(x)}=e^{\chi(H)}\|S\|_x, \qquad \|L_{f,x}^{(i)}S\|_{f^i(x)}=e^{i\chi(H)}\|S\|_x,\] for $i\in \Zbb$.
Set
\[\Ecal^- =\bigoplus_{\chi(H)<0}\Ecal^\chi, \quad \Ecal^0 =\bigoplus_{\chi(H)=0}\Ecal^\chi, \quad \Ecal^+ =\bigoplus_{\chi(H)>0}\mathcal E^\chi.\] The sum is orthogonal, and these are H\"older invariant subbundles. The finitely many negative and positive values give uniform exponential contraction on $\mathcal E^-$ and for the inverse on $\mathcal E^+$. On $\mathcal E^0$ the induced metric is exactly invariant. 
\end{proof}

\subsubsection{Regularizing $\nabla^m$} Now we are ready to improve the regularity of $\nabla^m$.

\begin{proposition}[A H\"older invariant torsion-free connection]
\label{prop:holconnect}
Let $\nabla^m$ be the measurable torsion-free $\alpha(\Gamma)$-invariant connection given by \cref{msblecon}. Then there is a H\"older torsion-free $\Gamma$-invariant connection $\nabla$ equal to $\nabla^m$ almost everywhere.
\end{proposition}
\begin{proof}
Write $\nabla^m=\nabla^0+S^m$, where $S^m$ is a measurable section of $\mathcal E$. Its invariance under the whole group is exactly \eqref{eq:consel} for every $\alpha(\gamma)$, almost everywhere. Choose the particular Anosov element from \Cref{propconntor}. Lemma~\ref{lem:conaff} gives a H\"older affine cocycle and the three tensor subbundles verify every linear assumption of \Cref{prop:conn}. Applying that proposition to the given measurable section $S^m$ produces a H\"older section $S$ with $S=S^m$ almost everywhere. Put $\nabla=\nabla^0+S$. It is H\"older and torsion free because $S$ is symmetric.

The only remaining duty is showing that $\nabla$ is $\alpha(\Gamma)$-invariant.  Fix $\gamma\in\Gamma$ and put $g=\alpha(\gamma)$. Since $g$ preserves $\mu$, on a conull subset, both $S(x)=S^m(x)$ and $S(g(x))=S^m(g(x))$ hold. The original measurable invariance therefore gives \[ S(g(x))-L_{g,x}S(x)-a_{g,x}=0 \quad\textrm{for almost every }x.\]
The left-hand side is a continuous section over $x$ with values in $\mathcal E_{g(x)}$. A continuous section vanishing $\mu$ almost everywhere for a fully support measure $\mu$,  vanishes everywhere. Thus, $g_*\nabla=\nabla$ at all points. This holds for each $\gamma\in\Gamma$. The countability of $\Gamma$ also allows for one common invariant conull subset.
\end{proof}

\section{Local Lipschitz and $C^{1+\beta}$ regularity of $P$}\label{sec5}

 Let us summary what we had so far: We are in \cref{setting}. Let $V=\Rbb^d$ be the model vector space. Under the preceding hypotheses, there exist 

\begin{enumerate}
\item a continuous representation with no zero weight $\rho:G\to \GL(V)$;
\item  a closed subgroup $C\subset \GL(V)$;
\item a H\"older principal $C$-reduction $P\subset \Fr(TM)$;
\item a H\"older metric $h$ on $TM$ such that \cref{prop:metric} holds;
\item a H\"older $\alpha(\Gamma)$-invariant connection $\nabla$ on $M$; 
\end{enumerate}
such that:
\[D_x\alpha(\gamma)(P)_x
       =(P)_{\alpha(\gamma)x}\rho(\gamma)\] for every
       $\gamma\in\Gamma$ and $x\in M$;

The main goal of this section is proving \cref{prop:reduction1}. 
Let us fix a closed subgroup $C'<C$ and H\"older $C'$-reduction $P'\subset P$ that satisfy \begin{equation}\label{eq:equiv'} D_{x}\alpha(\gamma)P'_x=P'_{\alpha(\gamma)(x)} \rho(\gamma)\end{equation} for every $\gamma\in\Gamma$ and $x\in M$. In this section, using the H\"older invariant connection $\nabla$, we show that $P'$ is locally Lipschitz first, \cref{lipP}. After that, we will show that $P'$ is indeed $C^{1+\beta}$, \cref{derivP}. Finally, we will find the final  reduction $(P_*,C_*)$ in \cref{prop:reduction1} which will be used in \cref{sec6}.

\subsection{Local Lipschitzness}
For each real irreducible components of $V$, we select a primitive weight. We first estimate the variation of the quotient section of $P'$ along curves tangent to the corresponding weight bundles, and then along their $\Gamma$-translates. Finitely many translated full weight spaces span $V$. A local estimate on the lengths of paths in these directions gives local Lipschitzness \cref{lipP}.

\subsubsection{Primitive weight direction} Let $\Sigma=\Sigma(G,A)\subset\alie^*$ be the restricted root system. Choose $H_+\in\alie$ such that $\alpha(H_+)\ne0$ for every $\alpha\in\Sigma$, and set \[ \Sigma^+=\{\alpha\in\Sigma:\alpha(H_+)>0\}. \]

Let $\Pi$ be the corresponding set of simple roots, that is, the positive roots that cannot be written as a sum of two positive roots.

Fix a Weyl group invariant positive definite inner product $(\cdot,\cdot)$ on $\alie^*$, for example the dual of the inner product induced by the Killing form on $\alie$. For each $\alpha\in\Sigma$, write $\alpha^\vee\in\alie$ for the coroot, normalized by
\[ \chi(\alpha^\vee)=\frac{2(\chi,\alpha)}{(\alpha,\alpha)} \] for $\chi\in\alie^*$. The weight lattice of the restricted root system is
\[ \Lambda_\Sigma=\{\chi\in\alie^*:\chi(\alpha^\vee)\in\Zbb\text{ for every }\alpha\in\Sigma\}.\]
A nonzero $\xi\in\Lambda_\Sigma$ is \emph{primitive} if $\Lambda_\Sigma\cap\Rbb\xi=\Zbb\xi$.

Let $\rho:G\to\GL(V)$ be the superrigidity representation and let $\Sigma_V\subset\alie^*$ be its set of restricted weights. Recall that $0\notin\Sigma_V$. Let $H$ be the algebraic model of $G$ as in \cref{Borel}. 

\begin{lemma}\label{lem:primitive}
Every real irreducible summand $U_i$ of $V$ contains a restricted weight $\xi$ such that
\[ \xi(\alpha^\vee)\in\{0,1\} \quad\text{for every }\alpha\in\Sigma^+, \qquad \xi(\alpha_0^\vee)=1 \quad\text{for some }\alpha_0\in\Pi.\]
In particular, $\xi$ is primitive in $\Lambda_\Sigma$. Furthermore, for any $\lambda,\nu\in\Sigma_V$,
\[ \lambda-\nu=t\xi,\qquad t>0\quad\Longrightarrow\quad t\in\Zbb_{\geq1}.\]
\end{lemma}

\begin{proof}
Fix an irreducible summand $U_i$. Among its restricted weights of minimum norm, choose $\xi$ so that $\xi(H_+)$ is maximal. Since $0\notin\Sigma_V$, we have $\xi\ne0$.

Fix $\alpha\in\Sigma^+$. The Weyl reflection \[ s_\alpha\xi=\xi-\xi(\alpha^\vee)\alpha\] is again a weight of $U_i$.  Since the inner product is Weyl group invariant, $s_\alpha\xi$ has the same norm as $\xi$. Thus \[ 0\leq \xi(H_+)-(s_\alpha\xi)(H_+)=\xi(\alpha^\vee)\alpha(H_+).\]
As $\alpha(H_+)>0$, this proves $\xi(\alpha^\vee)\geq0$.

For each $\lambda\in\Sigma_V$ and $\alpha\in\Sigma$, $\lambda(\alpha^\vee)$ is an eigenvalue of the standard diagonal element in the representation of the restricted-root $\mathfrak{sl}_2$-triple $(X_\alpha,\alpha^\vee,Y_\alpha)$. It is an integer, so $\Sigma_V\subset\Lambda_\Sigma$. Suppose $\xi(\alpha^\vee)\geq2$ for some $\alpha\in\Sigma^+$. Choose a nonzero vector $v$ in the $\xi$-weight space of $U_i$. On each irreducible $\mathfrak{sl}_2$-summand, the lowering operator is injective on every positive eigenspace of the diagonal element. Thus $d\rho(Y_\alpha)v\ne0$, and for every $H\in\alie$,
\[ d\rho(H)d\rho(Y_\alpha)v=d\rho([H,Y_\alpha])v+d\rho(Y_\alpha)d\rho(H)v=(\xi(H)-\alpha(H))\,d\rho(Y_\alpha)v. \]
Thus, $\xi-\alpha$ is a weight of $U_i$. Its norm satisfies 
\[ \|\xi-\alpha\|^2-\|\xi\|^2=(\alpha,\alpha)-2(\xi,\alpha)=(\alpha,\alpha)\left(1-\xi(\alpha^\vee)\right)<0,\]
contradicting the choice of $\xi$. Hence $\xi(\alpha^\vee)\in\{0,1\}$ for every $\alpha\in\Sigma^+$.

Finally, the simple coroots span $\alie$. Since $\xi\ne0$, it cannot vanish on all of them. There is therefore some $\alpha_0\in\Pi$ with $\xi(\alpha_0^\vee)=1$. If $t\xi\in\Lambda_\Sigma$ for some $t\in\Rbb$, then \[ t=(t\xi)(\alpha_0^\vee)\in\Zbb,\]
which proves primitivity. For $\lambda,\nu\in\Sigma_V$, integrality gives $\lambda-\nu\in\Lambda_\Sigma$, even if $\lambda$ and $\nu$ occur in different summands. Thus $\lambda-\nu=t\xi$ with $t>0$ implies that $t$ is a positive integer.
\end{proof}

\subsubsection{Transport along a primitive weight direction}

Retain the adapted norms $\|\cdot\|_V$ on $V$ and $\|\cdot\|_h$ on $TM$ from \cref{sec:adaptmetric}. Let $\|\cdot\|$ be the induced operator norm on $\textrm{End}(V)$.

As before, $\Sigma_V$ denotes the set of distinct restricted weights.
\[V=\bigoplus_{\lambda\in\Sigma_V}V_\lambda,\qquad TM=\bigoplus_{\lambda\in\Sigma_V}E^\lambda, \qquad E^\lambda_x=b_xV_\lambda
\quad(b_x\in P_x).\]
In the averaging construction of \cref{sec:adaptmetric}, start with an inner product making the $V_\lambda$ mutually orthogonal. Both $C$ and $\rho_H(T_c)$ preserve each summand, so the averaged inner product retains this orthogonality. Every element of $C$ preserves each $V_\lambda$, so $E^\lambda_x$ is independent of the chosen frame. Every $b_x\in P_x$ is an isometry for the resulting metrics.

Fix a compact subgroup $C'\subset C$ and a H\"older $C'$-reduction $P'\subset P$ satisfying \eqref{eq:equiv'}:
\[D_x\alpha(\gamma)P'_x =P'_{\alpha(\gamma)x}\rho(\gamma).\]
Recall that $C'$ commutes with $\rho(G)$. The constants below may depend on this fixed reduction.

Let $\Delta\subset q(\Gamma)$ be the subgroup of \cref{CartanGamma}. For each $\delta=a_\delta m_\delta\in\Delta$, choose $\gamma_\delta\in\Gamma$ with $q(\gamma_\delta)=\delta$, and put
\[
 f_\delta=\alpha(\gamma_\delta),\quad H_\delta=\log a_\delta,\quad
 D_\delta=\rho_H(a_\delta),\quad R_\delta=\rho_H(m_\delta).
\]
By \cref{Borel}, $\rho(\gamma_\delta)=\rho_H(\delta)=D_\delta R_\delta$. The vectors $H_\delta$ form a full lattice in $\alie$, and $D_\delta$ acts on $V_\lambda$ by $e^{\lambda(H_\delta)}I_{V_\lambda}$. The orthogonal map $R_\delta$ preserves every $V_\lambda$ and commutes with $D_\delta$ and $C'$. For a fixed $\delta$ we use the powers $f_\delta^j=\alpha(\gamma_\delta^j)$.

Equivariance and the adapted metric give
\begin{equation}\label{eq521}
\|D_x f_\delta\,v\|_h=e^{\lambda(H_\delta)}\|v\|_h \end{equation} for $v\in E^\lambda_x$.
In particular, if $c$ is tangent to $E^\xi$, then
\begin{equation}\label{eq522}
\ell_h(f_\delta^j\circ c)=e^{j\xi(H_\delta)}\ell_h(c) \end{equation} for every $j\in\Zbb$.

Fix a primitive weight $\xi$ supplied by \cref{lem:primitive}.
Define
\[\nlie_\xi = \left\{ U\in\textrm{End}(V) :  U(V_\lambda)\subset \bigoplus_{\substack{t>0\\\lambda+t\xi\in\Sigma_V}} V_{\lambda+t\xi} \quad\text{for every }\lambda\in\Sigma_V \right\}.\] Set $N_\xi=I_V+\nlie_\xi$.

\begin{lemma}\label{lem:primnil} The space $\nlie_\xi$ is a nilpotent associative algebra, hence a nilpotent Lie algebra under commutators. The set $N_\xi$ is a unipotent group normalized by $C'$.
\end{lemma}

\begin{proof}
Composition adds positive multiples of $\xi$ to the weight, so $\nlie_\xi$ is closed under composition. Choose $H\in\alie$ with $\xi(H)>0$.  Since there are finitely many weights, every sufficiently long product of elements of $\nlie_\xi$ vanishes.

Therefore \[(I+U)^{-1}=I-U+U^2-\cdots\] is a finite sum for $U\in\nlie_\xi$. Hence, $N_\xi$ is a unipotent group. Finally, $C'$ preserves each weight space, so conjugation by $C'$ preserves $\nlie_\xi$ and $N_\xi$.
\end{proof}

Let $c:[0,1]\to U$ be a piecewise $C^1$ curve from $x$ to $y$, and let $T_c:T_xM\to T_yM$ be its $\nabla$-parallel transport. Write $\nabla=d+A$ in a smooth local frame $e$ over $U$. A field $Y(t)=e_{c(t)}v(t)$ along $c$ is parallel if its coordinates satisfy \[v'(t)+A_{c(t)}(\dot c(t))v(t)=0\] on each $C^1$ piece. Let $S_c(t)$ solve \begin{equation}\label{pt_eq} S_c'(t)=-A_{c(t)}(\dot c(t))S_c(t),\qquad S_c(0)=I. \end{equation} The parallel transport along $c$ is \[T_c=e_{c(1)}\circ S_c(1)\circ e_{c(0)}^{-1}:T_{c(0)}M\longrightarrow T_{c(1)}M.\]
For a piecewise $C^1$ path, transport is the composition of the transports on its finitely many $C^1$ pieces.

For arbitrary $b_x\in P'_x$ and $b_y\in P'_y$, set $Q_c=b_y^{-1}T_cb_x\in\GL(V)$.

\begin{lemma}\label{lem:shtrans} There exist $r_0,K_0>0$ such that, for every piecewise $C^1$ curve $c$ with $\ell_h(c)\le r_0$ and every pair of adapted endpoint frames, \[ \|Q_c\|\le K_0.\]
There is a finite open cover of $M$ with a local H\"older section of $P'$ on each member, such that the following also holds. If $\ell_h(c_j)\to0$ and both endpoint frames of each sufficiently short $c_j$ come from the section of a cover member containing $c_j$, then
\[ Q_{c_j}\longrightarrow I.\]
\end{lemma}

\begin{proof}
Choose a finite open cover $\{U_j\}$ of $M$ with $U_j\Subset\Ocal_j$, where each $\Ocal_j$ is a coordinate neighborhood carrying a smooth reference frame $e_j$ and a H\"older local section $b_j^0$ of $P'$. In the frame $e_j$, write \[ b_{j,z}^0=e_{j,z}B_j(z),\qquad \nabla=d+\omega_j.\]
By compactness of the closures $\overline{U_j}$ and finiteness of the cover, there are constants $L,H,K>0$ and $0<\theta\le1$ such
that, for every $j$, all $z,w\in U_j$, and every $v\in T_zM$,
\begin{align*}
 \max\{\|B_j(z)\|,\|B_j(z)^{-1}\|\}&\le L, &  \|\omega_{j,z}(v)\|&\le K|v|_h,\\
 \|B_j(z)-B_j(w)\|&\le H\,d_h(z,w)^\theta .
\end{align*}

Choose $r_0>0$ smaller than a Lebesgue number of this cover for $d_h$, the length distance of $h$. Every curve with $\ell_h(c)\le r_0$ lies in some $U_j$. For calculations there, write $e,b^0,B_0,\omega$ for the corresponding local data. 

Let $c:[0,1]\to U_j$ be a sufficiently short piecewise $C^1$ curve from $x$ to $y$. Let $\ell=\ell_h(c)$. For $0\le t\le1$, set
\[ S(t)=e_{c(t)}^{-1}T_{c\restrict{[0,t]}}e_x,\qquad L(t)=\int_0^t|\dot c(s)|_h\,ds .\]
Then $S(0)=I$. Write $S_c=S(1)$ for the endpoint matrix. The parallel transport equation \eqref{pt_eq} gives
\[ S(t)=I-\int_0^t\omega_{c(s)}(\dot c(s))S(s)\,ds. \]
Taking operator norms gives \begin{equation}\label{eq541} \|S(t)\|\le1+\int_0^t K|\dot c(s)|_h\,\|S(s)\|\,ds.\end{equation} By \cref{gronwallineq}, $\|S(t)\|\le e^{KL(t)}$. Applying this bound to the integral expression for $S(t)-I$ gives
\[ \|S(t)-I\|\le\int_0^t K|\dot c(s)|_h\,e^{KL(s)}\,ds =e^{KL(t)}-1,\]
since $L'(s)=|\dot c(s)|_h$ almost everywhere. Taking $t=1$ and
writing $\ell=L(1)=\ell_h(c)$, we obtain
\begin{equation}\label{eq:primbd}  \|S_c\|\le e^{K\ell},  \qquad  \|S_c-I\|\le e^{K\ell}-1. \end{equation}
For endpoint frames from the common section $b^0$,
\[Q_c=B_0(y)^{-1}S_cB_0(x)\Rightarrow  \|Q_c\|\le L^2e^{Kr_0}.\]
Also,
\[Q_c-I=B_0(y)^{-1}(S_c-I)B_0(x)  +B_0(y)^{-1}\left(B_0(x)-B_0(y)\right),\]
so
\begin{equation}\label{eq:primQ} \|Q_c-I\|  \le L^2(e^{K\ell}-1)+LH\ell^\theta. \end{equation}
This tends to zero with $\ell\to0$.

Arbitrary adapted endpoint frames differ from those of $b^0$ by elements $k_x,k_y\in C'$. The corresponding map is $k_y^{-1}Q_ck_x$, which has the same operator norm because $k_x,k_y$ are orthogonal. This proves the uniform bound for arbitrary frames.
\end{proof}

For $T\in\textrm{End}(V)$ and $\nu,\lambda\in \Sigma_V$, define its weight block by
\[T_{\nu\to\lambda}\coloneqq\textrm{Pr}_\lambda\circ T\restrict{V_\nu} :V_\nu\longrightarrow V_\lambda.\]
Here $\textrm{Pr}_\lambda$ is the orthogonal projection onto $V_\lambda$. The input weight is $\nu$ and the output weight is $\lambda$.

\begin{lemma}\label{lem:QinCN}
Let $c:[0,1]\to M$ be a $C^1$ curve tangent to $E^\xi$.
Write $x=c(0)$, $y=c(1)$, and let $T_c:T_xM\to T_yM$ be its $\nabla$-parallel transport.
As before, for arbitrary $b_x\in P'_x$ and $b_y\in P'_y$, set $Q_c=b_y^{-1}T_cb_x$.
Then, $Q_c\in C'N_\xi$.
More precisely, the weight-diagonal part
\[D_0\coloneqq\sum_{\lambda\in\Sigma_V} \textrm{Pr}_\lambda Q_c\textrm{Pr}_\lambda \] belongs to $C'$, and $Q_c-D_0\in\nlie_\xi$.
\end{lemma}

\begin{proof}
The invariance of $\nabla$ means that \[\alpha(\gamma)_*(\nabla_XY)=\nabla_{\alpha(\gamma)_*X}(\alpha(\gamma)_*Y) .\]
Consequently, \begin{equation}\label{eq:Tfc} T_{\alpha(\gamma)\circ c}=D_y\alpha(\gamma)\circ T_c\circ(D_x\alpha(\gamma))^{-1}.\end{equation}

Fix $\delta\in\Delta$ and let $c_j=f_\delta^j\circ c$. Choose adapted frames at the endpoints of $c_j$. Equivariance for $\gamma_\delta^j$ gives $k_j(x),k_j(y)\in C'$ such that
\[D_zf_\delta^j b_z=b^{(j)}_{f_\delta^jz}\rho_H(\delta)^jk_j(z),\qquad z=x,y.\]
Substitution into \eqref{eq:Tfc} gives
\begin{equation}\label{eq:Q_C_nu} Q_{c_j} =\rho_H(\delta)^jk_j(y)Q_ck_j(x)^{-1}\rho_H(\delta)^{-j}. \end{equation}
All compact factors are orthogonal and preserve the weight spaces. Therefore
\begin{equation}\label{eq:primblock}
\|(Q_{c_j})_{\nu\to\lambda}\| =e^{j(\lambda-\nu)(H_\delta)}       \|(Q_c)_{\nu\to\lambda}\|.
\end{equation}
This holds for positive and negative integers $j$.

Suppose $\lambda\ne\nu$ and $\lambda-\nu$ is not positively proportional to $\xi$, i.e., $\lambda-\nu\notin\Rbb_{>0}\xi$.
The conditions $\xi(H)<0<(\lambda-\nu)(H)$ define a nonempty open cone in $\alie$: the two functionals are either linearly independent or negatively proportional. A full lattice meets every nonempty open cone, so we may choose $\delta\in\Delta$ with
\[\xi(H_\delta)<0<(\lambda-\nu)(H_\delta).\]
By \eqref{eq522},
$\ell_h(c_j)\longrightarrow0$ as $j\to \infty$. For large $j$, choose the endpoint frames of $c_j$ from one common local H\"older section in the finite collection of \cref{lem:shtrans}. Then \[ Q_{c_j}\longrightarrow I.\]
Its off-diagonal $(\nu\to\lambda)$ block tends to zero.
Since the exponential factor in \eqref{eq:primblock} grows, this forces
\[(Q_c)_{\nu\to\lambda}=0.\]
It follows that
\[Q_c=D_0+U, \qquad U\in\nlie_\xi,\] with $D_0$ block diagonal with respect to the weight decomposition.

It remains to show that $D_0\in C'$. Keep the original endpoint frames $b_x$ and $b_y$ fixed. Choose one $\delta\in\Delta$ with $\xi(H_\delta)<0$. By \eqref{eq522}, the lengths of the curves $c_j=f_\delta^j\circ c$ tend to zero. For each sufficiently large $j$, choose the two endpoint frames of $c_j$ from the same local H\"older section as in \cref{lem:shtrans}. Then $Q_{c_j}\to I$.
Taking the weight-diagonal part is a continuous linear map,
\[ D_j\coloneqq\sum_{\lambda\in\Sigma_V} \textrm{Pr}_\lambda Q_{c_j}\textrm{Pr}_\lambda \longrightarrow I. \]

We compute $D_j$ from \eqref{eq:Q_C_nu}. Every factor surrounding $Q_c$ in that equation preserves each weight space. Thus taking weight-diagonal parts replaces $Q_c$ by $D_0$. On the block $V_\lambda\to V_\lambda$, the positive factors from $\rho_H(\delta)^j$ and $\rho_H(\delta)^{-j}$ are the scalars $e^{j\lambda(H_\delta)}$ and $e^{-j\lambda(H_\delta)}$. They cancel. As $R_\delta$ is the orthogonal factor of $\rho_H(\delta)$, we obtain \[ D_j=R_\delta^j k_j(y)D_0 k_j(x)^{-1}R_\delta^{-j}. \]

On the other hand,
\[ k_j(y)^{-1}R_\delta^{-j}(D_j-I)R_\delta^j k_j(x) =D_0-k_j(y)^{-1}k_j(x). \]
Both outside factors on the left are orthogonal for the adapted
inner product. Hence, in the induced operator norm, as $j\to \infty$,
\[\|D_0-k_j(y)^{-1}k_j(x)\|=\|D_j-I\|\longrightarrow0. \]
The products $k_j(y)^{-1}k_j(x)$ lie in the compact group $C'$ and converge to $D_0$. Hence $D_0\in C'$. We have shown that $Q_c=D_0+U$ with $U\in\nlie_\xi$. Since $D_0^{-1}$ preserves every weight space, $D_0^{-1}U\in\nlie_\xi$. It follows that
\[ Q_c=D_0\bigl(I+D_0^{-1}U\bigr)\in C'N_\xi.\]
\end{proof}

\begin{remark}
Although $\nabla$ is only H\"older, its parallel transport along a prescribed $C^1$ curve is well-defined by a linear ordinary differential equation with continuous coefficients.
\end{remark}

\begin{lemma}\label{lem:LipQ}
There are $r_0,C_1>0$ such that, for every $C^1$ curve $c$ tangent to $E^\xi$ with $\ell_h(c)\le r_0$, and every choice of adapted endpoint frames,
\[ \|Q_c-D_0\|\le C_1\ell_h(c),\]
where $D_0\in C'$ is the weight-diagonal part of $Q_c$. Consequently,
\[\textrm{dist}(Q_c,C') \coloneqq\min_{k\in C'}\|Q_c-k\|\le C_1\ell_h(c).\]
In particular, \[\|(Q_c)_{\nu\to\lambda}\|\le C_1\ell_h(c)\qquad(\nu\ne\lambda).\]
\end{lemma}

\begin{proof}
By \cref{lem:QinCN}, every possibly nonzero off-diagonal block satisfies
\[ \lambda-\nu=t\xi,\qquad t>0.\]
By \cref{lem:primitive}, $t\in \Zbb$, $t\ge 1$. Choose $\delta\in\Delta$ with $a:=-\xi(H_\delta)>0$.
Decrease the short-curve constant $r_0$ from \cref{lem:shtrans} so that $r_0\le1$. Write $\ell=\ell_h(c)$.

If $\ell=0$, the curve is constant and $Q_c\in C'$. The assertion is immediate. For $0<\ell\le r_0$, find an integer $N$ such that $aN \le \log(r_0/\ell) \le a(N+1)$. Then \[r_0e^{-a}<e^{aN}\ell\le r_0.\]
Thus the expanded curve
$\wtd c= f_\delta^{-N}\circ c$ has length at most $r_0$ and satisfies $\|Q_{\widetilde c}\|\le K_0$.

For the $(\lambda,\nu)$-block $\lambda-\nu=t\xi$, \eqref{eq:primblock} with $j=-N$ gives
\[e^{taN}\|(Q_c)_{\nu\to\lambda}\|=\|(Q_{\widetilde c})_{\nu\to\lambda}\|\le K_0.\]
Therefore
\[\|(Q_c)_{\nu\to\lambda}\| \le K_0e^{-taN}\le K_0e^{ta}r_0^{-t}\ell^t.\]
Only finitely many weight pairs exist. Every corresponding exponent satisfies $t\ge1$, and $\ell\le1$, so $\ell^t\le\ell$.
Summing the finitely many block estimates proves
\[\|Q_c-D_0\|\le C_1\ell.\]
Since $D_0\in C'$, the distance estimate follows.
\end{proof}

For cosets in $\textrm{GL}(V)/C'$, define
\[
d_{C'}(AC',BC')
   :=\min_{k\in C'}\|A-Bk\|.
\]
Compactness gives attainment of the minimum.
Orthogonality of $C'$ gives
\[
d_{C'}(AC',BC')
   =\min_{k,k'\in C'}\|Ak-Bk'\|.
\]

\begin{corollary}[Local Lipschitzness along $E^\xi$]\label{cor:Lipalongprim}
Let $\Ocal$ be a sufficiently small coordinate neighborhood with a smooth reference frame $e_z:V\longrightarrow T_zM$.
On every smaller neighborhood with compact closure in $\Ocal$, there is $C>0$ such that every sufficiently short $C^1$ curve $c$ tangent to $E^\xi$ and starting in the smaller neighborhood satisfies
\[d_{C'}\left(B(x)C',B(y)C'\right)\le C\ell_h(c),\]
where $x=c(0)$, $y=c(1)$, and
\[b_x=e_xB(x)\in P'_x,\qquad b_y=e_yB(y)\in P'_y.\]
Equivalently, there is $k_c\in C'$ such that
\[\|B(y)k_c-B(x)\|\le C\ell_h(c).\] The quotient estimate holds independently of the chosen adapted endpoint frames.
\end{corollary}

\begin{proof}
Set $S_c=e_y^{-1}T_ce_x.$ Then \[Q_c=B(y)^{-1}S_cB(x) \Rightarrow S_cB(x)=B(y)Q_c.\] Let $D_0\in C'$ be as in \cref{lem:LipQ}.
Then \begin{equation}\label{eq:endpt} B(y)D_0-B(x)=(S_c-I)B(x)-B(y)(Q_c-D_0). \end{equation}

All sufficiently short curves starting in the smaller neighborhood stay in a fixed compact subset of $\Ocal$. There the endpoint matrices $B(x),B(y)$ are uniformly bounded, since they differ from matrices of a fixed local section by right multiplication by elements of $C'$.
Furthermore, \[\|S_c-I\|\le C_2\ell_h(c)\] by \eqref{eq:primbd}, and \[\|Q_c-D_0\|\le C_1\ell_h(c)\] by \cref{lem:LipQ}. Equation \eqref{eq:endpt} therefore gives\[ \|B(y)D_0-B(x)\|\le C\ell_h(c).\] Take $k_c=D_0$ and minimize over $C'$.
\end{proof}

\begin{remark} The estimate is in terms of the length of the chosen curve, not its endpoint distance. A bound by $d_h(x,y)$ will follow from \cref{joining}.\end{remark}

\subsubsection{Local Lipschitzness}

We transfer the directional estimate to lattice translates of the weight bundles.
\begin{lemma}\label{lem:primtrans}
Fix $\gamma\in\Gamma$ and an occurring primitive weight $\xi$. The estimate of \cref{cor:Lipalongprim} holds, with a possibly different constant, along curves tangent to
\[F^{\gamma,\xi}_z \coloneqq D_{\alpha(\gamma)^{-1}z}\alpha(\gamma)\left(E^\xi_{\alpha(\gamma)^{-1}z}\right).\]
For every $b_z\in P'_z$,
\[F^{\gamma,\xi}_z=b_z\rho(\gamma)V_\xi.\]
\end{lemma}

\begin{proof}
The bundle identity follows from equivariance of $P'$ and preservation of $V_\xi$ by $C'$. Given a curve $\widetilde c$ tangent to $F^{\gamma,\xi}$, its pullback $c=\alpha(\gamma)^{-1}\circ\widetilde c$ is tangent to $E^\xi$.
Choose endpoint frames by \[\wtd b_{\alpha(\gamma)z}=D_z\alpha(\gamma)b_z\rho(\gamma)^{-1}. \]
They belong to $P'$. Invariance of the connection gives
\[Q_{\widetilde c}=\rho(\gamma)Q_c\rho(\gamma)^{-1}.\]
Since $\rho(\gamma)$ centralizes $C'$,
\[\textrm{dist}(Q_{\widetilde c},C')\le\|\rho(\gamma)\|\,\|\rho(\gamma)^{-1}\|\textrm{dist}(Q_c,C').\]
Also,
\[\ell_h(c) \le L_\gamma\ell_h(\widetilde c),\]
because the derivative of the fixed diffeomorphism $\alpha(\gamma)^{-1}$ is bounded on compact $M$. If $\widetilde c$ is sufficiently short, \cref{lem:LipQ} applies to $c$. In the endpoint identity \eqref{eq:endpt} for $\widetilde c$, use a minimizing $k\in C'$ for $\textrm{dist}(Q_{\widetilde c},C')$ in place of $D_0$. This gives the claimed quotient estimate.
\end{proof}

\begin{lemma}\label{joining}
Let $X_1,\dots, X_d$ be continuous vector fields in a closed smooth manifold $M$ with dimension $d$. Assume that they are linearly independent (hence, span $TM$). Fix $x\in M$. There is an open neighborhood $\Ocal$ of $x$ and a constant $L$ such that for all $y\in\Ocal$, $x$ and $y$ can be joined by a concatenation of $C^1$ curves $c_1*\dots *c_l$ where $c_i$ tangents to $X_i$ with total length \[\sum_i \ell_h(c_i) <L d_h(x,y).\]\end{lemma}
\begin{proof}
The same proof for \cite[Lemma 3.4]{GHS} can be applied since we only ask the curves $c_i$ tangents to $X_i$. Note that even though the statement is written by stable distributions, the proof can be applied directly to any continuous transverse vector fields $X_1,\dots, X_d$ as well as the cited theorem \cite[Corollary 4.5]{KryszewskiPlaskacz}.  The only usage of stable distributions in \cite[Lemma 3.4.]{GHS} is the last line in the proof in order to guarantee that each curve is contained in the stable manifold using the unique integrability. However, in our case, we only ask that curves $c_i$ tangent to $X_i$, hence, the same proof shows \cref{joining}.
\end{proof}
Now we are ready to show the local Lipschitz regularity.
\begin{proposition}[Local Lipschitz regularity of $P'$]\label{lipP}
Let $C'\subset C$ be compact. For every H\"older $C'$-reduction $P'\subset P$ satisfying \eqref{eq:equiv'}, the quotient section
is locally Lipschitz, as a map\[M\longrightarrow\textrm{Fr}(TM)/C'. \] Equivalently, $P'$ is a locally Lipschitz reduction.
\end{proposition}

\begin{proof}
By \cref{lem:primitive}, choose one primitive weight per irreducible components and denote the resulting set by $\Xi$. The span of $\rho(\gamma)V_\xi$ over $\gamma\in\Gamma$ and $\xi\in\Xi$ is $G$-invariant by Borel density (\cref{Borel}). For each irreducible summand $U$, the span contains the nonzero space $U\cap V_\xi$ for its chosen $\xi$, hence all of $U$. Thus we can choose a basis $v_1,\ldots,v_d$ of $V$ with $v_j\in\rho(\gamma_j)V_{\xi_j}$.

In a local H\"older frame $b$ of $P'$, the fields $X_j(z)=b_zv_j$ form a continuous frame with $X_j(z)\in F_z^{\gamma_j,\xi_j}$. Let $s:M\to \GL(V)/C'$ represent $P'$ in a smooth local frame. By \cref{cor:Lipalongprim,lem:primtrans}, the change of $s$ along a sufficiently short segment tangent to any $X_j$ is bounded in $d_{C'}$ by $K$ times its $h$-length, with one $K$ for this finite family. For nearby $x,y$, apply \cref{joining} in the same chart. Summing over the resulting segments $c_i$ gives \[ d_{C'}(s(x),s(y))\le K\sum_i\ell_h(c_i)\le KL d_h(x,y), \] which proves local Lipschitz regularity.
\end{proof}

\begin{remark}
It does not assert that every choice of adapted frame is locally Lipschitz.  After this conclusion, smooth local sections of $\textrm{GL}(V)\to\textrm{GL}(V)/C'$ provide locally Lipschitz adapted frames. An arbitrary point dependent $C'$-change of such a frame can have lower regularity.
\end{remark}

\subsection{The weak derivative of the reduction}

Fix a compact subgroup $C'\subset C$ and a H\"older $C'$-reduction $P'\subset P$ satisfying \eqref{eq:equiv'}. By \cref{lipP}, the section $x\mapsto P'_x$ of $\Fr(TM)/C'$ is locally Lipschitz.

Fix a smooth chart $q=(q^1,\ldots,q^d):\Ocal\to U\subset\Rbb^d$. Write $e_x=(dq_x)^{-1}:V\to T_xM$ for its coordinate frame and the local Lipschitz map $\sigma:\Ocal\to\GL(V)/C'$ for the map representing $P'$ in this frame. If $p_{c'}:\GL(V)\to\GL(V)/C'$ is the quotient map, choose a smooth local section $r$ of $p_{c'}$ near $\sigma(x_0)$ for $x_0\in\Ocal$, after making $\Ocal$ smaller if necessary, so that \[ B=r\circ\sigma,\qquad b_x=e_x\circ B(x)\in P'_x,\qquad\sigma(x)=p_{c'}(B(x)). \]
For the standard basis $(\varepsilon_i)$ of $V$, put $E_i(x)=e_x\varepsilon_i=\partial/\partial q^i|_x$. The $a$th column of $B$ gives the coefficients of the vector field $b\varepsilon_a=\sum_i B^i{}_a E_i$ on $\Ocal$.

 In coordinates put $\widetilde B=B\circ q^{-1}$. Its entries have weak first derivatives in $L^\infty_{\mathrm{loc}}(U)$, characterized by \[ \int_U \widetilde B^i{}_a\,\partial_j\varphi\,dy  =-\int_U (\partial_j\widetilde B^i{}_a)\varphi\,dy\] for all $\varphi\in C_c^\infty(U)$. Since both $B$ and $B^{-1}$ are locally Lipschitz, weak derivatives agree almost everywhere with the classical derivatives by Rademacher's theorem. The smooth positive density $\mu$ has the same null sets as Lebesgue measure in this chart. Thus, for almost every $x$, $dB_x(Z)=D\widetilde B_{q(x)}(dq_x Z)$ is defined for every $Z\in T_xM$: When $Z=\sum_{j=1}^d Z^j\,\left.\frac{\partial}{\partial q^j}\right|_x$ and $dq_x(Z)=(Z^1,\ldots,Z^d)$, $dB_x(Z)$ is defined as
\[dB_x(Z)=\sum_{j=1}^d Z^j\,\partial_j\widetilde B(y)\] where the right hand side can be understood via weak derivatives (and same with classical derivatives almost everywhere).

Let $\omega$ denote the connection form in $e$ from \eqref{connd+A}, so $\nabla=d+\omega$. For a fixed
$v\in V$, the product rule applied to the local vector field $bv=e(Bv)$ gives, almost everywhere,
\[ (\nabla_Z(bv))_x=e_x\big((dB_x(Z)+\omega_x(Z)B(x))v\big).\]
Again, even though $b$ is only locally Lipschitz, $dB_x(Z)$ is the only term taking the derivatives of $b$. We already discussed that $dB_x(Z)$ is well-defined almost everywhere, hence $\nabla_Z(bv)$ is defined almost everywhere and locally essentially bounded.

In the frame $b$ this identity can be written
\[ \Omega_b^\nabla\coloneqq B^{-1}(dB+\omega B),\qquad \Omega_b^\nabla(x)(Z)v=b_x^{-1}(\nabla_Z(bv))_x .\] Especially, \[\Omega_b^\nabla(x)(Z)\in\mathfrak{gl}(V)\] is well-defined almost everywhere. Let us define
\begin{equation}\label{eq5111}
 W'=\Hom_{\Rbb}(V,\mathfrak{gl}(V)/\clie'),\qquad  a_b(x)(u)=\big[\Omega_b^\nabla(x)(b_xu)\big]_{\clie'}
\end{equation}
for $u\in V$. The brackets denote a class in a quotient of real vector spaces. Hence, $a_b(x)\in W'$ is well-defined almost every $x$ and (the coefficients of) $a_b$ are locally essentially bounded, i.e., $a_b\in L^\infty_{\mathrm{loc}}(\Ocal,W')$.

Now, let us see what happened if we change the adapted frame: On $W'$, we can define the natural actions of $G$ and $C'$ as
\[ (\tau_G(g)T)(u)=\textrm{Ad}_{\rho(g)} T(\rho(g)^{-1}u),\qquad (\tau_C(c)T)(u)=\textrm{Ad}_{c}T(c^{-1}u) \] for $T\in W'$.
Both conjugations preserve $\clie'$, and the actions commute because $C'$ centralizes $\rho(G)$. If $\widehat b=bc$ is another locally Lipschitz adapted frame, then $c(x)=b_x^{-1}\widehat b_x\in C'$ is locally Lipschitz. The weak product rule gives
\[ \Omega_{\widehat b}^\nabla=c^{-1}\Omega_b^\nabla c+c^{-1}dc.\]
At a differentiability point, $dc_x(Z)\in T_{c(x)}C'$, so
$c(x)^{-1}dc_x(Z)\in\clie'$. Passing to the quotient and
using $\widehat b_xu=b_xc(x)u$ for $u\in V$ gives
\begin{align}
 a_{\widehat b}(x)(u)&=\textrm{Ad}_{c(x)^{-1}}\big(a_b(x)(c(x)u)\big),\label{actionC}\\
 a_{\widehat b}(x)&=\tau_C(c(x)^{-1})a_b(x).\label{eq5112}
\end{align}
Thus $a(x)=[b_x,a_b(x)]$ is a measurable section of the associated bundle $P'\times_{C'}W'$. Here, $P'\times_{C'}W'=(P'\times W')/\!\sim$ where $(b,T)\sim(bc,\tau_C(c^{-1})T)$. After averaging an inner product on $W'$ over $C'$, the induced bundle norm $\|\cdot\|$ satisfies $\|a(x)\|=\|a_b(x)\|$ does not depend on the choice of $b$ almost everywhere. A finite cover by relatively compact frame neighborhoods therefore gives a global essential bound $\|a(x)\|\le K$ almost every $x$ and any choice of $b$.

\begin{proposition}\label{derivP}
Every H\"older reduction $P'$ as above is locally $C^{1+\beta}$ for some $0<\beta<1$. For every local $C^{1+\beta}$ section $b:\Ocal\to P'$, the coefficient $a_b$ of \eqref{eq5111} is classically defined and takes values in $(W')^G$ at every point of $\Ocal$, i.e., 
if we write \[a_b(x)(u)=\left[b_x^{-1}(\nabla_{b_xu}b)_x\right]_{\clie'}, \quad x\in\mathcal O,\quad u\in V,\]
then $a_b(x)\in(W')^G$ at every point, where
\[(g\cdot A)(u)=\textrm{Ad}_{\rho(g)}\left(A(\rho(g)^{-1}u)\right).\]

\end{proposition}

\begin{proof}[Proof of \cref{derivP}]

Fix $\gamma\in\Gamma$, set $f=\alpha(\gamma)$ and $R=\rho(\gamma)$, and choose locally Lipschitz adapted frames $b$ near $x$ and $b'$ near $f(x)$. By \eqref{eq:equiv'}, the frame
\[\wtd b_{f(x)}=D_xf\,b_xR^{-1}=b'_{f(x)}k(x)\]
is adapted, with $k(x)\in C'$. Both $\widetilde b$ and $k$ are locally Lipschitz because $f$ is $C^2$. For fixed $v\in V$ the target vector field $\widetilde bv$ is $f_*(bR^{-1}v)$.
Invariance of $\nabla$ and the weak chain rule in smooth coordinates give, for $Z\in T_xM$ and almost every $x$, 
\[ \Omega_{\widetilde b}^\nabla(f(x))(D_xf Z)v=R\,\Omega_b^\nabla(x)(Z)R^{-1}v .\]
Set $Z=b_xR^{-1}u$, so that $D_xfZ=\widetilde b_{f(x)}u$.
Taking classes and then applying \eqref{eq5112} at $f(x)$ gives
\[ \tau_G(\gamma)a_b(x) =a_{\widetilde b}(f(x)) =\tau_C(k(x)^{-1})a_{b'}(f(x)).\]
It follows that
\begin{equation}\label{eq:Pcov} a_{b'}(f(x))=\tau_C(k(x))\tau_G(\gamma)a_b(x) \qquad\text{for almost every }x .\end{equation}

Use the finite frame cover from the preceding subsection. Since $\Gamma$ is countable and preserves $\mu$, there is one invariant conull set on which the bound and all the frame-change and covariance identities for this cover hold along each orbit. For $x$ in this set and any $\gamma\in\Gamma$, choose a frame from the cover at $f(x)$. Isometricity of $\tau_C(k(x))$ gives
\[ \|\tau_G(\gamma)a_b(x)\|=\|a_{b'}(f(x))\|\le K .\]
By \cref{lem:bdd}, a vector with bounded lattice orbit in the finite-dimensional $G$-module $W'$ is fixed by $G$. Thus $a$ takes values almost everywhere in $E=P'\times_{C'}(W')^G$.

Now take $f=\alpha(\gamma_0)$ to be the given Anosov diffeomorphism. The map $p_x\mapsto D_xf\,p_x\rho(\gamma_0)^{-1}$ on $P'$ commutes with the right $C'$-action. It induces on $E$ the bundle map
\[ F_x[p_x,T]=[D_xf\,p_x\rho(\gamma_0)^{-1},T].\]
 In local frames its matrix is $\tau_C(k(x))|_{(W')^G}$. Thus $F$ is H\"older and isometric for the metric induced by the $C'$-invariant
inner product. \eqref{eq:Pcov} reduces to
$a(f(x))=F_xa(x)$ almost everywhere. Hence, we can consider $a$ as a measurable solution of the H\"older isometric cocycle $F_x$ of the H\"older bundle $E$. We use \cref{prop:sectionconf} so that there exists a H\"older representative of $a$, and write $a_b$ for its coordinates in the chosen lift.

Recall that  $p_{c'}:\GL(V)\to\GL(V)/C'$ be the quotient map, $p_{c'}(B)=BC'$. In the chosen coordinate frame, the reduction is represented by $\sigma=p_{c'}\circ B$.

Fix a point $x$ at which $B$ is differentiable and the defining identity for $a_b$ holds. Such points have full measure. For
$Z\in T_xM$, the chain rule gives
\[ d\sigma_x(Z)=(dp_{c'})_{B(x)}\bigl(dB_x(Z)\bigr).\]
The connection formula in the frame $b=eB$ gives
\[ dB_x(Z)=B(x)\,\Omega_b^\nabla(x)(Z)-\omega_x(Z)\,B(x).\]
We compute the image of each term under $(dp_{c'})_{B(x)}$.

For the first term, write $B=B(x)$ and $A=\Omega_b^\nabla(x)(Z)$. The matrix $BA$ is the velocity at $t=0$ of the curve $B\exp(tA)$. Thus
\[ (dp_{c'})_B(BA)=\left.\frac{d}{dt}\right|_{t=0}B\exp(tA)C'.\]
For $C_0\in\clie'$, the curve $B\exp(tC_0)C'$ is constant, since $\exp(tC_0)\in C'$. Hence $(dp_{c'})_B(BC_0)=0$. By linearity, $(dp_{c'})_B(BA)$ depends only on the class $[A]_{\clie'}$. We therefore write
\[ \iota_B([A]_{\clie'})\coloneqq (dp_{c'})_B(BA)\in T_{BC'}\bigl(\GL(V)/C'\bigr).\]

To identify the class of $A=\Omega_b^\nabla(x)(Z)$, recall that
\[a_b(x)(u)=\bigl[\Omega_b^\nabla(x)(b_xu)\bigr]_{\clie'},\qquad u\in V.\]
Taking $u=b_x^{-1}Z$ gives $b_xu=Z$, and hence
\[(dp_{c'})_B\bigl(B\,\Omega_b^\nabla(x)(Z)\bigr)=\iota_B\bigl(a_b(x)(b_x^{-1}Z)\bigr).\]
Here $b_x^{-1}Z=B(x)^{-1}dq_x(Z)$ is the column of coordinates of $Z$ in the adapted frame $b_x$.

For the second term, put $W=\omega_x(Z)$. The matrix $WB$ is the velocity at $t=0$ of $\exp(tW)B$. Denote its projection by
\[W^\#_{BC'}\coloneqq \left.\frac{d}{dt}\right|_{t=0}\exp(tW)BC'=(dp_{c'})_B(WB).\]
Substituting these two expressions into the chain rule gives
\begin{equation}\label{eq:Pquot} d\sigma_x(Z)=\iota_{B(x)}\bigl(a_b(x)(b_x^{-1}Z)\bigr)-(\omega_x(Z))^\#_{\sigma(x)}.\end{equation}
This identity holds for every $Z\in T_xM$ at almost every $x$.

In smooth local coordinates, the right-hand side of \eqref{eq:Pquot} is $C^\beta$ for some $\beta>0$: 
Indeed, the two terms in \eqref{eq:Pquot} are obtained by applying $(dp_{c'})_{B(x)}$ to matrix products. In local coordinates the entries of $(dp_{c'})_B$ are smooth functions of $B$, and a representative of each class in $\mathfrak{gl}(V)/\clie'$ can be chosen linearly. Thus, for each coordinate vector field $Z$, substituting the locally Lipschitz matrix $B(x)$ and the H\"older coefficients $a_b(x)(b_x^{-1}Z)$ and $\omega_x(Z)$
gives a H\"older right-hand side.

Since $\sigma$ is locally Lipschitz, \eqref{eq:Pquot} therefore gives $C^\beta$ representatives of its weak first derivatives. Continuous weak first derivatives are classical, so $\sigma$ is $C^{1+\beta}$. Composing $\sigma$ with the smooth local section $r$ of $p_{c'}$ gives a $C^{1+\beta}$ lift $B=r\circ\sigma$ and hence a $C^{1+\beta}$ adapted frame $b=eB$.

In any $C^{1+\beta}$ adapted frame, \eqref{eq5111} is now classical and continuous. By \eqref{actionC} its values lie in the closed subspace $(W')^G$ almost everywhere, hence everywhere by full support of $\mu$.

\end{proof}

\subsection{Reduction to one compact stabilizer}

We explain first the purpose of this step. The preceding regularity
argument gives local $C^{1+\beta}$ adapted frames $b$. In such a frame,
the connection $b^{-1}\nabla b$ has a coefficient with values in
$\mathfrak{gl}(V)$. Its class $[b^{-1}\nabla b]$ modulo the Lie algebra of the compact
structure group is independent of the derivatives of the choice of adapted frame. \Cref{derivP} says that $[b^{-1}\nabla b]$ is fixed by $G$. It does not, however, say that coordinates of $[b^{-1}\nabla b]$  are constant when we change the base point or the adapted frame, see \eqref{actionC}.

We will choose a smaller compact reduction on which $[b^{-1}\nabla b]$ is constant. The smaller compact structure group will also fix $[b^{-1}\nabla b]$.
These two facts permit us to subtract a globally defined tensor from
$\nabla$. The resulting connection $\wtd\nabla$ will preserve the smaller reduction
and will have parallel torsion, see \cref{prop:newconnred}.

Before stating the main statement in this section recall what we have shown so far: We are in \cref{setting}. Let $V=\Rbb^d$ be the vector space of the model. Combining previous discussions, there exist 

\begin{enumerate}
\item a continuous representation with no zero weight $\rho:G\to \GL(V)$;
\item  a closed subgroup $C\subset \SL(V)\cap \mathcal{Z}_{\SL(V)}(\rho(G))$ ;
\item a $C^{1+\beta}$ principal $C$-reduction $P\subset \Fr(TM)$ such that:
\[D_x\alpha(\gamma)(P)_x
       =(P)_{\alpha(\gamma)x}\rho(\gamma)\] for every
       $\gamma\in\Gamma$ and $x\in M$.
;
\item any $C'$-reduction $P'\subset P$ with a compact subgroup $C'\subset C$ satisfies $D\alpha(\gamma)P'=P'\rho$ is $C^{1+\beta}$;
\item a H\"older metric $h$ on $TM$ such that \cref{prop:metric} holds;
\item a H\"older $\alpha(\Gamma)$-invariant connection $\nabla$ on $M$; 
\end{enumerate}

\begin{theorem}
\label{prop:reduction1}
Under the preceding hypotheses, there exist a closed subgroup
$C_*\subset C$, a $C^{1+\beta}$ principal $C_*$-reduction
$P_*\subset P$ for some $\beta>0$, and an element
\[
 a_*\in W_*:=\Hom_{\Rbb}
       \left(V,\mathfrak{gl}(V)/\clie_*\right),
 \qquad \clie_*:=\textrm{Lie}(C_*),
\]
such that:
\begin{enumerate}
\item $D_x\alpha(\gamma)(P_*)_x
       =(P_*)_{\alpha(\gamma)x}\rho(\gamma)$ for every
       $\gamma\in\Gamma$ and $x\in M$;
\item $a_*$ is fixed by $G$ and $C_*$, acting by
\[
 (g\cdot T)(u)=\textrm{Ad}_{\rho(g)}
                    \left(T(\rho(g)^{-1}u)\right),\qquad
 (c\cdot T)(u)=\textrm{Ad}_{c}\left(T(c^{-1}u)\right);
\]
\item for every local $C^{1+\beta}$ section $b$ of $P_*$,
\[
 \left[b_x^{-1}(\nabla_{b_xu}b)_x\right]_{\clie_*}
       =a_*(u)
 \qquad \forall x\text{ in the domain of }b,\ \forall u\in V.
\]
\end{enumerate}
Here $[A]_{\clie_*}$ means the class of the endomorphism
$A\in\mathfrak{gl}(V)$ in the vector-space $\mathfrak{gl}(V)/\clie_*$. 
\end{theorem}
\begin{remark}
    Indeed, for the proof of \cref{prop:reduction1}, we only need to assume that there is no trivial representation component of $\rho$, i.e., $V$ is a direct sum of non-trivial $G$-irreducible representations, instead of no zero restricted weight. 
\end{remark}

\begin{proof}

For $K \subset C$ a compact subgroup of $C$, we denote $W_K = \Hom(V,\mathfrak{gl}(V)/\text{Lie}(K))$.

We first use \cref{derivP} on the $C$-reduction $P$, which is allowed by \cref{thm:reduction}.
For $p \in P_x$, let $b$ be a local $C^{1+\beta}$ section of $P$ around $x$ such that $b_x = p$ and define
\[
a_C(p) =a_b(x),
\]
where $a_b$ is defined in \cref{derivP}.
The above is well defined since it does not depend on the choice of $b$.
Fix $\gamma \in \Gamma$ and $p \in P_x$.
Write $q = D\alpha(\gamma)p\rho(\gamma)^{-1} \in P_{\alpha(\gamma)x}$.
Choose a local section $b$ of $P$ around $x$ satisfying $b_x = p$.
By equivariance of $P$ under the $\Gamma$ action, the expression
\[
\tilde{b}_{\alpha(\gamma)(y)} = D_y\alpha(\gamma) b_y \rho(\gamma)^{-1}
\]
defines a local section of $P$ nearby $\alpha(\gamma)(x)$ which satisfies $\tilde{b}_{\alpha(\gamma)(x)} = q$.
For $u,v \in V$, we have the vector $qu = D_x\alpha(\gamma)p\rho(\gamma)^{-1}$ and the vector field $\tilde{b}v$ on a neighborhood of $\alpha(\gamma)(x)$ defined by $z \mapsto \tilde{b}_zv$.
By invariance of $\nabla$,
\[
(\nabla_{qu}(\tilde{b}v))_{\alpha(\gamma)(x)} = D_x\alpha(\gamma)((\nabla_{p\rho(\gamma)^{-1}u}(b(\rho(\gamma)^{-1}v))_x).
\]
Applying $q^{-1}$ above and using that $q^{-1}D_x\alpha(\gamma) = \rho(\gamma)p^{-1}$, we obtain the equality of endomorphisms
\[
q^{-1}(\nabla_{qu}\tilde{b})_{\alpha(\gamma)(x)}
= \rho(\gamma)p^{-1}((\nabla_{p\rho(\gamma)^{-1}u}b)_x)\rho(\gamma)^{-1}.
\]
Considering the above in the quotient $\mathfrak{gl}(V)/\mathfrak{c}$, while using that $\rho(\gamma)$ centralizes $\mathfrak{c}$,
\[
a_C(D_x\alpha(\gamma)p\rho(\gamma)^{-1})(u) =
Ad_{\rho(\gamma)}(a_C(p)(\rho(\gamma)^{-1}u)).
\]
But by \cref{derivP},
\[
Ad_{\rho(\gamma)}(a_C(p)(\rho(\gamma)^{-1}u)) = a_C(p)(u).
\]
Therefore,
\begin{equation}
    \label{eq:aCergodic}
    a_C(p) = a_C(D_x\alpha(\gamma)p\rho(\gamma)^{-1}).
\end{equation}

Now fix $p_0 \in P$ and set $w = a_C(p_0)$.
Consider a $C$-invariant norm on $W_C$ by averaging an inner product along its orbit of $C$.
Define the function
\[
x \mapsto d(a_C(p),C \cdot w),
\]
where $p \in P_x$.
It is well defined since such a choice of $p$ is unique up to an element of $C$ and the distance is $C$-invariant.
Moreover, it is continuous since we can choose $p$ continuously and $\Gamma$-invariant by \eqref{eq:aCergodic}.
By ergodicity, it must be constant almost everywhere and, by continuity and full support of the measure, it must be constant.
Since it vanishes at the basepoint of $p_0$, it vanishes.
By compactness of $C$, we must have $a_C(p) \in C\cdot w$ for all $p \in P$.

Recall that $C$ acts on $W_C$ by
\[
(c\cdot A)(u) = Ad_c(A(c^{-1}u)).
\]
Moreover, for $c \in C$, $p \in P$ at the basepoint $x$, $u \in V$ and $b$ local section of $P$ with $b_x = p$,
\begin{align*}
    a_C(pc)(u) =& [(pc)^{-1}(\nabla_{(pc)u}(bc))_x]_{\mathfrak{c}} \\
    =& [c^{-1}p^{-1}(\nabla_{p(cu)}b)_xc]_{\mathfrak{c}} \\
    =& Ad_{c^{-1}}(a_C(p)(cu)).
\end{align*}
Therefore,
\begin{equation}
\label{aCequivariant}
    a_C(pc)= c^{-1}\cdot a_C(p)
\end{equation}
Let $C_*=\operatorname{Stab}_C(w)$ and $P_*=a_C^{-1}(w)$.
Since $a_C(P_x)=C\cdot w$, each fiber $(P_*)_x$ is nonempty.
Additionally, for $p\in(P_*)_x$, every frame in $P_x$ has the form $pc$, and
\eqref{aCequivariant} gives
\[
    pc\in(P_*)_x \iff c^{-1}\cdot w=w \iff c\in C_*.
\]
Consequently, $(P_*)_x=pC_*$, making $P_*$ a $C_*$-reduction.

In order to verify it is H\"older, let $b$ be a local H\"older section of $P$ around some point $x_0$, we wish to find a local H\"older section $c$ of $C$ which satisfies $bc \in P_*$.
By \eqref{aCequivariant}, this is equivalent to solving
\begin{equation}
\label{eq:solveforcx}
a_C(b_x) = c_x \cdot w.
\end{equation}
Let $s \colon C\cdot w\to C$ be a local smooth section of $C \to C \cdot w \cong C/C_*$ in a neighborhood of $a_C(b_{x_0})$.
Taking a smaller neighborhood of $x_0$ if necessary, we define $c_x = s(a_C(b_x))$, which satisfies \eqref{eq:solveforcx}.
Since $s$ is smooth and $b$ is H\"older, then $c$ is H\"older, which shows that the $C_*$-reduction is H\"older.

We now wish to apply \cref{derivP} to the $C_*$-reduction $P_*$.
To do so, we must verify \eqref{eq:equiv'}.
This follows from the fact that $p \in P_*$ if, and only if, $a_C(p) = w$ together with \eqref{eq:aCergodic}.
In particular, this shows the first property of the \cref{prop:reduction1}.

Now \cref{derivP} implies that the $C_*$ reduction is $C^{1+\beta}$ and constructs $a_{C_*} \colon P_* \to W_{C_*}^G$ which we will show is constant.
Let $p \in (P_*)_x$ and $b$ a local $C^{1+\beta}$ section of $P_*$ satisfying $b_x = p$.
It is also a section of $P$, which implies for $u \in V$
\[
[a_{C_*}(p)(u)]_{\mathfrak{c}} = [[b_x^{-1}(\nabla_{b_xu}b)]_{\mathfrak{c}_*}]_{\mathfrak{c}} = a_{C}(p)(u) = w(u).
\]
Denoting by $\pi \colon W_{C_*} \to W_C$ the natural projection, we just showed $\pi(a_{C_*}(p)) = w$ for any $p \in P_*$.
In particular, for $p,q \in P_*$,
\[
a_{C_*}(p)-a_{C_*}(q) \in \ker \pi = \Hom(V,\mathfrak{c}/\mathfrak{c}_*).
\]
However, \cref{derivP} also implies that $G$ fixes $A_{p,q} = a_{C_*}(p)-a_{C_*}(q)$.
This means that for any $u \in V$ and $g \in G$,
\[
A_{p,q}(u) = Ad_{\rho(g)}A(\rho(g)^{-1}u) = A(\rho(g)^{-1}u),
\]
since $\rho(G)$ commutes with $\mathfrak{c}$.
But $V$ is a direct sum of non-trivial $G$-irreducible representations, which implies that $A_{p,q} = 0$ and that $a_{C_*}$ is equal to a constant value which we denote by $a_*$.
It satisfies the third property of \cref{prop:reduction1} by definition.

\cref{derivP} implies that $a_*$ is fixed by $G$.
For $c \in C_*$, we use \eqref{aCequivariant} which also holds for $C_*$ to obtain
\[
c^{-1} \cdot a_* = c^{-1} \cdot a_{C_*}(p) = a_{C_*}(pc) = a_*.
\]
This shows that $a_*$ is $C_*$-invariant, verifying the second condition of \eqref{prop:reduction1} and finishes its proof.
\end{proof}

\section{A H\"older connection $\wtd\nabla$ and a nilpotent Lie algebra}\label{sec6}
Let us summary what we had so far: We are in \cref{setting}. Let $V=\Rbb^d$ be the model vector space. Under the preceding hypotheses, there exist 

\begin{enumerate}
\item a continuous representation with no zero weight $\rho:G\to \GL(V)$;
\item  a closed subgroup $C_* \subset \GL(V)$ with  $\clie_*:=\textrm{Lie}(C_*)$;
\item a $C^{1+\beta}$ principal $C_*$-reduction $P_*\subset \Fr(TM)$ for some $\beta>0$;
\item a H\"older metric $h$ on $TM$ such that \cref{prop:metric} holds;
\item a H\"older $\alpha(\Gamma)$-invariant connection $\nabla$ on $M$; 
\item an element
\[
 a_*\in W_*:=\Hom_{\Rbb}
       \left(V,\mathfrak{gl}(V)/\clie_*\right),
\]
\end{enumerate}
such that:
\begin{enumerate}
\item $D_x\alpha(\gamma)(P_*)_x
       =(P_*)_{\alpha(\gamma)x}\rho(\gamma)$ for every
       $\gamma\in\Gamma$ and $x\in M$;
\item $a_*$ is fixed by $G$ and $C_*$, acting by
\[
 (g\cdot T)(u)=\textrm{Ad}_{\rho(g)}
                    \left(T(\rho(g)^{-1}u)\right),\qquad
 (c\cdot T)(u)=\textrm{Ad}_{c}\left(T(c^{-1}u)\right);
\]
\item for every local $C^{1+\beta}$ section $b$ of $P_*$,
\[
 \left[b_x^{-1}(\nabla_{b_xu}b)_x\right]_{\clie_*}
       =a_*(u)
 \qquad(x\text{ in the domain of }b,\ u\in V).
\]
\end{enumerate}
Here $[A]_{\clie_*}$ means the class of the endomorphism
$A\in\mathfrak{gl}(V)$ in the vector-space quotient by
$\clie_*$.

In this section, we will show \cref{sec6conclusion}. 
The torsion-free connection $\nabla$ was useful for improving regularity, but its parallel transport need not preserve the final compact reduction $P_*$. $a_*$ is exactly the defect, modulo $\clie_*$. We lift $a_*$ to a constant map $A_0:V\to\mathfrak{gl}(V)$ and subtract the corresponding tensor $S$. The resulting connection, denoted $\wtd\nabla$, will preserve $P_*$, \cref{prop:newconnred}.

Its torsion will be the tensor whose coefficients in adapted frames are \[T_0(u,v)=-A_0(u)v+A_0(v)u.\] Because these coefficients are constant and fixed by $C_*$, and
because $\wtd\nabla$ preserves $P_*$, the torsion will be parallel. Later, we prove that $-T_0$ satisfies the Jacobi identity and is nilpotent.

\subsection{The unique equivariant lift of the constant coefficient}

\begin{lemma}[The lift $A_0$]
\label{lem:unique-A0-lift}
There exists a unique $A_0$ lift of $a_*$
\[
\begin{tikzcd}
                                             & \mathfrak{gl}(V) \arrow[d]      \\
V \arrow[r, "a_*"] \arrow[ru, "A_0", dashed] & \mathfrak{gl}(V)/\mathfrak{c}_*
\end{tikzcd}
\]
satisfying for every $g \in G$ and $u,v \in V$,
\begin{equation}
\label{eq:G-equiv}
	A_0(\rho(g)u)\rho(g)v = \rho(g)A_0(u)v.
\end{equation}
Additionally, it satisfies for every $c \in C_*$ and $u,v \in V$,
\begin{equation}
\label{eq:C-equiv}
	A_0(cu)cv = cA_0(u)v.
\end{equation}
\end{lemma}

\begin{proof}
We define a $G$-module structure on $\Hom_{\mathbb{R}}(V,W)$ for $W = \mathfrak{c}_*, \mathfrak{gl}(V)$ and $\mathfrak{gl}(V)/\mathfrak{c}_*$ by defining the $G$ action
\begin{equation*}
	(g \cdot A)(u) = \rho(g) A(\rho(g)^{-1}u) \rho(g)^{-1}.
\end{equation*}
Note that $A \in \Hom(V,\mathfrak{gl}(V))$ is a fixed element for the $G$ action if, and only if, it satisfies \eqref{eq:G-equiv}.

Consider the short exact sequence of $G$-modules
\[
 0\rightarrow\Hom(V,\clie_*)
 \rightarrow\Hom(V,\mathfrak{gl}(V))
 \xrightarrow{p} \Hom(V,\mathfrak{gl}(V)/\mathfrak{c}_*)\rightarrow0.
\]
Since $G$ is semisimple, the short exact sequence splits and there exists $q \colon \Hom(V,\mathfrak{gl}(V)/\mathfrak{c}_*) \to \Hom(V,\mathfrak{gl}(V))$ such that $p \circ q = \Id_{\Hom(V,\mathfrak{gl}(V)/\mathfrak{c}_*)}$.
Define $A_0 = q(a_*)$.

Since $\text{Im}(q)$ is $G$-invariant
 \[
 g \cdot A_0 \in \text{Im}(q), \ \text{for all} \ g \in G.
 \]
On the other hand,
\[
p(g \cdot A_0) = g \cdot p(A_0) = g \cdot a_* = a_*, \ \text{for all} \ g \in G.
\]
But $A_0$ is the unique element in $\text{Im}(q)$ that projects to $a_*$.
Therefore, $g \cdot A_0 = A_0$, showing that $A_0$ is a fixed element of the $G$ action.

Any other lift of $a_*$ is of the form $A_0 + B$ where $B \in \Hom(V,\mathfrak{c}_*)$.
If it additionally fixed by the $G$ action, $B$ must also be fixed by the $G$ action.
However, since $\mathfrak{c}_*$ commutes with $\rho(G)$, we must have
\[
B(u) = (g \cdot B)(u) = \rho(g)B(\rho(g)^{-1}u)\rho(g)^{-1} = B(\rho(g)^{-1}u), \ \text{for all} \ g \in G \ \text{and}  \ u \in V.
\]
For $U \subset V$ an irreducible $G$-space, the kernel of $B|_U$ is invariant, making it either trivial of equal to $U$.
In the trivial case, $B|_U$ is injective.
Since $B(\rho(g)^{-1}u) =B(u)$ for all $g \in G$ and $u \in U$, injectivity implies $\rho(g)^{-1}u = u$, making $U$ a trivial irreducible $G$ space, contradicting our hypothesis.
Therefore, $B$ vanishes on every irreducible $G$-space, implying that $B$ vanishes and showing uniqueness of the $G$-equivariant lift.

For $c\in C_*$, let
\[
 (c\cdot A)(u)=cA(c^{-1}u)c^{-1}.
\]
Since $C_*$ commutes with $\rho(G)$, this action commutes with that of $G$.
In particular $c \cdot A_0$ is fixed by $G$ and projects to $a_*$.
By uniqueness, $c \cdot A_0 = A_0$, for all $c \in C_*$.
\end{proof}

\subsection{Constructing the tensor and the modified connection}

\begin{proposition}[A connection preserving the final reduction]
\label{prop:newconnred}
For $x\in M$, $p\in(P_*)_x$ and $u,v\in V$, set
\begin{equation}
 S_x^{(p)}(pu,pv):=pA_0(u)v.
 \label{eq:S-from-A0}
\end{equation}
Then $S_x^{(p)} = S_x$ does not depend on $p$ and $x\mapsto S_x$ defines a $C^{1+\beta}$ section of $\Hom(TM \otimes TM,TM)$ satisfying for all $x \in M$, $X,Y \in T_xM$ and $\gamma \in \Gamma$,
\[
S_{\alpha(\gamma)x}(D_x\alpha(\gamma)X,D_x\alpha(\gamma)Y) = D_x\alpha(\gamma)S_x(X,Y).
\]
Moreover, the H\"older connection defined by
\begin{equation}
 \wtd\nabla_XY:=\nabla_XY-S(X,Y)
 \label{eq:newconn}
\end{equation}
is $\Gamma$-invariant and $P_*, S$ and its torsion are parallel with respect to it.
\end{proposition}

\begin{proof}
If $c \in \mathfrak{c}_*$, $x \in M$ and $u,v \in V$, then by $C$-equivariance of $A_0$ \eqref{eq:C-equiv}
\[
S^{(pc)}_x(pcu,pcv) = pcA_0(u)v = pA_0(cu)cv =S^{(p)}_x(pcu,pcv).
\]
Therefore, $S^{(p)}_x$ does not depend on the choice of $p \in (P_*)_x$ and we may denote it by $S_x$.
For $\gamma \in \Gamma$ and $X,Y \in T_xM$, choose $p \in (P_*)_x$ and $u,v \in V$ such that $pu = X$ and $pv = Y$.
Since $D_x\alpha(\gamma) p = p' \rho(\gamma)$ for some $p' \in (P_*)_{\alpha(\gamma)x}$, 
\begin{align*}
S_{\alpha(\gamma)x}^{(p')}(D_x\alpha(\gamma)(X),D_x\alpha(\gamma)(Y)))=&   
S_{\alpha(\gamma)x}^{(p')}(D_x\alpha(\gamma)(pu),D_x\alpha(\gamma)(pv))\\
=&  S_{\alpha(\gamma)x}^{(p')}((D_x\alpha(\gamma)p\rho(\gamma)^{-1})\rho(\gamma)u,(D_x\alpha(\gamma)p\rho(\gamma)^{-1})\rho(\gamma)v)\\
=& D_x\alpha(\gamma)p \rho(\gamma)^{-1}A_0(\rho(\gamma)u)\rho(\gamma)v  \\
=& D_x\alpha(\gamma)pA_0(u)v \\
=& D_x\alpha(\gamma)S^{(p)}_x(X,Y).
\end{align*}
This shows $\Gamma$-equivariance of $S$.
Since $\nabla$ is also $\Gamma$-equivariant, then so is $\tilde{\nabla}$.

Let $b$ be a local section of $P_*$ defined on a small neighborhood and $c \colon [0,1] \to M$ a $C^1$ curve mapping to that neighborhood.
Denoting by $\tilde{T}_{c|_{[0,t]}} \colon T_{c(0)}M \to T_{c(t)}M$ the parallel transport for $\tilde{\nabla}$, we will show
\[
	\tilde{T}_{c|_{[0,t]}} b_{c(0)} \in (P_*)_{c(t)}.
\]
Indeed, let $B(t) \in \GL(V)$ be such that
\[
	b_{c(t)} B(t) = \tilde{T}_{c|_{[0,t]}} b_{c(0)}.
\]
Note that $B(0) = \text{Id}$ and we must show $B(t) \in C_*$ for all $t$.
Let $v(t) \in V$ be such that $c'(t) = b_{c(t)}(v(t))$.
Since $b_{c(t)}B(t)$ is being parallel transported along $c(t)$,
\begin{align*}
0 = \tilde{\nabla}_{c'(t)}b_{c(t)}B(t) =& b_{c(t)}B'(t) + (\tilde{\nabla}_{c'(t)}b_{c(t)})B(t)\\
\end{align*}
By construction,
\[
[b_{c(t)}^{-1}\tilde{\nabla}_{c'(t)}b_{c(t)}]_{\mathfrak{c}_*} = [b_{c(t)}^{-1}\nabla_{b_{c(t)}v(t)}b_{c(t)} - b_{c(t)}^{-1}S(b_{c(t)}v(t),b_{c(t)}\cdot)]_{\mathfrak{c}_*} = 0,
\]
meaning 
\begin{equation}
\label{eq:Bincstar}
	b_{c(t)}^{-1}\tilde{\nabla}_{c'(t)}b_{c(t)} \in \mathfrak{c}_* \ \text{for all} \ t.
\end{equation}
Since $B(t)$ satisfies the ODE
\[
B'(t) =- b_{c(t)}^{-1}\tilde{\nabla}_{c'(t)}b_{c(t)} B(t),
\]
and $B(0) \in C_*$, then we must have $B(t) \in C_*$ for all $t$, as we wanted to show.

We now show $S$ is parallel with respect to $\tilde{\nabla}$.
Let $X,Y \in T_xM$ and parallel transport them to $\tilde{T}(X),\tilde{T}(Y) \in T_yM$ to some nearby $y$.
Choosing $b \in (P_*)_x$ we may write $X = bu$ and $Y = bv$ for some $u,v \in V$.
By what we just showed, there exists $b' \in (P_*)_y$ such that $\tilde{T}b = b'$.
In particular,
\[
\tilde{T}(X) = \tilde{T}bu = b'u \ \text{and} \ \tilde{T}(Y) = \tilde{T}bv = b'v.
\]
Therefore,
\[
S(\tilde{T}X,\tilde{T}Y) = S^{(b')}(b'u,b'v) = b'A_0(u,v) = \tilde{T}bA_0(u,v) = \tilde{T}S(X,Y).
\]

Finally, since $\nabla$ is torsion free, the torsion $T$ of $\tilde{\nabla}$ is given by
\[
T(X,Y) = -S(X,Y) + S(Y,X).
\]
Since $S$ is parallel, so is $T$.
\end{proof}

\subsection{Vector bundle valued currents}

We take a small detour to recall the theory of currents taking values on a vector bundle.
In this text we only apply it to local charts, so we assume an orientable manifold for simplicity.

Let $E \to M$ be a smooth vector bundle over a smooth $n$-dimensional oriented manifold $M$.
For $0 \leq k \leq n$, a \emph{$E$-valued $k$-current} is an element of $(\Omega^{n-k}_c(M,E^*))'$.
If $\omega \in L^1_{loc}(M,\bigwedge^kT^*M \otimes E)$ is a locally $L^1$ $E$-valued $k$-form, then it corresponds to the $E$-valued $k$-current defined by
\[
\eta \mapsto \int \text{ev}(\omega \wedge \eta), \ \text{for} \ \eta \in \Omega_c^{n-k}(M,E^*).
\]
where $\text{ev}$ is the evaluation map $\text{ev} \colon E \otimes E^*\otimes \bigwedge^n T^*M  \to \bigwedge^n T^*M$.

In the case $E = \text{End}(TM)$,
$E$ is isomorphic to $E^*$ by identifying $A \in \text{End}(TM)$ with the map $\text{End}(TM) \to \text{End}(TM)$ defined by $B \mapsto \text{tr}(AB)$.
With this identification, the evaluation pairing between elements in $A,B$ in $E$ is given by $\text{tr}(AB)$.
 
\subsection{Curvature and the Bianchi identities}

We defined a H\"older connection $\tilde{\nabla}$.
In this section we define its curvature in the distributional sense and show that it satisfies the Bianchi identities under some regularity assumptions.

Fix a smooth local chart in $U \subset M$ inducing a frame $e_x \colon V \to T_xM$ and denoting the coordinate vector fields of this frame by $\{\partial_1,...,\partial_n\}$.
Let $A \in C^{\beta}(U,T^*U \otimes \text{End}(V))$ be the $\text{End}(V)$-valued $1$-form characterized by
\[
\tilde{\nabla}_X(ev) = e(X(v)+A(X)v),
\]
for every $v \colon U \to V$ and vector field $X$ on $U$.

In these coordinate charts, the curvature of $\tilde{\nabla}$ is the $\text{End}(V)$-valued distributional $2$-current
\[
F = dA + A\wedge A.
\]
Here $dA$ is the exterior derivative of the distributional $\text{End}(V)$-valued $1$-current, making it a $\text{End}(V)$-valued $2$-current.
For a compactly supported smooth $\text{End}(V)$-valued $n-2$ form $\phi \in \Omega^{n-2}(U,\text{End}(V))$,
\[
\langle dA,\phi \rangle = \int_U \text{tr}(A \wedge d\phi).
\]
\begin{remark}
\label{rem:C1test}
	Note that $dA$ also acts on $C^1$ test forms. Consequently, so does $F$.
\end{remark}
For $X,Y,Z$ smooth vector fields on $U$, we define the distributional curvature form by
\[
R(X,Y)Z = e(F(X,Y)e^{-1}Z).
\]
\subsubsection{First Bianchi identity}

We first verify the first Bianchi identity when the torsion is $C^1$.
\begin{lemma}
\label{lem:1stbianchi}
In the above notation and assuming that the torsion of $T$ is $C^1$, the curvature satisfies the first Bianchi identity, that is, for $X,Y,Z$ smooth vectors fields on $U$,
	\begin{equation}
		\label{eq:1stbianchi}
		\sum_{\text{cyc}} R(X,Y)Z = \sum_{\text{cyc}} (\tilde{\nabla}_XT)(Y,Z) + T(T(X,Y),Z).
	\end{equation}
\end{lemma}
\begin{proof}
Denote by $\theta = e^{-1}$, which satisfies $d\theta = 0$.
Define the torsion $V$-valued $2$-form by
\[
\Theta(X,Y) = \theta T(X,Y).
\]
	By the definition of torsion,
	\begin{align*}
	\Theta(X,Y) =& X(\theta(Y)) - Y(\theta(X)) - \theta([X,Y]) + A(X)\theta(Y) - A(Y)\theta(X)\\
	=& d\theta(X,Y) + (A\wedge \theta)(X,Y) = (A \wedge \theta)(X,Y).
	\end{align*}
	Now wedging the definition of $F$ with $\theta$ and substituting $A \wedge \theta = \Theta$ gives us
	\begin{align*}
	F \wedge \theta =& dA \wedge \theta + A \wedge A \wedge \theta\\
	=& d(A\wedge \theta) + A \wedge (A \wedge \theta) \\
	=& d\Theta + A \wedge \Theta.
	\end{align*}
	
	For $X,Y,Z$ smooth vector fields on $U$,
	\[
	(F \wedge \theta)(X,Y,Z) = \sum_{\text{cyc}} F(X,Y)\theta(Z),
	\]
	which is a well defined distribution since $\theta$ is smooth.
	On the other hand,
	\[
	d\Theta(X,Y,Z) = \sum_{\text{cyc}} (X(\Theta(Y,Z)) - \Theta([X,Y],Z))
	\]
	and
	\[
	(A\wedge \Theta)(X,Y,Z) = \sum_{\text{cyc}} A(X)\Theta(Y,Z).
	\]
	The connection formula gives us
	\begin{equation}
	\label{eq:connectionformula1}
	e((d\Theta+A\wedge\Theta)(X,Y,Z)) = \sum_{\text{cyc}} \tilde{\nabla}_X(T(Y,Z)) - T([X,Y],Z).
	\end{equation}
	The derivative of $T(Y,Z)$ in the direction of $X$ gives us
	\[
	\tilde{\nabla}_X(T(Y,Z)) = (\tilde{\nabla}_XT)(Y,Z) + T(\tilde{\nabla}_XY,Z) + T(Y,\tilde{\nabla}_XZ).
	\]
	Substituting in \eqref{eq:connectionformula1} gives us
	\begin{equation}
	\label{eq:1stbianchiproof}
	e((d\Theta+A\wedge\Theta)(X,Y,Z)) =
	 \sum_{\text{cyc}} (\tilde{\nabla}_XT)(Y,Z) +
	  \sum_{\text{cyc}} T(\tilde{\nabla}_XY,Z) + T(Y,\tilde{\nabla}_XZ) - T([X,Y],Z).
	\end{equation}
	Focusing on the second summand on the RHS of the above, we cycle the middle summand and use anti-symmetry of $T$ to obtain that it is equal to
	\[
	\sum_{\text{cyc}} T(\tilde{\nabla}_XY,Z) - T(\tilde{\nabla}_YX,Z) - T([X,Y],Z)
	= \sum_{\text{cyc}}T(T(X,Y),Z).
	\]
	Substituting in \eqref{eq:1stbianchiproof}, we obtain the first Bianchi identity.
\end{proof}

\subsubsection{Second Bianchi identity}
The second Bianchi identity will now be verified under the assumption of $C^1$ curvature tensor.

\begin{lemma}
\label{lem:2ndbianchi}
		If the curvature $R$ has a $C^1$ representative, then the second Bianchi equation holds, meaning
	\begin{equation}
		\label{eq:2ndbianchi}
		\sum_{\text{cyc}} (\tilde{\nabla}_XR)(Y,Z) + R(T(X,Y),Z)  =0.
	\end{equation}
\end{lemma}
\begin{proof}
	$R$ having a $C^1$ representative implies that $F$ is $C^1$.
	In particular, $dA =  F - A \wedge A$ is continuous.
	Convolve $A$ with a mollifier $\rho_{\varepsilon}$ defining the smooth forms
	\[
	A_{\varepsilon} = A * \rho_{\varepsilon} \ \text{and} \ F_{\varepsilon} = dA_{\varepsilon} + A_{\varepsilon} \wedge A_{\varepsilon}.
	\]
	Taking the exterior derivative of $F_{\varepsilon}$ gives us
	\[
	dF_\varepsilon = dA_\varepsilon \wedge A_\varepsilon - A_{\varepsilon} \wedge dA_\varepsilon.
	\]
	Substituting $dA_{\varepsilon} = F_{\varepsilon} - A_{\varepsilon} \wedge A_{\varepsilon}$,
	\begin{align*}
	dF_{\varepsilon} =& (F_\varepsilon \wedge A_\varepsilon - A_\varepsilon \wedge A_\varepsilon \wedge A_\varepsilon) - (A_\varepsilon \wedge F_\varepsilon - A_\varepsilon \wedge A_\varepsilon \wedge A_\varepsilon) \\
	=& F_\varepsilon \wedge A_\varepsilon - A_\varepsilon \wedge F_\varepsilon.
	\end{align*}
	Since the exterior derivative commutes with convolution, meaning $dA_{\varepsilon} = (dA) * \rho_{\varepsilon}$, then $A_{\varepsilon} \rightarrow A$ and $dA_{\varepsilon} \rightarrow dA$ uniformly on compact sets as $\varepsilon \rightarrow 0$.
	In particular, $F_\varepsilon \rightarrow F$ uniformly on compact sets and $dF_{\varepsilon} \rightarrow dF$ distributionally.
	Letting $\varepsilon \rightarrow 0$ in the above, we obtain the following equation distributionally
	\begin{equation}
		\label{eq:2ndbianchitensor}
		dF + A \wedge F - F \wedge A = 0,
	\end{equation}
	which also holds pointwise since it involves only continuous tensors.
	
	Let $X,Y,Z$ be smooth vector fields.
	Evaluating \eqref{eq:2ndbianchitensor} on them gives us
	\begin{equation}
	\label{eq:2ndbianchialmost}
	\sum_{\text{cyc}} (X(F(Y,Z))-F([X,Y],Z)) + [A(X),F(Y,Z)] = 0,
		\end{equation}
	where the bracket is the bracket of endomorphisms.
	For $Q$ a $C^1$ section of $\text{End}(TM)$ and $q = e^{-1} Q e$,
	\[
	e^{-1}(\tilde{\nabla}_XQ)e = X(q) + [A(X),q].
	\]
	Therefore, conjugating \eqref{eq:2ndbianchialmost} by $e$,
	\[
	\sum_{\text{cyc}} \tilde{\nabla}_X(R(Y,Z)) - R([X,Y],Z) = 0
	\]
	Substituting $\tilde{\nabla}_X(R(Y,Z)) = (\tilde{\nabla}_XR)(Y,Z) + R(\tilde{\nabla}_XY,Z) + R(Y,\tilde{\nabla}_XZ)$ in the above and cycling the $R(Y,\tilde{\nabla}_XZ)$ term inside the sum,
	\begin{align*}
	\sum_{\text{cyc}} (\tilde{\nabla}_XR)(Y,Z) + R(\tilde{\nabla}_XY,Z) + R(Z,\tilde{\nabla}_YX) - R([X,Y],Z) \\
	= \sum_{\text{cyc}} (\tilde{\nabla}_XR)(Y,Z) + R(\tilde{\nabla}_XY - \tilde{\nabla}_YX - [X,Y],Z)\\
	= \sum_{\text{cyc}}  (\tilde{\nabla}_XR)(Y,Z)  + R(T(X,Y),Z) = 0.
	\end{align*}
	This shows the Lemma.
	\end{proof}
	
	\subsection{Injectivity of cyclic sum}
	
	\begin{lemma}
		\label{lem:algebralemma}
		The linear map $\delta \colon \Hom(\bigwedge^2 V, \mathfrak{c}_*)  \to \Hom(\bigwedge^3V,V) $
		\begin{equation*}
			(\delta Q)(u,v,w) = Q(u,v)w + Q(v,w)u + Q(w,u)v
		\end{equation*}
		is injective.
	\end{lemma}
\begin{proof}
 Recall that $C_*$ commutes with $\rho(G)$. Thus, every $S\in\clie_*$ preserves each restricted weight space. Also, $V^G=\{0\}$, since $V$ has no zero restricted weight.

Suppose that $\delta Q=0$. Fix $\lambda, \nu\in \Sigma_V$ and fix $u\in V_\lambda$ and $v\in V_\nu$.
For every $w\in V$, the cyclic identity gives \[ Q(u,v)w=-Q(v,w)u-Q(w,u)v.\] Since $Q(v,w)$ and $Q(w,u)$ belong to $\clie_*$, $Q(v,w)u\in V_\lambda$ and $Q(w,u)v\in V_\nu$.  Consequently, \[ W\coloneqq\textrm{Im}Q(u,v)\subset V_\lambda \oplus V_\nu.\]
Moreover, $W$ is $G$-invariant: $\rho(g)Q(u,v)w=Q(u,v)\rho(g)w\in W$ for $g\in G$ and $w\in V$.

\begin{claim}\label{W=0} $W=0$ \end{claim}
\begin{proof}[Proof of \cref{W=0}] Suppose $W\neq 0$ for the sake of the contradiction. Since $V^G=\{0\}$, the representation $\rho\restrict{W}$ of $G$ on $W$ is nontrivial. As $G$ is connected and its real Lie algebra $\glie$ is simple, the differential $d(\rho\restrict{W})$ is injective. Hence, the set $\Sigma_W$ of restricted weights of $W$ span $\alie^*$. Indeed, if $H\in\alie$ is vanished by all these weights, then $d(\rho\restrict{W})(H)$ vanishes on every weight space of $W$, hence, $d(\rho\restrict{W})(H)=0$ so that $H=0$. 

On the other hand, $\glie=[\glie,\glie]$ and the trace of a commutator is zero. Therefore, we have
\[ 0=\tr (d(\rho\restrict{W})(H))=\sum_{\chi\in\Sigma_W}(\dim_{\Rbb}W_\chi)\,\chi(H)\] for all $H\in \alie$.

Therefore, $\Sigma_W$ span $\alie^*$ and satisfy a nontrivial linear relation. Hence, there must be at least $\dim_\Rbb\alie+1\ge3$ distinct weights. This contradicts $W\subset V_\lambda\oplus V_\mu$, which gives $\Sigma_W\subset\{\lambda,\mu\}$.
\end{proof}
Hence, $W=0$, i.e., $Q(u,v)=0$ for every pair of weight vectors $u$ and $v$. Since $V=\bigoplus_{\nu\in\Sigma_V} V_\nu$, bilinearity gives $Q=0$.
\end{proof}
	
	\subsection{Improving curvature regularity}
	
	Since $\tilde{\nabla}T = 0$ in our case, \eqref{eq:1stbianchi} implies
	\begin{equation}
	\label{eq:1stbianchi0derivative}
	\sum_{\text{cyc}} R(X,Y)Z  = \sum_{\text{cyc}} T(T(X,Y),Z).
	\end{equation}
	Since $R$ is only a distribution, the above is a distributional equation.
	However, the right hand since is a continuous tensor, meaning that so it the left hand side.
	We will now  show that $R$ itself is $C^{1+\beta}$.
	
	\begin{lemma}
    \label{lem:regular-curvature}
		The curvature $R$ is $C^{1+\beta}$ and is parallel.
	\end{lemma}
	\begin{proof}
	Let $b$ be a $C^{1+\beta}$ local section of $P_*$ on a small oriented neighborhood $U$.
	With respect to this local section, we may write the connection as
	\[
	\tilde{\nabla}_X(bv) = bB(X)v,
	\]
	for $X$ vector field and constant $v \in V$, where $B$ is a H\"older $\text{End}(V)$-valued $1$-form.
	We previously showed that $B(X) \in \mathfrak{c}_*$ in \eqref{eq:Bincstar}.
	In this frame, the curvature is given by
	\[
	R_b = dB + B \wedge B.
	\]
	Fixing a basis $\{E_i\}_i$ of $\mathfrak{c}_*$, we may write
	\[
	B = \sum_i \omega_i E_i,
	\]
	where $\omega_i \in C^{\beta}(U,T^*U)$.
	Taking the exterior derivative, the following holds distributionally
	\[
	dB = \sum_i d\omega_i E_i,
	\]
	which still takes values in $\mathfrak{c}_*$.
	Meanwhile, for $X,Y$ vector fields in $U$
	\[
	B\wedge B(X,Y) = B(X)B(Y) - B(Y)B(X) = [B(X),B(Y)] \in \mathfrak{c}_*,
	\]
	since $\mathfrak{c}_*$ is a Lie algebra.
	Therefore, $R_b$ is a $\mathfrak{c}_*$-valued $2$-current.
	 
	Consequently, expressing the inputs of $R_b$ using the frame $b$, define the $\Hom(\bigwedge^2 V,\mathfrak{c}_*)$-valued distribution by
	\begin{equation}
	\label{eq:hatRdef}
	\hat{R}(u,v)w = b^{-1}R(bu,bv)bw,
	\end{equation}
	for $u,v,w \in V$.
	Although $b$ is only $C^1$, the above is a well defined distribution by \cref{rem:C1test}.
	
	We now define $T_0 \colon \bigwedge^2V \to V$ by
	\[
	T_0(u,v) = A_0(v)u - A_0(u)v.
	\]
	Since $A$ is $C_*$-equivariant, then so is $T_0$, meaning that for $c \in C_*$,
	\begin{equation}
		\label{eq:T0equivariant}
		T_0(cu,cv) = cT_0(u,v).
	\end{equation}
	By the definition of torsion of $\tilde{\nabla} = \nabla - S$ and the definition of $S$,
	\[
	T(bu,bv) = bT_0(u,v).
	\]
	Consequently, the map $\hat{T} \colon \bigwedge^3V \to V$ defined by
	\[
	\hat{T}(u,v,w) = \sum_{\text{cyc}}T_0(T_0(u,v),w)
	\]	
	is constant and satisfies
	\[
	\sum_{\text{cyc}} b^{-1}T(T(bu,bv),bw) = \hat{T}(u,v,w).
	\]
	Writing \eqref{eq:1stbianchi0derivative} in the frame $b$ gives us
	\begin{equation}
	\label{eq:distribution1stbianchi}
		\sum_{\text{cyc}}\hat{R}(u,v)w = \hat{T}(u,v,w),
	\end{equation}
	as distributions.
	Let $\mu$ be a compactly supported $n$-form on $U$ with $\int_U \mu = 1$.
	Define $R_0 \in \Hom(\bigwedge^2V,\mathfrak{c}_*)$ by $R_0 = \hat{R}(\mu)$.
	Applying $\mu$ to \eqref{eq:distribution1stbianchi},
	\begin{equation*}
		\delta R_0 = \hat{T} \int_U \mu = \hat{T},
	\end{equation*}
	where $\delta$ is the map from \cref{lem:algebralemma}.
	For $\nu$ a compactly supported $n$-form on $U$,
	\begin{equation*}
		\delta(\hat{R}(\nu)) = \hat{T}\int_U \nu = \delta R_0 \int_U \nu.
	\end{equation*}
	By injectivity of $\delta$,
	\[
	\hat{R}(\nu) = R_0 \int_U \nu, \ \text{for all} \ \nu \in \Omega_c^n(U).
	\]
	Therefore, $\hat{R}$ is represented by the constant map $R_0$.
	Consequently, $R$ has representative
	\[
	R_x(X,Y)Z = b_xR_0(b_x^{-1}X,b_x^{-1}Y)b_x^{-1}Z.
	\]
	Since $b$ is $C^{1+\beta}$, so is $R$.
	
	In order to show $R$ is parallel, we start by showing that $R_0$ is $C_*$-equivariant.
	Indeed, for $c \in C_*$ define
	\begin{equation*}
		R'_0(u,v) = c^{-1}R_0(cu,cv)c.
	\end{equation*}
	It takes values in $\mathfrak{c}_*$ since conjugation by elements of $C_*$ preserves $\mathfrak{c}_*$.
	Since $T_0$ is $C_*$-equivariant, so is $\hat{T}$.
	Therefore,
	\begin{equation*}
		(\delta R'_0)(u,v,w) = c^{-1}(\delta R_0)(cu,cv,cw)
		= c^{-1}\hat{T}(cu,cv,cw) = \hat{T}(u,v,w) = (\delta R_0)(u,v,w).
	\end{equation*}
	By injectivity of $\delta$, $R'_0 = R_0$, meaning
	\begin{equation}
		\label{eq:R0Cequivariant}
		R_0(cu,cv)cw = cR_0(u,v)w, 
	\end{equation}	
	showing equivariance of $R_0$.
	
	Now to show $R$ is parallel, let $X,Y,Z \in T_xM$ for some $x \in U$ and parallel transport them to $\mathcal{P}(X),\mathcal{P}(Y), \mathcal{P}(Z) \in T_yM$ to some nearby $y \in U$.
	We may write $X = bu$, $Y = bv$ and $Z = bw$ for some $u,v,w \in V$.
	Recall from \cref{prop:newconnred} that there exists $c \in C_*$ such that $\mathcal{P}b_x = b_yc$.
	In particular,
	\[
	\mathcal{P}(X) = b_ycu,  \mathcal{P}(Y) = b_ycv \ \text{and} \ \mathcal{P}(Z) = b_ycw.
	\]
	Therefore,
	\begin{equation}
	\label{eq:Rparallel1}
	R(\mathcal{P}(X),\mathcal{P}(Y))\mathcal{P}(Z)
	= R(b_ycu,b_ycv)b_ycw
	= b_yR_0(cu,cv)cw.
	\end{equation}
	On the other hand,
	\begin{equation}
		\label{eq:Rparallel2}
	\mathcal{P}(R(X,Y)Z) = \mathcal{P}(b_xR_0(u,v)w) 
	= b_ycR_0(u,v)w.
	\end{equation}
	By \eqref{eq:R0Cequivariant}, we conclude that \eqref{eq:Rparallel1} and \eqref{eq:Rparallel2} are equal, meaning $R$ is parallel and finishing the proof.
	\end{proof}
	
	\subsection{An auxiliary Lie algebra}

	\begin{lemma}
	\label{lem:auxiliaryLie}
		On the vector space $\mathfrak{h} = \mathfrak{c}_* \oplus V$, define a bilinear antisymmetric map $[\cdot,\cdot] \colon \mathfrak{h} \times \mathfrak{h} \to \mathfrak{h}$ by
		\begin{equation*}
			[A,B] = AB-BA, \ [A,u] = Au, \ [u,v] = -R_0(u,v) -T_0(u,v).
		\end{equation*}
		Then $\mathfrak{h}$ equipped with $[\cdot,\cdot]$ is a Lie algebra.
	\end{lemma}
	\begin{proof}
	We will verify that the Jacobi expression
	\[
	J(x,y,z) = \sum_{\text{cyc}} [[x,y],z]
	\]
	vanishes. Since $J$ is trilinear and alternating, it suffices to verify its vanishing when each $x,y,z$ is an element of $\mathfrak{c}_*$ or $V$.
	
	If $x,y,z \in \mathfrak{c}_*$, then the vanishing of $J(x,y,z)$ follows from the Jacobi identity for the commutator of endomorphisms.
	For $x,y = A,B \in \mathfrak{c}_*$ and $z = w \in V$,
	\begin{align*}
	J(x,y,z) =& [A,B]w + [Bw,A] - [Aw,B] \\
	=& (AB-BA)w - ABw + BAw = 0.
	\end{align*}
	
	For $x = A \in \mathfrak{c}_*$ and $y,z = v,w \in V$, we will use equivariance of $R_0$ and $T_0$.
	Let $c(t) = \exp(tA)$ which satisfies $c'(0) = A$.
	Since $c(t) \in C_*$, \eqref{eq:T0equivariant} implies
	\begin{equation*}
		T_0(c(t)v,c(t)w) = c(t)T_0(v,w).
	\end{equation*}
	Differentiating at $t=0$ we obtain
	\begin{equation}
		\label{eq:vectorterm}
		T_0(Av,w) + T_0(v,Aw) = AT_0(v,w).
	\end{equation}
	Similarly using \eqref{eq:R0Cequivariant},
	\begin{equation*}
		R_0(c(t)v,c(t)w)c(t) = c(t)R_0(v,w).
	\end{equation*}
	Differentiating at $t=0$ gives us
	\begin{equation}
		\label{eq:matrixterm}
		R_0(Av,w) + R_0(v,Aw) + R_0(v,w)A = AR_0(v,w).
	\end{equation}
	Expanding $J(A,v,w)$ gives 
	\begin{equation*}
		J(A,v,w) = ([A,R_0(v,w)] - R_0(Av,w)-R_0(v,Aw)) + (AT_0(v,w)-T_0(Av,w)-T_0(v,Aw)).
	\end{equation*}
	The first term above is $0$ by \eqref{eq:matrixterm} while the second term vanishes due to \eqref{eq:vectorterm}.
	
	If $x,y,z = u,v,w \in V$, then
	\[
	[[u,v],w] = (R_0(T_0(u,v),w)) + (T_0(T_0(u,v),w)-R_0(u,v)w).
	\]
	The cyclic sum of the first sum above is $0$ by \eqref{eq:2ndbianchi} together with the fact that the curvature is parallel and the second is $0$ due to \eqref{eq:1stbianchi} and torsion being parallel.
	Therefore, $J$ vanishes and $[\cdot,\cdot]$ defines a Lie bracket.
	\end{proof}

	\subsection{Nilpotent Lie subalgebra lemma}
	
	\begin{lemma}
	\label{lem:nilpotentV}
		Let $\mathfrak{h} = \mathfrak{c}_* \oplus V$ be a finite dimensional Lie algebra written as a direct sum of vector spaces.
		Suppose that
		\begin{enumerate}
			\item $\mathfrak{c}_*$ is a Lie subalgebra, $[\mathfrak{c}_*,V] \subset V$ and its action on $V$ is faithful;
			\item there exists an inner product on $\mathfrak{h}$ making $\text{ad} (A)$ skew-adjoint for every $A \in \mathfrak{c}_*$;
			\item there exists a derivation $L$ of $\mathfrak{h}$ diagonalizable over $\mathbb{R}$ satisfying $\mathfrak{c}_* \subset \ker L$ and $L|_V \colon V \to V$ invertible.
		\end{enumerate}
		Then $V$ is a nilpotent Lie subalgebra.
	\end{lemma}
	
	\begin{proof}
		Consider the Killing form $\kappa(x,y) = \text{tr}(ad(x)ad(y))$. Its kernel 
		\[
		\ker \kappa = \{ x \in \mathfrak{h} \ | \ \kappa(x,y) = 0, \forall y \in \mathfrak{h} \}
		\]
		is an ideal.
		We will first show $\ker \kappa = V$.
		
		To show $V \subset \ker \kappa$, we use $V = L(\mathfrak{h})$.
		This means we need to verify $\kappa(Lx,y) = 0$, for all choices of $x, y \in \mathfrak{h}$.
		Since $L$ is a derivation, for $z \in \mathfrak{h}$,
		\begin{equation}
		\label{eq:adL}
		ad(Lx)z = [Lx,z] = L[x,z]-[x,Lz] = (L ad(x) - ad(x) L)(z).
		\end{equation}
		Substituting in the definition of the Killing form,
		\begin{align*}
		\kappa(Lx,y) = & \text{tr}(ad(Lx)ad(y)) \\
		= & \text{tr}(Lad(x)ad(y)) - \text{tr}(ad(x)Lad(y)).
		\end{align*}
		Changing the last term to $\text{tr}(Lad(y)ad(x))$ using cyclicity of trace, regrouping the terms and using Jacobi to write $[ad(x),ad(y)] = ad([x,y])$,
		\begin{equation}
			\label{eq:kappaL}
			\kappa(Lx,y) = \text{tr}(Lad([x,y])).
		\end{equation}
		Write $f(z) = \text{tr}(Lad(z))$.
		We will verify that $f$ vanishes on $\mathfrak{h}$ by verifying it vanishes separately on $V$ and $\mathfrak{c}_*$.
		
		For $v \in V$, write $Lx = v$ for some $x \in \mathfrak{h}$ and use \eqref{eq:adL} with the cyclicity of trace to obtain
		\begin{align*}
			f(v) = \text{tr}(Lad(Lx)) =& \text{tr}(L(L ad(x) - ad(x) L)) \\
			=& \text{tr}(L^2ad(x))-\text{tr}(Lad(x)L) \\
			=& \text{tr}(L^2ad(x)) - \text{tr}(L^2ad(x)) = 0.
		\end{align*}
		Therefore, $f$ vanishes on $V$.
		
		For $A \in \mathfrak{c}_*$, we have $LA = 0$.
		Combining it with \eqref{eq:adL} gives
		\[
		L,ad(A) = ad(A)L.	
		\]
		Since $L$ is diagonalizable, we may write $\mathfrak{h} = \bigoplus_{\lambda} E_\lambda$ as a sum of eigenspaces of $L$ and $ad(A)$ preserves each of these eigenspaces.
		Moreover,
		\[
		Lad(A)|_{E_{\lambda}} = ad(A)L|_{E_{\lambda}} = \lambda ad(A)|_{E_{\lambda}}.
		\]
		Since $ad(A)$ is skew-symmetric, it has trace $0$ and so does $Lad(A)|_{E_{\lambda}}$.
		The trace being $0$ on each eigenspace implies that the trace is $0$, that is,
		\[
		f(A) = \text{tr}(Lad(A)) = 0.
		\]
		This shows $f$ vanishes on $\mathfrak{c}_*$.
		Consequently, $f$ vanishes on all of $\mathfrak{h}$.
		By \eqref{eq:kappaL}, $\kappa(Lx,y) = 0$, for all $x, y \in \mathfrak{h}$, which shows $\text{Im}(L) = V \subset \ker \kappa$.
		
		For the other direction, suppose by that $A+v \in \ker \kappa$ for $A \in \mathfrak{c}_* \setminus \{0\}$ and $v \in V$.
		Then, using $v \in \ker \kappa$ and skew-adjointness of $ad(A)$,
		\[
		0 = \kappa(A+v,A) = \kappa(A,A) = \text{tr}(ad(A)^2) = -\text{tr}(ad(A)^*ad(A)).
		\]
		Therefore, $ad(A) = 0$, which implies that $A =0$ by faithfulness and that $V = \ker \kappa$ is an ideal.
		
		In particular, $[V,V] \subset V$, making $V$ a finite dimensional Lie subalgebra admitting an invertible derivation.
		By Jacobson's theorem (\cite{Jacobson}), $V$ is a nilpotent Lie algebra.		
	\end{proof}

	\begin{lemma}
	\label{lem:verify}
		The Lie algebra defined in \cref{lem:auxiliaryLie} satisfies the hypothesis of \cref{lem:nilpotentV}.
	\end{lemma}
	
	\begin{proof}
		We showed in \cref{lem:auxiliaryLie} that $\mathfrak{h} = \mathfrak{c}_* \oplus V$ is a Lie algebra such that $[\mathfrak{c}_*,V] \subset V$.
		Moreover, for $A \in \mathfrak{c}_*$, if $[A,v] = Av = 0$, for all $v \in V$, then $A = 0$, showing that the action of $\mathfrak{c}_*$ on $V$ is faithful.
		
		We now define an inner product making $ad(A)$ skew-adjoint.
		First, consider the $C_*$ action on $\mathfrak{h}$ defined by $c(B+v) = cBc^{-1} +cv$.
		For $A \in \mathfrak{c}_*$, let $c(t)$ be a smooth curve on $C_*$ satisfying $c(0) = \text{id}$ and $c'(0) = A$.
		For an element $x = B+v \in \mathfrak{h}$, define
		\[
		x(t) = c(t)x = c(t)Bc(t)^{-1} + c(t)v.
		\]
		Differentiating the above at $0$ we obtain
		\begin{equation}
		\label{eq:xprimegivesad}
		x'(0) = AB-BA +Av = [A,B]+[A,v] = ad(A)x.
		\end{equation}
		Now take any inner product on $\mathfrak{h}$ and average it over the action of the compact group $C_*$ to obtain a $C_*$-invariant inner product $\langle \cdot, \cdot \rangle$.
		For $x,y \in \mathfrak{h}$, it satisfies
		\[
		\langle c(t)x,c(t)y \rangle = \langle x, y \rangle.
		\]
		Differentiating the above and using \eqref{eq:xprimegivesad},
		\[
		\langle ad(A)x,y \rangle + \langle x, ad(A) y \rangle = 0,
		\]
		which is the desired skew-adjointness property.
		
		Now in order to construct the derivation, we first show that the bracket is preserved under the $G$-action $g(A+v) = A+\rho(g)v$.
		For elements $A+u,B+v \in \mathfrak{h}$ and $g \in G$,
		\begin{align*}
		[g(A+u),g(B+v)] = & [A+\rho(g)u,B+\rho(g)v] \\
		=&  [A,B]+\rho(g)(Av-Bu) -R_0(\rho(g)u,\rho(g)v)-T_0(\rho(g)u,\rho(g)v).
		\end{align*}
		On the other hand,
		\begin{equation*}
			g[A+u,B+v] =  [A,B] + \rho(g)(Av-Bu) - R_0(u,v)-\rho(g)T_0(u,v).
		\end{equation*}
		Therefore, it suffices to verify
		\[
		R_0(\rho(g)u,\rho(g)v) = R_0(u,v) \ \text{and} \ T_0(\rho(g)u,\rho(g)v) = \rho(g)T_0(u,v).
		\]
		For $T_0$, we use $G$-equivariance of $A_0$ from \cref{lem:unique-A0-lift}
		\begin{align*}
		T_0(\rho(g)u,\rho(g)v) = & A_0(\rho(g)v)\rho(g)u - A_0(\rho(g)u)\rho(g)v\\
		= & \rho(g)(A_0(v)u-A_0(u)v) = \rho(g)T_0(u,v).
		\end{align*}
		Note that this not only shows $G$-equivariance for $T_0$, but also for $\hat{T}$.
		
		For $R_0$, define $R_g \in \Hom(\bigwedge^2 V, \mathfrak{c}_*)$ by
		\[
		R_g(u,v) = \rho(g)^{-1}R_0(\rho(g)u,\rho(g)v)\rho(g).
		\]
		It indeed maps to $\mathfrak{c}_*$ since $\rho(G)$ commutes with $C_*$ and, consequently, with $\mathfrak{c}_*$.
		Moreover, for the linear map $\delta$ defined in \cref{lem:algebralemma},
		\begin{align*}
			(\delta R_g)(u,v,w) = & \rho(g)^{-1}(\delta R_0)(\rho(g)u,\rho(g)v,\rho(g)w)\\
			=& \rho(g)^{-1}\hat{T}(\rho(g)u,\rho(g)v,\rho(g)w) \\
			=& \hat{T}(u,v,w) = (\delta R_0)(u,v,w).
		\end{align*}
		By injectivity of $\delta$ from \cref{lem:algebralemma}, we conclude that $R_g = R_0$, which proves $G$-equivariance for $R_0$ and, consequently, for the Lie bracket.
		
		Let $\mathfrak{a} = \text{Lie}(A)$ where $A \subset G$ is a maximal $\mathbb{R}$-split torus.
        Recall from \cref{sec:adaptmetric} that we write
        \[
        V = \bigoplus_{\lambda \in \Sigma_V} V_\lambda,
        \]
        where $\Sigma_V$ is the set of restricted weights of the representation $\rho$.
        By \cref{prop:nozero}, there exists $H \in \mathfrak{a}$ that satisfies $\lambda(H) \neq 0$ for all $\lambda \in \Sigma_V$.
		
		Let $g(t) = \exp(tH) \in G$.
		We define $L \colon \mathfrak{h} \to \mathfrak{h}$ by
		\[
		L(A+u) = \frac{d}{dt}\Big|_{t=0}g(t)(A+u) =  \frac{d}{dt}\Big|_{t=0}(A+\rho(g(t))u) = d\rho(H)u.
		\]
		For $x,y \in \mathfrak{h}$, $G$-equivariance of the Lie bracket implies
		\[
		g(t)[x,y] = [g(t)x,g(t)y].
		\]
		Differentiating the above at $t=0$ gives
		\[
		L[x,y] = [Lx,y] + [x,Ly],
		\]
		making $L$ a derivation on $\mathfrak{h}$.
		Moreover, from its definition we have $L(A) = 0$, for $A \in \mathfrak{c}_*$ and, for $v \in V_{\lambda}$,
		\[
		L(v) = d\rho(H)v = \lambda(H)v \neq 0,
		\]
		meaning $L|_V$ is invertible and showing the Lemma.
	\end{proof}
	
	\begin{corollary}
		\label{cor:curvature0}
		$R_0 = 0$ and $-T_0$ defines a nilpotent Lie algebra on $V$.
	\end{corollary}
	\begin{proof}
		By \cref{lem:verify}, the Lie algebra $\mathfrak{h} = \mathfrak{c}_* \oplus V$ satisfies \cref{lem:nilpotentV}, making $V$ a nilpotent Lie subalgebra.
		In particular, for $u,v \in V$
		\[
		[u,v] = -R_0(u,v) - T_0(u,v) \in V,
		\]
		making the $\mathfrak{c}_*$ of the above $0$, that is, $R_0(u,v) =0$.
		It also follows from \cref{lem:nilpotentV} that $-T_0$ defines a nilpotent Lie bracket on $V$.
	\end{proof}

\subsection{Adapted Metric}

\begin{lemma}
\label{sec:adapted-metric}
    There exists a $C^{1+\beta}$ metric which is parallel with respect to $\tilde{\nabla}$.
\end{lemma}
\begin{proof}
Recall that in \cref{sec:adaptmetric} we defined an inner product $\langle\cdot,\cdot\rangle_V$ on $V$ which is $C$-invariant.
By \cref{prop:reduction1}, $C_* \subset C$.
Therefore, $\langle\cdot,\cdot\rangle_V$ is $C_*$ invariant.

For $p \in P_*$ a local frame and $X,Y$ local vector fields, define
\[
g(X,Y) = \langle p^{-1}X,p^{-1}Y\rangle_V.
\]
Since a different choice of local frame in $P_*$ differs by elements in $C_*$, then $g$ does not depend on the choice of $p$ and is well defined.
Moreover, since $p$ is $C^{1+\beta}$, then so is $g$.
Let $\mathcal{P}$ be the parallel transport with respect to $\tilde{\nabla}$ between two points $x,y$ along a $C^1$ curve.
By \cref{prop:newconnred}, $\mathcal{P}p_x =  p'\in (P_*)_y$.
Therefore, for $p_xu,p_xv \in T_xM$,
\[
g(\mathcal{P}p_xu,\mathcal{P}p_xv) = g(p'u,p'v)
= \langle u,v\rangle_V = g(p_xu,p_xv).
\]
This shows that $g$ is parallel.
\end{proof}

\subsection{The resulting flat connection $\wtd\nabla$}

\begin{theorem}[A flat metric connection with parallel torsion]
\label{sec6conclusion}
Under the standing hypotheses of this section, the reduction $P_*$
from Proposition~\ref{prop:reduction1} and the connection
$\wtd\nabla=\nabla-S$ have the following properties, for a
common exponent $\beta>0$:
\begin{enumerate}
\item $P_*$ has $C^{1+\beta}$ local sections and is equivariant:
\[
 D_x\alpha(\gamma)(P_*)_x=(P_*)_{\alpha(\gamma)x}\rho(\gamma).
\]
\item $\wtd\nabla$ is a H\"older $\Gamma$-invariant affine
connection. Its parallel transport preserves $P_*$ and its associated
positive definite $C^{1+\beta}$ metric $g$; equivalently,
$\wtd\nabla g=0$.
\item Its torsion $T^{\wtd\nabla}$ is $C^{1+\beta}$,
has constant $G$- and $C_*$-equivariant adapted coefficients $T_0$,
and satisfies $\wtd\nabla T^{\wtd\nabla}=0$.
\item Its distributional curvature is zero.
\item The alternating map $B:V\times V\to V$ defined by
$B(u,v)=-T_0(u,v)$ is a nilpotent Lie bracket.
\end{enumerate}
The regularity of the connection itself is only asserted to be
H\"older. No integration to a global nilmanifold is asserted in
this proposition.
\end{theorem}

\begin{proof}
\Cref{prop:reduction1} gives the regularity and equivariance of $P_*$.
By \cref{prop:newconnred}, the connection $\tilde{\nabla}$ is H\"older and $\Gamma$-invariant, its parallel transport preserves $P_*$, and its torsion is parallel.
\cref{sec:adapted-metric} gives the associated $C^{1+\beta}$ metric $g$ and shows that $\tilde{\nabla} g=0$.

In every adapted frame, the torsion has the constant coefficient
\[
  T_0(u,v)=A_0(v)u-A_0(u)v.
\]
It is $G$ and $C_*$-equivariant by \cref{lem:unique-A0-lift}.
Since the adapted frames are $C^{1+\beta}$, so is the torsion.

By \cref{lem:regular-curvature}, the curvature distribution is
represented in adapted frames by the constant coefficient $R_0$.
Finally, \cref{cor:curvature0} gives $R_0=0$ and shows that
$B=-T_0$ is a nilpotent Lie bracket on $V$.  Hence the curvature
vanishes as a distribution.
\end{proof}

\section{Nilpotent group, deck group, and conjugacy}\label{sec7}

We retain the setting of \cref{setting}.  In this section, we show that $M$ is diffeomorphic to an infra-nilmanifold and the action is conjugated to an affine action on it.

Fix $V=\Rbb^{\dim M}$. By \cref{sec5}, there are a continuous representation $\rho:G\to\GL(V)$, a compact subgroup $C_*\subset\GL(V)$ commuting with $\rho(G)$, and a
$C^{1+\beta}$ principal $C_*$-reduction $P_*\subset\Fr(TM)$ such that, for all $\gamma\in\Gamma$ and $x\in M$,
\begin{equation}\label{eq61} D_x\alpha(\gamma)(P_*)_x  =(P_*)_{\alpha(\gamma)x}\rho(\gamma).\end{equation}

Fix $x_0\in M$ and a frame $b_0:V\to T_{x_0}M$ belonging to $(P_*)_{x_0}$. We use the following properties of $\wtd\nabla$ established in \cref{sec6}.

\begin{enumerate}
\item For some $0<\beta<1$, the coefficients of $\wtd\nabla$ in smooth local frames are $C^\beta$, and the curvature of $\wtd\nabla$ vanishes in the sense of distributions.

\item Parallel transport of $\wtd\nabla$ preserves the positive definite $C^{1+\beta}$ Riemannian metric $g$ on $M$.

\item Parallel transport of $\wtd\nabla$ preserves its torsion $T^{\wtd\nabla}$. The alternating map
\[ B(u,v)\coloneqq -b_0^{-1}\left(T^{\wtd\nabla}_{x_0}(b_0u,b_0v)\right), \qquad u,v\in V,\] is a nilpotent Lie bracket on $V$.

\item For every $\gamma\in\Gamma$, the diffeomorphism $\alpha(\gamma)$ preserves $\wtd\nabla$.

\item Parallel transport of $\wtd\nabla$ preserves the compact reduction $P_*$. 
\end{enumerate}

By (1), the lifted connection on $\wtd M$ admits a global parallel frame. By (3), the Lie brackets of its parallel vector fields are given by $B$. By (2) and compactness of $M$, these fields are complete. We therefore identify $\wtd M$ with the simply connected nilpotent Lie group $N$ with Lie algebra $(V,B)$.

In this identification, the deck transformations act affinely. By (5), their linear parts belong to $C_*$. Since the deck action is free, properly discontinuous, and cocompact, the generalized first Bieberbach theorem gives finite linear holonomy. Thus $M$ is an infranilmanifold.

By (4), every lift of $\alpha(\gamma)$ also acts affinely on $N$ and normalizes the deck group. Hence the same identification conjugates $\alpha$ to an affine action on $M$ in its infranilmanifold structure. Finally, (5) and \eqref{eq61} imply that the linear part of each such lift belongs to $C_*\rho(\gamma)$.

\subsection{Nilpotent lie algebra structure on $T\wtd{M}$}

\subsubsection{Local parallel frames} Let $\Ocal$ be a sufficiently small, simply connected coordinate ball. Let $q:\Ocal\to\Omega\subset\Rbb^d$ be a smooth coordinate chart, $q(x)=\left(u^1(x),\ldots,u^d(x)\right)$, whose image is a ball.
 In the coordinate frame $e_x=(dq_x)^{-1}:\Rbb^d\to T_xM$, write \[ \wtd\nabla=d+A,\qquad A=\sum_{i=1}^d A_i\,du^i.\] Note that $A_i$ are locally $\beta$-H\"older, for some $0<\beta<1$, and that $dA+A\wedge A=0$ in distributions. 

\begin{lemma}[Parallel frame]\label{pframe} 
$\Ocal$ admits a locally $C^{1+\beta}$ parallel frame $\bar b$, $d\bar b=-A \bar b$. More precisely, for any two points $x,y\in\Ocal$ and any piecewise $C^1$ path $c$ in $\Ocal$ from $x$ to $y$, parallel transport $\Tcal_c=\bar b_y \circ \bar b_x^{-1}: T_xM\to T_yM$ only depends on the end points.
\end{lemma}
We will use the following theorem: 
\begin{theorem}[{\cite[Theorem 6.4, Lemma 6.1]{Mardare}}]
\label{thm:mardare}
Let $\Omega\subset\Rbb^d$ be a connected, simply connected open set, where $d\ge2$. Fix an integer $m\ge1$ and $d<p<\infty$. Suppose that
\[C_i\in L^p_{\textrm{loc}}     \left(\Omega,\textrm{Mat}_{m\times m}(\Rbb)\right),\qquad 1\le i\le d,\] satisfy \[\partial_jC_i-\partial_iC_j=C_iC_j-C_jC_i\qquad (1\le i,j\le d)\] in the sense of distributions.

For every $x_0\in\Omega$ and $Y_0\in\GL(m,\Rbb)$, there is a unique \[Y\in W^{1,p}_{\textrm{loc}}  \left(\Omega,\textrm{Mat}_{m\times m}(\Rbb)\right)\] such that \[\partial_iY=YC_i\quad\text{in distributions},\qquad Y(x_0)=Y_0.\] Moreover, $Y(x)$ is invertible for every $x\in\Omega$, and
\[Y^{-1}\in W^{1,p}_{\textrm{loc}}\left(\Omega,\textrm{Mat}_{m\times m}(\Rbb)\right),\qquad\partial_i(Y^{-1})=-C_iY^{-1}\] in distributions.
\end{theorem}

\begin{proof}[Proof of \cref{pframe}]
Fix $x_0\in\Ocal$ and put $u_0=q(x_0)$. We want to find the frame in the form
\[ \bar b_x=e_x\circ B(q(x)),\qquad B:\Omega\longrightarrow\textrm{GL}(d,\Rbb),  \qquad B(u_0)=I. \]
To determine the equation for $B$, fix a constant vector $v\in\Rbb^d$. The vector field $x\mapsto\bar b_xv$ has coordinate column $B(u)v$, where $u=q(x)$. Thus its covariant derivative in the $i$th coordinate direction has coordinate column \[ \partial_i(Bv)+A_iBv=(\partial_iB+A_iB)v.\]
Consequently all columns of the frame are parallel precisely when $\partial_iB=-A_iB$ for every $i$, i.e., 
$d\bar b=-A\bar b$ in the chosen coordinate frame.

Assume first that $d\ge2$. Choose a finite $p>d/(1-\beta)$. Each $A_i$ is continuous and hence belongs to $L^p_{\textrm{loc}}$. To put the equation into the right-multiplication convention of \cref{thm:mardare}, set $C_i=-A_i^{\tr}$. Distributional flatness means, in coordinates, that
\[ \partial_iA_j-\partial_jA_i=A_jA_i-A_iA_j.\]
Since transposition commutes with distributional differentiation and reverses the order of matrix multiplication, we obtain
\begin{align*}  \partial_jC_i-\partial_iC_j  &= (\partial_iA_j-\partial_jA_i)^{\tr}\\
 &= (A_jA_i-A_iA_j)^{\tr}\\
 &= A_i^{\tr}A_j^{\tr} -A_j^{\tr}A_i^{\tr} = C_iC_j-C_jC_i.
\end{align*}
Hence, we can apply \cref{thm:mardare}. The domain $\Omega$ is connected and simply connected, and $p>d$. Applying \cref{thm:mardare} with matrix size $m=d$ and initial value $I$ therefore gives $Y\in W^{1,p}_{\textrm{loc}}(\Omega)$ satisfying
\[ \partial_iY=YC_i,\qquad Y(u_0)=I. \]
Its continuous representative is invertible at every point.   Put $B=Y^{\tr}$. Then $B(u_0)=I$ and, initially in distributions,
\[ \partial_iB=(YC_i)^{\tr}=C_i^{\tr}B=-A_iB.\]

It remains to prove that this weak solution gives a classical parallel frame. Morrey's embedding gives $B\in C^\alpha_{\textrm{loc}}$ with $\alpha=1-d/p>\beta$. On each smaller ball the factors $A_i$ and $B$ are bounded, so the identity
\[ (A_iB)(u)-(A_iB)(w)=(A_i(u)-A_i(w))B(u)+A_i(w)(B(u)-B(w)) \]
shows that $A_iB$ is $\beta$-H\"older there. Since $B$ is continuous and its distributional partial derivatives $\partial_iB=-A_iB$ are locally $\beta$-H\"older, $B\in C^1_{\textrm{loc}}$ with these same classical derivatives. Hence $B\in C^{1+\beta}_{\textrm{loc}}$, and $\bar b_x=e_x\circ B(q(x))$ is a $C^{1+\beta}_{\textrm{loc}}$ parallel frame.
\end{proof}

We may choose the parallel frame $\bar b$ to be a local section of
$P_*$. Indeed, let $e$ be the parallel frame given by \cref{pframe}
on a coordinate ball $U$. Fix $y\in U$ and choose
$p\in(P_*)_y$. Set
\[
 \bar b_x=e_x(e_y^{-1}p),\qquad x\in U.
\]
Since $e_y^{-1}p$ is a constant in $x$, $\bar b$ is parallel and
$C^{1+\beta}$, and $\bar b_y=p$. For any path $c$ in $U$ from
$y$ to $x$, parallelism gives
\[
 \bar b_x=\mathcal T_c p\in(P_*)_x,
\]
because parallel transport preserves $P_*$. Thus $\bar b$ is the
unique parallel frame on $U$ with initial value $p$.

\subsubsection{Parallel fields on the universal cover}
We lift the connection and the metric we build to the universal cover in order to define the nilpotent Lie algebra structure on the tangent space. 
\begin{definition}
Let $\pi:\wtd M\to M$ be the universal covering map. Write \[ \widehat\nabla=\pi^*\wtd\nabla,\qquad \widehat g=\pi^*g,\] and define the lifted reduction by \[ (\widehat P_*)_x =\{b:V\to T_x\wtd M:\ d_x\pi\circ b\in(P_*)_{\pi(x)}\}. \] Fix $x_0\in\wtd M$ and $b_0\in(\widehat P_*)_{x_0}$. Let $\langle\cdot,\cdot\rangle_*$ (and $|\cdot|_*$) be the fixed $C_*$-invariant inner product (and norm, respectively) on $V$ defining the associated metric $g$. 
\end{definition}

\begin{proposition}\label{globfr}
There is a unique global $C^{1+\beta}$ frame \[ b_x:V\longrightarrow T_x\wtd M \] which is $\widehat\nabla$-parallel and satisfies $b_{x_0}=b_0$. It belongs to $\widehat P_*$ at every point and satisfies \[ \widehat g_x(b_xu,b_xv)=\langle u,v\rangle_*.\]
For $u\in V$, define the vector field $X_u$ on $\wtd M$ by $X_u(x)=b_xu$. Then
\begin{equation} \label{sec7brac}  [X_u,X_v]=X_{B(u,v)},\qquad B=-T_0.\end{equation}
\end{proposition}

\begin{proof}
Let $\pi:\wtd M\to M$ be the universal covering map, and write $\widehat\nabla=\pi^*\wtd\nabla$ and $\widehat g=\pi^*g$. Fix $x_0\in\wtd M$ and $b_0\in(\widehat P_*)_{x_0}$ as in the statement. Since $g$ is the metric associated with $P_*$ and $\langle\cdot,\cdot\rangle_*$, we have
\[ \widehat g_{x_0}(b_0u,b_0v)=\langle u,v\rangle_* \qquad(u,v\in V).\]

On each lift of a coordinate ball in \cref{pframe}, the lifted connection admits a $C^{1+\beta}$ parallel frame $\bar b$. Consequently, transport along a path in that ball from $y$ to $z$ is $\bar b_z\bar b_y^{-1}$ and depends only on $y$ and $z$. It follows that parallel transport is unchanged under a piecewise $C^1$ homotopy fixing the endpoints. Since $\wtd M$ is simply connected, we can define
\[ b_x:=\mathcal T_c\circ b_0:V\longrightarrow T_x\wtd M,\qquad c\text{ any piecewise $C^1$ path from $x_0$ to $x$}.\]

To check regularity and parallelism, fix a parallel-frame
neighborhood $U$, a point $y\in U$, and a local parallel frame
$\bar b$ on $U$. For $x\in U$, first take a path from $x_0$ to $y$
and then a path inside $U$ from $y$ to $x$. The local transport
formula gives
\[
 b_x=\bar b_x\bigl(\bar b_y^{-1}b_y\bigr).
\]
The factor in parentheses is a single constant matrix. Thus $b$
is $C^{1+\beta}$ and parallel on $U$, and hence on $\wtd M$.
Any parallel frame with value $b_0$ at $x_0$ must have value
$\mathcal T_c\circ b_0$ at $x$, which proves uniqueness.

Parallel transport preserves $\what P_*$ and $\what g$, so $b_x\in(\what P_*)_x$ and 
$\what g_x(b_xu,b_xv)=\langle u,v\rangle_*$ for all $x,u,v$. By the torsion formula in \cref{sec6conclusion},
\[ T^{\what\nabla}_x(b_xu,b_xv)=b_xT_0(u,v).\] Since the fields $X_u,X_v$ are parallel,  \[T^{\what\nabla}(X_u,X_v)=\what\nabla_{X_u}X_v-\what\nabla_{X_v}X_u-[X_u,X_v]=-[X_u,X_v].\] Consequently $[X_u,X_v]=-bT_0(u,v)=X_{B(u,v)}$. 
\end{proof}

Fix the global $C^{1+\beta}$ frame $b$ and vector fields $X_u$ for $u\in V$ as in \cref{globfr}.

\begin{lemma}\label{completi}
Let $T>0$, let $u:[0,T]\to V$ be continuous, and fix $x_0\in\wtd M$. There is a unique $C^1$ curve $x:[0,T]\to\wtd M$ satisfying 
\begin{equation}\label{eq721} \dot x(t)=b_{x(t)}u(t),\qquad x(0)=x_0.\end{equation}
In particular, for each fixed $v\in V$, the vector field $X_v(x)=b_xv$ on $\wtd M$ is complete.
\end{lemma}

\begin{proof}
Extend continuously $u$ to $\Rbb$ by setting $u(t)=u(0)$ for $t<0$ and $u(t)=u(T)$ for $t>T$. Write $K=\max_{0\leq t\leq T}|u(t)|_*$. Since the frame $b$ is $C^1$, the vector $b_xu(t)$ depends continuously on $t$ and locally Lipschitz continuously on $x$, uniformly for $t$ in bounded intervals. \cref{ode:existence}  gives a unique solution through $x_0$ on the maximal interval.

Choose a smooth Riemannian metric $g_0$ on $M$. Since $M$ is compact, there is a constant $C\geq1$ such that $|w|_{g_0}\leq C|w|_g$ for every $w\in TM$. Hence, $|\wtd w|_{\what g_0} \le C|\wtd w|_{\what g}$ for every $\wtd w\in T\wtd M$ where $\what g_0=\pi^*g_0$. Note that $\what g_0$ is complete because it is the lift of a complete Riemannian metric. Since $b_x:(V,|\cdot|_*)\to(T_x\wtd M,|\cdot|_{\widehat g})$ is an isometry, \eqref{eq721} gives \[ |\dot x(t)|_{\widehat g_0} \le C|\dot x(t)|_{\widehat g} =C|u(t)|_* \le CK.\] Thus, for any two times $s<t$ in the interval of existence, \[ d_{\what g_0}\left(x(s),x(t)\right)\le CK(t-s). \] Hence,  the solution is Cauchy in a complete metric.

If the maximal interval had a finite endpoint $t_*$, this estimate would make $x(t)$ Cauchy as $t$ approaches $t_*$. Completeness gives a limit $x_*\in\wtd M$. The local existence theorem at $(t_*,x_*)$ then continues the solution past $t_*$, contradicting maximality. Hence the solution exists for every real time. Restricting it to $[0,T]$ proves the first assertion. Taking $u(t)\equiv v$ proves that every integral curve of $X_v$ exists for all $t\in\Rbb$, which is precisely completeness of $X_v$.
\end{proof}

\subsection{A developing map and the nilpotent Lie group structure on $\wtd M$}

Let $N$ be the connected simply connected Lie group with Lie algebra $(V,B)$, where $B=-T_0$, and let $\id_N$ denote its identity. Since $(V,B)$ is nilpotent, $N$ is nilpotent. We identify $T_{\id_N}N=\textrm{Lie}(N)$ with $V$. For $a\in N$, let \[L_a:N\to N\] be left
translation by $a$. 

Define the $V$-valued one-forms
\[ \theta_x=b_x^{-1}:T_x\wtd M\longrightarrow V,\qquad (\theta_N)_a=(dL_{a^{-1}})_a:T_aN\longrightarrow V. \]

The first is the coframe dual to the parallel frame of \cref{globfr}. The second is the left Maurer--Cartan form of $N$. We will identify $\wtd M$ with $N$ by finding the map that identifies these two forms and sends $x_0$ to $\id_N$.

We follow the notation in \cite[Chapter~2, Lemma~2.7]{McKayCartan}. A \emph{Maurer--Cartan system} on $M$ is a collection of $C^1$ one-forms $\theta^1,\ldots,\theta^d$ with the following properties:
\begin{enumerate}
\item For every $x\in Q$, the covectors $\theta^1_x,\ldots,\theta^d_x$ form a basis of $T_x^*M$.
\item There are real constants $c^k_{ij}$, independent of $x$, such that 
\[ d\theta^k =-\sum_{1\le i<j\le d}c^k_{ij}\,\theta^i\wedge\theta^j, \qquad 1\le k\le d.\]
\end{enumerate}
A vector field $X$ on $M$ is called \emph{constant with respect to the system} if each function $\theta^i(X)$ is constant on $Q$. The system is \emph{complete} if every such vector field has a flow defined for all $t\in\mathbb R$.

\begin{lemma} \label{MCgeomeric}
Fix a basis $e_1,\ldots,e_d$ of $V$. Write
\[ B(e_i,e_j)=\sum_{k=1}^d c^k_{ij}e_k, \qquad \theta=\sum_{k=1}^d\theta^k e_k.\]
Then $\theta^1,\ldots,\theta^d$ form a $C^{1+\beta}$ coframe on $\wtd M$ satisfying, for  $1\le k\le d$,
 \[ d\theta^k=-\sum_{i<j}c^k_{ij}\,\theta^i\wedge\theta^j,\]
Especially, they form a Maurer--Cartan system with the structure constants of $(V,B)$. This system is complete, i.e., for every $u\in V$, the unique vector field $X$ on $\wtd M$ satisfying $\theta(X)=u$ is the complete parallel field $X_u$.
\end{lemma}

\begin{proof}
Put $X_i=X_{e_i}$. Since $X_i(x)=b_xe_i$, we have $\theta^k(X_i)=\delta^k_i$. The formula for the exterior derivative of a one-form and \eqref{sec7brac} give
\[ \begin{aligned} d\theta^k(X_i,X_j)&=X_i\left(\theta^k(X_j)\right) -X_j\left(\theta^k(X_i)\right)-\theta^k([X_i,X_j])\\
 &=-\theta^k\left(X_{B(e_i,e_j)}\right) =-c^k_{ij}.
 \end{aligned}\]
As the $X_i$ form a frame, these values determine $d\theta^k$.
Equivalently, in $V$-valued notation, \[ d\theta(Y,Z)=-B\left(\theta(Y),\theta(Z)\right)\] for tangent vectors $Y,Z$ at the same point. The identical calculation with left-invariant fields gives this equation for $\theta_N$ on $N$.

Finally, $\theta(X)=u$ if and only if $X(x)=b_xu$ at every $x$. This is exactly the field $X_u$, which is complete by \cref{completi}. This verifies completeness in the sense of \cite[Chapter~2, Theorem~2.9]{McKayCartan}. \end{proof}

\begin{proposition}[Developing map]\label{dev.map}
There is a unique $C^1$ diffeomorphism $\Phi:\wtd M\to N$ such that $\Phi(x_0)=\id_N$ and
\begin{equation} \label{devmapdif}
 \Phi^*\theta_N=\theta, \qquad\text{equivalently}\qquad  d_x\Phi\circ b_x=(dL_{\Phi(x)})_{\id_N}. \end{equation}
\end{proposition}
\begin{proof}
By \cref{MCgeomeric}, $\theta$ is a complete Maurer--Cartan system on $\wtd M$ with Lie algebra $(V,B)$. Since $\wtd M$ is simply connected, \cite[Theorem~2.8]{McKayCartan} gives a unique developing local diffeomorphism $\Phi:\wtd M\to N$ satisfying $\Phi(x_0)=\id_N$ and $\Phi^*\theta_N=\theta$. By completeness, \cite[Theorem~2.9]{McKayCartan} makes $\Phi$ a covering of $N$. Since $N$ is simply connected and $\wtd M$ is connected, $\Phi$ is a diffeomorphism. The arguments apply to our setting as $\theta$ is $C^{1+\beta}$. Indeed, we only need $C^1$ regularity. 
\end{proof}

\begin{corollary}\label{dev.map2}
Both $\Phi$ and $\Phi^{-1}$ are $C^{2+\beta}$.
\end{corollary}
\begin{proof}
The identity $\Phi^*\theta_N=\theta$ gives \begin{equation}\label{eq:781}d_x\Phi=(dL_{\Phi(x)})_{\id_N}\circ\theta_x.\end{equation}
Since $\Phi$ is $C^1$, $\theta$ is $C^{1+\beta}$, and $a\mapsto(dL_a)_{\id_N}$ is smooth, the right-hand side is first $C^1$, giving $\Phi\in C^2$. 
Using \eqref{eq:781} again with the fact that $\Phi\in C^2$, the right-hand side of \eqref{eq:781} is $C^{1+\beta}$, so $\Phi\in C^{2+\beta}$. The inverse function theorem gives the same regularity for $\Phi^{-1}$.
\end{proof}

\begin{corollary} \label{Nilpgp}
The operations
\[ x\star y=\Phi^{-1}\left(\Phi(x)\Phi(y)\right),\qquad x^{-1}_{\star}=\Phi^{-1}\left(\Phi(x)^{-1}\right) \]
give $\wtd M$ a group structure with identity $x_0$. With respect to the smooth structure transported from $N$ (via $\Phi$), this is a connected simply connected nilpotent Lie group with Lie algebra $(V,B)$. In the original smooth structure on $\wtd M$, multiplication and inversion are $C^{2+\beta}$. For each $u\in V$, the field $X_u$ corresponds under $\Phi$ to the left-invariant field $u^L$ on $N$.
\end{corollary}

\begin{proof}
The formulas transport the group operations of $N$ by the diffeomorphism $\Phi$. The regularity follows from
\cref{dev.map}. Finally, for every $x\in\wtd M$,
\[ d_x\Phi(X_u(x))=d_x\Phi(b_xu)=(dL_{\Phi(x)})_{\id_N}u=u^L(\Phi(x)),\] as required.
\end{proof}

\subsection{The deck group and the infranilmanifold structure on $M$}
Retain the notations and settings in the previous subsections. The following lemma says that the deck transformation is given by affine since it should preserve $\what\nabla$.
\begin{lemma}\label{connpresaff}
Let $F:\wtd M\to\wtd M$ be a $C^2$ diffeomorphism preserving $\what\nabla$, meaning that
\begin{equation}\label{sec7conninv} \widehat\nabla_{F_*Y}(F_*Z)=F_*(\widehat\nabla_YZ) \end{equation}
for $C^1$ vector fields $Y,Z$ on open subsets of $\wtd M$. Then there is a constant $A_F\in\textrm{Aut}(V,B)$ such that
\[ d_xF\circ b_x=b_{F(x)}A_F.\] If $\sigma_{A_F}\in\textrm{Aut}(N)$ is the automorphism with derivative $A_F$ at $\id_N$, then
\[ \Phi\circ F\circ\Phi^{-1}=L_{a_F}\circ\sigma_{A_F}, \qquad a_F=\Phi(F(x_0)).\]
\end{lemma}

\begin{proof}
Equation \eqref{sec7conninv} implies that $F_*X_u$ is parallel.

For $y\in\wtd M$, put $x=F^{-1}(y)$ and define \[ A(y)\coloneqq b_y^{-1}\circ d_xF\circ b_x\in\GL(V).\]
For each fixed $u\in V$, \[ (F_*X_u)(y)=d_xF(X_u(x))=d_xF(b_xu)=b_yA(y)u.\] 
The map $A$ is $C^1$ because $F$ is $C^2$ and $b$ is $C^{1+\beta}$. Let $Z$ be a $C^1$ vector field on $\wtd M$. Applying
\eqref{sec7conninv} to $(F^{-1})_*Z$ and $X_u$ gives
\[\what\nabla_Z(F_*X_u)=F_*\left(\widehat\nabla_{(F^{-1})_*Z}X_u\right)=0,\]
since $X_u$ is parallel. We now compute the same derivative in the frame $b$. Fix a basis
$e_1,\ldots,e_d$ of $V$, and put $E_i(y)=b_ye_i$. Write
\[ u=\sum_j u^j e_j,\qquad A(y)e_j=\sum_i A^i{}_j(y)e_i.\] Then $F_*X_u=\sum_{i,j}A^i{}_j u^j E_i$. The numbers $u^j$ are constant, and each vector field $E_i$ is parallel. Hence the product rule gives, at $y$,
\[\begin{aligned}
(\widehat\nabla_Z(F_*X_u))_y &=\sum_{i,j} Z(A^i{}_j)(y)\,u^j E_i(y)+\sum_{i,j} A^i{}_j(y)\,u^j(\widehat\nabla_Z E_i)_y\\
 &=\sum_{i,j} d(A^i{}_j)_y(Z_y)\,u^j\,b_y e_i\\
 &=b_y\!\left(dA_y(Z_y)u\right).
 \end{aligned}\]
Combining the two calculations yields $b_y(dA_y(Z_y)u)=0$. Since $b_y:V\to T_y\wtd M$ is an isomorphism, we obtain $dA_y(Z_y)u=0$. This holds for every $u\in V$ and every $Z_y\in T_y\wtd M$, so $dA_y=0$ for all $y$. As $\wtd M$ is connected, $A$ is constant. Let us denote the constant by $A_F$. We also have \[ X_{A_FB(u,v)}=F_*[X_u,X_v]=[F_*X_u,F_*X_v]=X_{B(A_Fu,A_Fv)}.\]
Therefore $A_FB(u,v)=B(A_Fu,A_Fv)$. The corresponding group automorphism $\sigma_{A_F}$ is
\[ \sigma_{A_F}(\exp_Nu)=\exp_N(A_Fu). \]
It is a group automorphism because $A_F$ preserves the bracket and hence preserves group operation via Baker--Campbell--Hausdorff formula.

Put $\widehat F=\Phi\circ F\circ\Phi^{-1}$. By \eqref{devmapdif}, it sends the left-invariant field $u^L$ to $(A_Fu)^L$. Uniqueness of integral curves gives
\[ \what F(a\exp_N(tu))=\widehat F(a)\exp_N(tA_Fu). \] Taking $a=\id_N$ and $t=1$, we proved 
\[\what F=L_{a_F}\circ\sigma_{A_F}\] using surjectivity of $\exp_N$.
\end{proof}

\subsubsection{The deck group and the infranilmanifold structure}
Now, we are ready to show that $M$ is an infra-nilmanifold. This shows the infra-nilmanifold structure part of \cref{thm:main}. 

We use the following form of the generalized first Bieberbach theorem, for instance, see   \cite[Proposition~3.2 and Theorem~3.4]{DekimpeGuide}.
\begin{theorem}[Generalized first Bieberbach theorem]\label{Bieber}
Let $N$ be a connected, simply connected nilpotent Lie group, let $C$ be a maximal compact subgroup of $\Aut(N)$, and let $\Delta$ be a discrete cocompact subgroup of $N\rtimes C$. Identify $N$ with its subgroup of left translations. Then $\Delta\cap N$ is a cocompact lattice in $N$, and the linear-part homomorphism $\ell(a,A)=A$ induces an isomorphism \[ \Delta/(\Delta\cap N)\simeq\ell(\Delta)\]
onto a finite group. The action of $\Delta$ on $N$ is free if and only if $\Delta$ is torsion-free. In that case $\Delta\backslash N$ is an infranilmanifold.
\end{theorem}

\begin{proposition}\label{deckgp}
Let $\Dcal=\operatorname{Deck}(\pi)$ be the group of deck transformations for the universal covering $\pi:\wtd M\to M$. Using the developing diffeomorphism $\Phi:\wtd M\to N$, normalized by $\Phi(x_0)=\id_N$, set \[ \Pi=\Phi\Dcal\Phi^{-1},\qquad K=C_*\cap\Aut(V,B).\]
Identify each $A\in K$ with the automorphism $\sigma_A$ of $N$ whose derivative at $\id_N$ is $A$. Then $\Pi$ is a discrete torsion-free subgroup of $N\rtimes K$ acting freely, properly discontinuously, and cocompactly on $N$. The map 
\[ \varphi:M\longrightarrow X\coloneqq\Pi\backslash N,\qquad \varphi(\pi(x))=\Pi\cdot\Phi(x),\] is a $C^{2+\beta}$ diffeomorphism onto an infranilmanifold.

Moreover, $\Lambda:=\Pi\cap N$ is a cocompact lattice in $N$, and $\Lambda\backslash N\to X$ is a finite regular covering with deck group $\Pi/\Lambda$.
\end{proposition}
\begin{proof}
Every $\delta\in\Dcal$ preserves $\widehat\nabla$ and $\widehat P_*$ because both are pulled back from $M$. By \cref{connpresaff},  
$\Phi\circ\delta\circ\Phi^{-1}$ is affine with linear part
\[ A_\delta=b_{\delta(x)}^{-1}\circ d_x\delta\circ b_x\in\Aut(V,B).\] Since $d_x\delta\circ b_x$ and $b_{\delta(x)}$ are in the same fiber $(\what P_*)_{\delta(x)}$,  $A_\delta$ lies in $C_*$. Thus $\Pi\subset N\rtimes K$, where $K=C_*\cap\Aut(V,B)$ is compact because $\Aut(V,B)$ is closed in $\GL(V)$.

Conjugation by $\Phi$ identifies the deck action on $\wtd M$ with the action of $\Pi$ on $N$. Hence $\Pi$ acts freely and properly discontinuously, and $\Phi$ descends in covering charts to the $C^1$ diffeomorphism
\[ \varphi:M\longrightarrow\Pi\backslash N,\qquad \varphi(\pi(x))=\Pi\cdot\Phi(x). \]
In particular, $\Pi\backslash N$ is compact.

Choose a maximal compact subgroup $C\subset\Aut(N)$ containing $K$. Proper discontinuity implies that $\Pi$ is discrete in $N\rtimes C$. Since $\Pi\backslash N$ and $C$ are compact, $\Pi$ is also cocompact in $N\rtimes C$.

We now apply \cref{Bieber}. Since $\Pi$ acts freely on $N$, $\Pi$ is torsion-free and $X=\Pi\backslash N$ is an infranilmanifold. Furthermore, $\Lambda=\Pi\cap N$ is a cocompact lattice in $N$ of finite index in $\Pi$. Finally, $\Lambda=\ker(\ell|_\Pi)$ is normal in $\Pi$, so 
$\Lambda\backslash N\to X$ is a finite regular covering with deck group $\Pi/\Lambda$.
\end{proof}

\subsubsection{Conjugacy}
Finally, we show that the action $\alpha$ is conjugate to an affine action through the map $\varphi$ we constructed. This shows the conjugacy part of  \cref{thm:main}.
\begin{proposition}\label{affconj}
The $C^{2+\beta}$-diffeomorphism $\varphi:M\to X$ we constructed in \cref{deckgp} conjugates the action of the
entire lattice $\Gamma$ to an affine action on $X$.
\end{proposition}

\begin{proof}
For each $\gamma\in\Gamma$, choose a lift $F_\gamma:\wtd M\to\wtd M$ of $\alpha(\gamma)$. $F_\gamma$ is a $C^{2+\eta}$ diffeomorphism and $\pi\circ F_\gamma=\alpha(\gamma)\circ\pi$. Since $\alpha(\gamma)$ preserves $\wtd\nabla$, $F_\gamma$ preserves $\widehat\nabla$. By
\cref{connpresaff},
\[ \widehat F_\gamma:=\Phi\circ F_\gamma\circ\Phi^{-1}=L_{a_\gamma}\circ\sigma_{A_\gamma},\qquad a_\gamma=\Phi(F_\gamma(x_0)).\]
Moreover, $d_xF_\gamma\circ b_x=b_{F_\gamma(x)}A_\gamma$. By the lifted equivariance of $P_*$ in \eqref{eq61}, the left side belongs to $(\what P_*)_{F_\gamma(x)}\rho(\gamma)$. It follows that $A_\gamma\in C_*\rho(\gamma)$. Since $C_*$ commutes with $\rho(G)$, we may write $A_\gamma=\rho(\gamma)k_\gamma$ for a constant $k_\gamma\in C_*$ depending on the chosen lift. For each $\delta\in\Dcal$,
\[ \pi\circ F_\gamma\circ\delta\circ F_\gamma^{-1}=\alpha(\gamma)\circ\pi\circ\delta\circ F_\gamma^{-1}=\pi.\]
Applying the same argument to $F_\gamma^{-1}$ gives $F_\gamma\Dcal F_\gamma^{-1}=\Dcal$. Conjugating by $\Phi$ shows that $\widehat F_\gamma$ normalizes $\Pi$. Its affine action on $N$ therefore induces the affine diffeomorphism $\rho_\gamma: \Pi\backslash N\to \Pi\backslash N$
\[\rho_\gamma( \Pi\cdot n) = \Pi\cdot\widehat F_\gamma(n)\] for $n\in N$.
of $X$. Using $\phi:M\to \Pi\backslash N$,   $\rho_\gamma=\varphi\circ\alpha(\gamma)\circ\varphi^{-1}$. This identity gives independence of both the chosen lift and the choice of $\gamma\in\Gamma$.
\end{proof}

\newpage

\bibliography{volanosov_references.bib}

\bibliographystyle{alpha}
\end{document}